\documentclass[a4paper,twoside,11pt]{book}
\usepackage[german,english]{babel}
\usepackage{amsmath, amsxtra, latexsym,amscd, amsthm, amssymb, indentfirst}
\usepackage[mathscr]{eucal}
\usepackage{graphics, graphpap}
\usepackage{graphicx}
\usepackage{makeidx}
\usepackage[all]{xy}
\usepackage[utf8]{inputenc}
\usepackage{pgf,tikz}
\usetikzlibrary{arrows}
\usetikzlibrary{matrix}
\usetikzlibrary{calc}
\usepackage{lipsum}
\usepackage{array, tabularx, longtable}
\usepackage{multicol}

\newtheorem{theorem}{Theorem}[section]

\newtheorem{lemma}[theorem]{Lemma}
\newtheorem{proposition}[theorem]{Proposition}
\newtheorem{corollary}[theorem]{Corollary}
\newtheorem{claim}{Claim}

\newtheorem{problem}{Problem}
\newtheorem{inproblem}{Problem}
\newtheorem{intheorem}{Theorem}

\newtheoremstyle{definition}
  {6pt}
  {6pt}
  {}
  {}
  {\bfseries}
  {.}
  {.5em}
  {}%
\theoremstyle{definition}
\newtheorem{definition}[theorem]{Definition}
\newtheorem{indefinition}{Definition}
\newtheorem{example}[theorem]{Example}
\theoremstyle{remark}
\newtheorem{remark}[theorem]{Remark}
\def\rem{\begin{remark}}
\def\reme{\end{remark}}
\def\ex{\begin{example}}
\def\exe{\end{example}}
\def\tr{\begin{theorem}}
\def\tre{\end{theorem}}
\def\prop{\begin{proposition}}
\def\prope{\end{proposition}}
\def\df{\begin{definition}}
\def\dfe{\end{definition}}
\def\cor{\begin{corollary}}
\def\core{\end{corollary}}
\def\lem{\begin{lemma}}
\def\leme{\end{lemma}}
\def\pr{\begin{proof}}
\def\pre{\end{proof}}
\def\prob{\begin{problem}}
\def\probe{\end{problem}}

\def\it{\begin{itemize}}
\def\ite{\end{itemize}}
\def\intr{\begin{intheorem}}
\def\intre{\end{intheorem}}
\def\inprob{\begin{inproblem}}
\def\inprobe{\end{inproblem}}
\def\indf{\begin{indefinition}}
\def\indfe{\end{indefinition}}
\def\qnn{$^qN$}
\DeclareMathOperator{\hkf}{HKF}

\newcommand{\diff}{\operatorname{diff}}
\newcommand{\qn}{{^qN}}

\newcommand{\qp}{ ^qP}
\newcommand{\nn}{\hat{N}}
\newcommand{\spec}{\operatorname{Spec}}
\newcommand{\rad}{\operatorname{rad}}
\newcommand{\nil}{\operatorname{nil}}
\newcommand{\im}{\operatorname{im}}
\newcommand{\sep}{\operatorname{sep}}
\newcommand{\tf}{\operatorname{tf}}
\newcommand{\red}{\operatorname{red}}
\newcommand{\ann}{\operatorname{Ann}}
\newcommand{\krull}{\operatorname{Krull}}
\newcommand{\supp}{\operatorname{supp}}
\newcommand{\trian}{\operatorname{\bigtriangleup}}
\newcommand{\len}{\operatorname{length}}

\newcommand{\gen}{\operatorname{Gen}}
\newcommand{\rk}{\operatorname{rank}}

\newcommand{\vol}{\operatorname{vol}}
\newcommand{\chara}{\operatorname{char}}

\newcommand*{\suchthat}[1]{\left|\vphantom{#1}\right.}

\def\t{\mathbb T}
\def\r{\mathbb R}

\def\n{\mathbb N}
\def\q{\mathbb Q}
\def\c{\mathbb C}

\def\z{\mathbb Z}

\DeclareRobustCommand\longtwoheadrightarrow {\relbar\joinrel\twoheadrightarrow}
\newcommand*{\longhookrightarrow}{\ensuremath{\lhook\joinrel\relbar\joinrel\rightarrow}}
\DeclareFontFamily{U}  {MnSymbolC}{}
\DeclareSymbolFont{MnSyC}         {U}  {MnSymbolC}{m}{n}
\SetSymbolFont{MnSyC}       {bold}{U}  {MnSymbolC}{b}{n}
\DeclareFontShape{U}{MnSymbolC}{m}{n}{
<-6> MnSymbolC5
<6-7> MnSymbolC6
<7-8> MnSymbolC7
<8-9> MnSymbolC8
<9-10> MnSymbolC9
<10-12> MnSymbolC10
<12-> MnSymbolC12}{}
\DeclareFontShape{U}{MnSymbolC}{b}{n}{
<-6> MnSymbolC-Bold5
<6-7> MnSymbolC-Bold6
<7-8> MnSymbolC-Bold7
<8-9> MnSymbolC-Bold8
<9-10> MnSymbolC-Bold9
<10-12> MnSymbolC-Bold10
<12-> MnSymbolC-Bold12}{}

\DeclareMathSymbol{\cupdot}{\mathbin}{MnSyC}{60}

\DeclareFontFamily{U} {MnSymbolF}{}
\DeclareSymbolFont{mnsymbols} {U} {MnSymbolF}{m}{n}
\DeclareFontShape{U}{MnSymbolF}{m}{n}{
<-6> MnSymbolF5
<6-7> MnSymbolF6
<7-8> MnSymbolF7
<8-9> MnSymbolF8
<9-10> MnSymbolF9
<10-12> MnSymbolF10
<12-> MnSymbolF12}{}
\DeclareFontShape{U}{MnSymbolF}{b}{n}{
<-6> MnSymbolF-Bold5
<6-7> MnSymbolF-Bold6
<7-8> MnSymbolF-Bold7
<8-9> MnSymbolF-Bold8
<9-10> MnSymbolF-Bold9
<10-12> MnSymbolF-Bold10
<12-> MnSymbolF-Bold12}{}

\DeclareMathSymbol{\bigcupdot}{\mathop}{mnsymbols}{34}

\begin{document}
\thispagestyle{empty}

\vspace{2cm}
\begin{center}
\huge{\textbf{Hilbert-Kunz theory for binoids}}\\[3.5cm]
\LARGE {BAYARJARGAL BATSUKH}\\[3.5cm]
\Large {Dissertation\\[0.2cm]
zur Erlangung des Doktorgrades (Dr.~rer.~nat.)} \\[3cm]
\Large {FB Mathematik/Informatik \\Universität Osnabrück}\\[0.5cm]
\Large {December 2014}
\end{center}

\chapter*{Acknowledgements}
\thispagestyle{empty}

First and foremost I wish to thank my advisor, professor Holger Brenner for his guidance, trust, support and all the help. 
Without his assistance and dedicated involvement in every step throughout the process, this paper would have never been accomplished. 
I would like to thank you very much for your friendly, kind, sociable relation and understanding over these past three years.  

I would also like to thank the Mongolian Government and National University of Mongolia for the financial support for these three years and gave me a chance to finish my thesis here. 

I could not have finished this thesis without my family, father Batsukh, mother Chuluunbileg, my wife Enkhnasan and my children, who gave me the power, help and support. To my family, thank you. 

This thesis would also not be possible without the support of my dear comrades and friends Le Van Dinh, Christof Söger, Simone Böttger, Daniel Brinkmann, Sadiq Al Nassir and all the professors and students of University of Osnabrück. 
Special thank to Davide Alberelli, who read my thesis and made corrections and many recommendations. Also I thank Hoang Le Truong, who also gave helpful comments. 

\clearpage

\tableofcontents
\pagenumbering{arabic}

\chapter*{Introduction}
\markboth{\MakeUppercase{Introduction}}{}
\addcontentsline{toc}{chapter}{Introduction}

Let $(R,\mathfrak{m})$ be a commutative Noetherian local ring of dimension $d$ and containing a field $K$ of positive characteristic $p$. For an ideal $I$ and a prime power $q=p^e$ we define the ideal $I^{[q]}=\langle a^q|a\in I\rangle$ which is the ideal generated by the $q$th power of the elements of $I$.
Let $I$ be an $\mathfrak{m}$-primary ideal of $R$ and $M$ a finite $R$-module. Then the $R$-modules $M/I^{[q]}M$ have finite length. The \emph{Hilbert-Kunz function} of $M$ with respect to $I$ is $$HKF(I,M)(q)=\len(M/I^{[q]}M).$$ 
If $M=R,I=\mathfrak{m}$ then we have the classical Hilbert-Kunz function $HKF(\mathfrak{m},R)(q)=HK_R(q)$, introduced by Kunz \cite{Kunz}. He also showed that $R$ is regular if and only if $HK_R(q)=q^d$.
In \cite{Monsky}, P.Monsky proved that there is a real constant $c(M)$ such that $$\len(M/I^{[q]}M)=c(M)q^d+O(q^{d-1}).$$
The \emph{Hilbert-Kunz multiplicity} $e_{HK}(I,M)$ of $M$ with respect to $I$ is $$e_{HK}(I,M):=\lim_{q\rightarrow\infty}\dfrac{\len(M/I^{[q]}M)}{q^d}.$$

There are many questions related to Hilbert-Kunz function and multiplicity.
\inprob
 Is the Hilbert-Kunz multiplicity always a rational number?
\inprobe
\inprob
 Is there any interpretation in characteristic 0?
\inprobe
For the following problems the ring comes from a finitely generated $\z$-algebra by reduction modulo $p$. 
\inprob
When and how does the Hilbert-Kunz multiplicity depend on characteristic $p$?
\inprobe
\inprob
Does the limit $$\lim_{p\rightarrow\infty} e_{HK}(I_p,R_p)$$ exist?
\inprobe
\inprob[C.Miller]
 Does the limit $$\lim_{p\rightarrow \infty} \dfrac{\len(R_p/I_p^{[p]})}{p^d}$$ exist?
\inprobe
In most known cases the Hilbert-Kunz multiplicity is a rational number, for example for monoid rings (\cite{Eto}), toric rings (\cite{Watanabe},\cite{Bruns}), monomial ideals and binomial hypersurfaces (\cite{Conca}), rings of finite Cohen-Macaulay type (\cite{Seibert}), two-dimensional graded rings (\cite{Bre06},\cite{Trivedi2}).
Also let $G$ be a finite group acting linearly on a polynomial ring $B=K[x_1,\ldots ,x_n]$ with invariant ring $A=B^G$ and let $I = A_+ B$ be the Hilbert ideal in $B$. Let $R$ be the localization of $A$ at the irrelevant ideal. If $K$ has positive characteristic, then
$$e_{HK}(R) = \frac{ \operatorname{dim}_K \, (B/I ) }{\#( G )}$$ and this is a rational number and depends only on the invariant ring, not on its representation.
But H. Brenner proved in \cite{Brenner} that there exist three-dimensional quartic hypersurface domains and Artinian modules with irrational Hilbert-Kunz multiplicity.

There are many situations where the Hilbert-Kunz multiplicity is independent of the characteristic $p$. For example, for toric rings(\cite{Watanabe}) or invariant rings as above the Hilbert-Kunz multiplicity is independent of the characteristic of the base field at least for almost all prime characteristics. But there are also examples where the Hilbert-Kunz multiplicity is depending on the characteristic (See Section 2.1).

We can ask when the limit of the Hilbert-Kunz multiplicity exists for $p\rightarrow\infty$. If it exists then this limit is a candidate for the Hilbert-Kunz multiplicity in characteristic zero. This  leads us to the question of whether a characteristic zero Hilbert-Kunz multiplicity could be defined directly. In all known cases this limit exists.

H. Brenner, J. Li and C. Miller (\cite{BreLiMil}) have observed that in all known cases where $$\lim_{p\rightarrow\infty}e_{HK}(R_p)$$ exists then this double limit can be replaced by the limit of Problem 5.

If the rings are of a more combinatorial nature, like for example monoid rings (\cite{Eto}, \cite{Watanabe}, \cite{Bruns}), Stanley-Reisner rings and binomial hypersurfaces (\cite{Conca}), we have positive answers for all these problems. 
Also, the proofs in these cases are easier compared to the methods of P. Monsky, C. Han, P. Teixeira (\cite{MonHan}, \cite{MonTei04}, \cite{MonTei06}) or the geometric methods of H. Brenner, V. Trivedi (\cite{Bre06}, \cite{Bre07}, \cite{TriFak03}, \cite{Trivedi1}, \cite{Trivedi2}, \cite{Trivedi07}). 
So we want to generalize these results to a broad and unified concept of a combinatorial ring. For that we work with a new combinatorial structure namely binoids which were introduced in the thesis of S. Boettger \cite{Simone}.
A binoid $(N,+,0,\infty)$ is a monoid with an absorbing element $\infty$ which means that for every $a\in N$ we have $a+\infty=\infty+a=\infty$. In the first chapter we will give the basic properties of binoids and related objects, namely $N$-sets and its partitions, homomorphisms, smash products, exact sequences. Also we will give the definition and some properties on the dimension of binoids and their binoid algebras.

In the second chapter we will define the Hilbert-Kunz function and multiplicity of binoids. This function is given by counting the elements in certain residue class binoids, and not the K-dimension of residue class rings. 
We define the Hilbert-Kunz function and Hilbert-Kunz multiplicity not only for binoids but also for $N_+$-primary ideals of $N$ and a finitely generated $N$-set in the following way:

Let $N$ be a finitely generated, semipositive binoid, $T$ a finitely generated $N$-set and $\mathfrak{n}$ an $N_+$-primary ideal of $N$. Then we call the number $$\hkf^N(\mathfrak{n},T,q)=\hkf(\mathfrak{n},T,q):=\# T/([q]\mathfrak{n}+T)$$ (where for a finite binoid we do not count $\infty$) the Hilbert-Kunz function of $\mathfrak{n}$ on the $N$-set $T$ at $q$. In particular, for $T=N$ and $\mathfrak{n}=N_+$ we have $$\hkf(N,q):=\hkf(N_+,N,q)=\#N/[q]N_+.$$
Note that this function is defined for all $q\in\n$, not only for powers of a fixed prime number. Also note that this residue construction is possible in the category of binoids, not in the category of monoids. The combinatorial Hilbert-Kunz multiplicity is the limit of this function divided by $q^{\dim N}$, provided this limit exists and provided that there is a reasonable notion of dimension.

It turns out that the Hilbert-Kunz function of a binoid for $q=p^e$ is the same as the Hilbert-Kunz function of its binoid algebra over a field of characteristic $p$. 
Hence from here we have a chance to study the above mentioned five problems, by just studying the binoid case. 

Since the Hilbert-Kunz function of binoids is given by just counting elements, not vector space dimensions, this reveals clearer the combinatorial nature of the problem.
For example we have $$\dim_K K[N]/K[N_+]^{[q]}=\dim_K K[N/[q]N_+]=\# N/[q]N_+. $$ 
So the computation of the Hilbert-Kunz function and multiplicity of a binoid is in some sense easier than the  standard Hilbert-Kunz function and multiplicity. 

Also in Chapter 2, we develop the Hilbert-Kunz theory of binoids and prove the following main theorem.
\intr
Let $N$ be a finitely generated, semipositive, cancellative, reduced binoid and $\mathfrak{n}$ be an $N_+$-primary ideal of $N$. Then $e_{HK}(\mathfrak{n},N)$ exists and is a rational number.
\intre
Note that this is a characteristic free statement and that the existence does not follow from Monsky's theorem. From this result we can deduce a positive answers to our five problems.

To prove this theorem we use several reduction methods: If all less dimensional finitely generated, cancellative, semipositive (torsion-free) and integral binoids have a rational Hilbert-Kunz multiplicity then
\begin{enumerate}
 \item Non integral but reduced binoids have a rational Hilbert-Kunz multiplicity
 \item for non reduced binoids there is at least a bound for the Hilbert-Kunz multiplicity.
\end{enumerate}
To prove these reduction methods, we develop the notion of exact and strongly exact sequences of $N$-sets and their properties.

In Chapter 3, we will study the partitions of the $N$-sets $\qn$. K.-I.Watanabe showed in his paper (\cite{Watanabe}), that for toric rings there are finitely many non-isomorphic irreducible components in these partitions. So we will prove the same result for finitely generated, positive, integral, torsion-free and cancellative binoids. For that we will study irreducible $N$-sets, partition of $N$-sets and their normalization. 

Almost all our results are for cancellative binoids. One problem in the non-cancellative case is already the dimension of a binoid. The right candidate for the dimension of a finitely generated binoid $N$ is $$\max\{\rk \Gamma_\mathfrak{p}\mid \mathfrak{p}\text{ prime ideal of } N\}$$ where $\Gamma_\mathfrak{p}$ is the difference group of $N/\mathfrak{p}$, and the right candidate for doing Hilbert-Kunz theory in a semipositive, finitely generated binoid is $$\dim N=\max\{\rk \Gamma_\mathfrak{p}\mid \mathfrak{p}\text{ prime ideal}, N_+ \text{ lies on the cancellative component of $\mathfrak{p}$}\}.$$
We would like to prove the theorem also in this case. From our result we get at least, by arguing in the ring setting, that the Hilbert-Kunz multiplicity for reduced binoid algebras in fixed positive characteristic is rational.

For non reduced binoids we have only a bound. Moreover, if we can compute the Hilbert-Kunz function of any non-cancellative and reduced binoid then we can compute it for non-cancellative and non-reduced binoids (see Proposition \ref{make integral}) as well.

\chapter{Basics}

In this chapter we will introduce binoids, describe basic properties of binoids and binoid sets and their properties. For a general introduction to binoids we refer to \cite{Simone}, where most of the basic concepts were developed. 
We focus on material which is relevant for Hilbert-Kunz theory.

\section{Binoids and their properties}
\df
Let $N$ be a monoid (an additive semigroup with identity 0), then we call an element $\infty\in N$  \emph{absorbing} if $\infty+n=n+\infty=\infty$ for all $n\in N$.
\dfe
\df
A \emph{binoid} $(N,+,0,\infty)$ is a monoid with an absorbing element $\infty$. We write $N^\bullet$ for the set $N\setminus \{\infty\}$.
\dfe
\ex
The binoid $\{\infty\}$, which means $0=\infty$, is called the \emph{zero} binoid and $\{0,\infty\}$ the \emph{trivial} binoid.
\exe
\df
Let $N$ be a binoid and $M$ be submonoid of $N$ with $\infty$ then we say that $M$ is a \emph{subbinoid} of $N$.
\dfe
\df An element $u$ in a binoid $N$ is a \emph{unit} if there exists an element $a\in N$ such that $a+u=u+a=0$. The set of all units, denoted by $N^\times$, is a submonoid of $N$ which is a group, the \emph{unit group} of $N$. The set of all non-units $N\setminus N^\times$ will be denoted by $N_+$.
\dfe
\df
Let $N$ be a nonzero binoid then we say that $N$ is \emph{semipositive}, if $|N^\times|$ is finite and \emph{positive}, if $N^\times=\{0\}$.
\dfe
\df
Let $M$ and $N$ be binoids. A map $\phi:M\rightarrow N$ is a binoid $homomorphism$ if it is a monoid homomorphism which sends $\infty_M$ to $\infty_N$.
The set $\phi(M)=:\im \phi$ is the \emph{image} and $\{m\in M\mid \phi(m)=\infty_N\}=:\ker \phi$ is the \emph{kernel} of $\phi$. The set of all binoid homomorphisms from $M$ to $N$ is denoted by $\hom (M,N)$.
\dfe
\df A nonempty subset I of a binoid N is called an \emph{ideal} of N if both $I+N$ and $N+I$ are contained in $I$.
\dfe
\ex
Let $N$ be a binoid then $N_+$ is an ideal of $N$. Because if $a\in N_+$ and $n\in N$ then $a+n$ can not be a unit. So $a+n\in N_+$.
\exe
\df
Let $N$ be a commutative binoid and $I$ an ideal of $N$. Then the \emph{radical} of $I$, denoted by $\rad(I)$, is the set of all $a\in N$ such that $na\in I$, for some $n\in \n$. An ideal $I\subseteq N$ with $\rad(I)=I$ is called a \emph{radical ideal}.
\dfe
\rem
If $N$ is not a commutative binoid then $\rad(I)$ is not always an ideal. But if $N$ is a commutative binoid then it is an ideal, because if $x\in \rad(I)$ and $y\in N$ then there exists $n\in \n $ such that $nx\in I$, so $n(x+y)=nx+ny\in I$.
\reme

\df
Let $N$ be a binoid and $I$ be an ideal of $N$. The relation $\sim$ in $N$ defined by $$x\sim y \Leftrightarrow x=y \text{ or } x,y\in I$$ is called the \emph{Rees congruence} on $N$. The quotient $N/ \sim$ is called the \emph{quotient} of $N$ modulo $I$ and is denoted by $N/I$.
\dfe
\df For every ideal $I$ and $n\in \n_+$ we will denote $$nI:=\{a_1+a_2+\cdots+a_n\mid a_i \in I\}$$ and $$[n]I:=\langle na\mid a\in I\rangle$$ the minimal ideal including the set $\{na\mid a\in I\}$.

It is easy to see that $nI$ is an ideal of $N$. We call $nI$ the $n$th \emph{sum} of the ideal and $[n]I$ the $n$th \emph{Frobenius sum} of the ideal.
\dfe
\prop
\label{ordinaryforbinuos}
Let $N$ be a commutative binoid and $I$ be a finitely generated ideal of $N$. Then
$$\bigcap_n nI=\bigcap_q [q]I.$$
\prope
\pr
The inclusion $qI\supseteq [q]I$ is always true. So $$\bigcap_n nI \supseteq \bigcap_q [q]I.$$
Take $x \in \bigcap_n nI$ and fix $q$. If $m$ is the number of generators of $I$ then $x\in qmI$ implies that $x=s_1e_1+\cdots+s_me_m$ and $s_1+\cdots+s_m\geqslant qm$. So by the Pigeonhole principle there exists $s_i\geqslant q$ and from here we get $x\in [q]I$. Hence  $\bigcap_n nI \subseteq \bigcap_q [q]I.$
\pre
For a binoid homomorphism $\varphi: N\rightarrow M$ and an ideal $I\subseteq N$ we denote by $I+M$ the ideal generated $\varphi(I)$ and call it the \emph{extended ideal}. Since (for $n\in \n$) $[n]:N\rightarrow N, x\mapsto nx$, is a binoid homomorphism, the ideal $[n]I$ can be considered as the extended ideal under this homomorphism.
\lem
\label{extended ideal commutes with q}
Let $N$, $M$ be commutative binoids, $I$ an ideal of $N$ and $\varphi:N\rightarrow M$ a binoid homomorphism. Then $$[q]I+M=[q](I+M).$$
\leme
\pr
This follows from the commutative diagram
\[\begin{tikzpicture}[node distance=1.8cm, auto]
  \node (N) {$N$};
  \node (M) [right of=N] {$M$};
  \node (qN) [below of=N] {$N$};
  \node (qM) [below of=M] {$M$};
  \draw[->] (N) to node {$\varphi$} (M);
  \draw[->] (qN) to node {$\varphi$} (qM);
  \draw[->] (N) to node  [swap]{$[q]$} (qN);
  \draw[->] (M) to node [swap] {$[q]$} (qM);
\end{tikzpicture}
.\qedhere\]
\pre
\df
Let $N$ be a commutative binoid. Then we say that $\mathfrak{n}$ is an $N_+$-\emph{primary ideal} if $\rad(\mathfrak{n})=N_+$.
\dfe
\prop
Let $N$ be a commutative binoid and $\mathfrak{n}$ an $N_+$-primary ideal. Then $[q]\mathfrak{n}$ is an $N_+$-primary ideal.
\prope
\pr
Let $a\in N_+$. Then there exists $k\in \n$ such that $ka\in \mathfrak{n}$, so $qka\in [q]\mathfrak{n}$.
\pre
\prop
\label{primary ideal quotient}
Let $N$ be a commutative, finitely generated, semipositive binoid and $\mathfrak{n}$ an $N_+$-primary ideal. Then $N/\mathfrak{n}$ is a finite set.
\prope
\pr
Let $a_1,\dots,a_k$ be binoid generators of $N$. Then we can say that $a_1,\dots,a_s\in N_+$ and $a_{s+1},\dots,a_k\in N^\times$. So there exist $n_1,\dots,n_s\in\n_+$ such that $n_ia_i\in \mathfrak{n}$ for every $1\leq i\leq s$. Let $a\in N/\mathfrak{n}$ be represented by $a\in N$. There are only finitely many $a\in N^\times$. So let $a\in N_+$.
Then $a=\sum_{i=1}^s m_ia_i+b$, where $b\in N^\times$.
If $m_i\geq n_i$ for some $i$ then $a\in \mathfrak{n}$ and $a=\infty$ in $N/\mathfrak{n}$. So $m_i< n_i$ for all $1\leq i\leq s$. Hence there are only finitely many different $a\in N/\mathfrak{n}$.
\pre

\df Let $N$ be a commutative binoid. An ideal $\mathfrak{p}\neq N$ is called \emph{prime} if for every $a,b\in N$ and $a+b \in \mathfrak{p}$ implies that $a\in \mathfrak{p}$ or $b\in \mathfrak{p}$. The \emph{spectrum} of $N$, denoted by $\spec N$, is the set of all prime ideals of $N$.
\dfe
We have a bijective correspondence between $\spec N$ and $\hom (N,\{0,\infty\}).$ Indeed, if $\mathfrak{p}\in \spec N$ then we can define a homomorphism $\chi:N \rightarrow \{0,\infty\}$ by $$a\longmapsto \begin{cases} \infty, & \operatorname{if} a\in \mathfrak{p}, \\ 0, & \operatorname{if} a\notin \mathfrak{p}, \end{cases}$$
and if there is a binoid homomorphism $\phi:N\rightarrow \{0,\infty\}$ then it is easy to check that $\ker \phi$ is a prime ideal of $N$.

\df
Let $N$ be a commutative binoid and $I$ an ideal of $N$. A prime ideal containing $I$ is a \emph{minimal prime} of $I$ if there is no prime ideal $\mathfrak{q}$ with $I\subseteq \mathfrak{q}\subsetneq \mathfrak{p}$. The minimal primes of the zero ideal $\{\infty\}$ are called \emph{minimal prime ideals} of $N$ and the set of all of them is denoted by $\min N$.
\dfe
\df
Let $N$ be a binoid, then we say that $N$ is a \emph{separated} binoid if $$\bigcap_{n\in\n} nN_+=\{\infty\}.$$
\dfe

\prop
If a binoid $N$ has a positive $\n$-grading then it is separated.
\prope
\pr
Every $f\in N_+,f\neq \infty,$ has by positivity $\deg(f)=k>0$. If we take $n>k$ and $g\in nN_+$ then $\deg(g)=\deg(f_1+\cdots+f_n)=\sum_{i=1}^n \deg(f_i)\geqslant n>k$. This implies that $f\not\in nN_+$ hence $\bigcap nN_+=\infty$.
\pre
\df
Let $N$ be a commutative binoid, then we set $$N^{\sep}:=N/\bigcap_n nN_+$$ and call it the \emph{separation} of $N$.
\dfe
The separation of $N$ is clearly a separated binoid.
\df
Let $N$ be a commutative binoid. An element $a\in N$ is \emph{nilpotent} if $na=a+\cdots+a=\infty$ for some $n\in\n$. The set of all nilpotent elements will be denoted by $\nil(N)$. We say that $N$ is  \emph{reduced} if $\nil(N)=\{\infty\}$.
\dfe
It is easy to see that $\nil(N)$ is an ideal of $N$.
\df
Let $N$ be a commutative binoid. We call the quotient binoid $N_{\red}:=N/\nil(N)$ the \emph{reduction} of $N$.
\dfe
The reduction of $N$ is clearly a reduced binoid.
\df
An element $a$ in a binoid $N$ is a \emph{torsion} element in case $a=\infty$ or $na = nb$ for some $b\in N, b\neq a, n \geqslant 2$. 
We say $N$ is \emph{torsion-free} if there are no other torsion elements in M besides $\infty$, i.e. $na = nb$ implies $a = b$ for every $a, b \in N$ and $n\geqslant 1$.
A binoid is called \emph{torsion-free up to nilpotence} if $na = nb \neq \infty$ implies $a = b$ for
every $a, b \in N$ and $n \geqslant 1$.
\dfe
\df
An \emph{integral} element $a\in N^\bullet$ satisfies the property that $a+b=\infty$ or $b+a=\infty$ implies $b=\infty$. The set of all integral elements will be denoted by $\operatorname{int}(N)$ and the complement $N\setminus \operatorname{int}(N)$ by $\operatorname{int}^c(N)$. Let $N$ be a commutative binoid and every $a\in N^\bullet$ is integral then we call it \emph{integral binoid}.
\dfe
For the definition of a localization we refer to \cite[Subsection 1.13]{Simone}.
\df
Let $N$ be a commutative binoid and $M$ the submonoid of $N$ that is generated by $N^\bullet$. The binoid $N_M=:\diff N$ will be called the \emph{difference} binoid of $N$. Here $N_M$ is the localization of $N$ at $M$.
\dfe
\rem
If $N$ is a non-integral binoid then $\infty\in M$, so $\diff N=\{\infty\}$. If a binoid is integral then $\diff N=N_{N^\bullet}$.
\reme
\df
Let $N\subseteq M$ be commutative binoids. If $M\subseteq \diff N$ then we say that $M$ is \emph{birational} over $N$.
\dfe
\df
Let $N$ be a integral binoid. The \emph{rank} of a binoid $N$ is the rank of the difference group $(\diff N)^\bullet$, i.e. the vector space dimension of $\q\otimes_\z (\diff N)^\bullet$ over $\q$.
\dfe

\section{$N$-sets}
In this section arbitrary binoids are assumed to be commutative. 
\df
A \emph{pointed set} $(S,p)$ is a set $S$ with a distinguished element $p\in S$.
\dfe
Let $N$ be a binoid.
\df
An \emph{operation} of $N$ on a pointed set $(S,p)$ is a map $$+:N\times S \longrightarrow S,~~ (n,s)\longmapsto n+s,$$ such that the following conditions are fulfilled:
\begin{enumerate}
  \item For all $s\in S: 0+s=s$.
  \item For all $s\in S: \infty+s=p$.
  \item For all $n\in N: n+p=p$.
  \item For all $n,m\in N$ and $s\in S:(n+m)+s=n+(m+s)$.
\end{enumerate}
Then $S$ is called an $N$-\emph{set}.
\dfe
\ex
Let $N$ be a binoid, then we can think of $N=(N,\infty)$ as an $N$-set with the addition as the operation.
\exe
\ex
\label{n set}
Let $N=\n^\infty:=\n\cup\{\infty\}$, then an $N$-set $S$ is the same as the set $S$ together with a map $f:S\rightarrow S$. The operation is given by $n+s=f^n(s)$ and the map is given by $f(s)=1+s$. More generally if $N=(\n^k)^\infty$ then an $N$-set $S$ is a set $S$ together with $k$ commuting maps $$f_1,\dots,f_k:S\rightarrow S,$$ meaning that $f_i\circ f_j=f_j\circ f_i$ for every $i,j$. The operation is $(n_1,\dots,n_k)+s=f_k^{n_k}\circ\cdots \circ f_1^{n_1}(s)$ and the $i$th map is $f_i(s)=e_i+s$.
\exe
\df
Let $N$ be a binoid and $(S,p)$ be an $N$-set. We say that $a\in N$ is an \emph{annihilator} if $a+s=p$ for every $s\in S$. We denote the set of all annihilators by $\ann S$.
\dfe
\lem
Let $N$ be a binoid and $(S,p)$ be an $N$-set, then $\ann S$ is an ideal of $N$.
\leme
\pr
It is clear that $\infty \in \ann S$ and if $a\in\ann S$ then for every $n\in N,s\in S$ we have $(a+n)+s=a+(n+s)=p$, so $a+n\in \ann S$.
\pre
\lem
\label{N/ann-set}
Let $N$ be a binoid and $(S,p)$ be an $N$-set. If $\mathfrak{a}\subseteq \ann S$ is an ideal of $N$ then $(S,p)$ is also an $(N/\mathfrak{a})$-set.
\leme
\pr
  We can define an operation $(N/\mathfrak{a})\times S\rightarrow S$ in the following way. If $n\notin \mathfrak{a}$ then $[n]+s=n+s$, otherwise $[\infty]+s=p$. So then it is easy to check that $(S,p)$ is an $(N/\mathfrak{a})$-set.
\pre
\lem
\label{homomorphism gives N set}
Let $N,M$ be binoids and let $\phi:N\rightarrow M$ be a binoid homomorphism. Then we can consider $M$ as an $N$-set.
\leme
\pr
We define the operation of $N$ on $M$ by $$+:N\times M \longrightarrow M,~~ (n,m)\longmapsto \phi(n)+m.$$ Then $\phi(0_N)+m=0_M+m=m$, $\phi(\infty_N)+m=\infty_M+m=\infty_M$, $\phi(n)+\infty_M=\infty_M$, $\phi(n_1+n_2)+m=(\phi(n_1)+\phi(n_2))+m=\phi(n_1)+(\phi(n_2)+m).$
\pre
\ex
Let $N$ be an integral binoid, $T=\{0,\infty\}$ the trivial binoid. Then we have the following binoid homomorphisms
\begin{enumerate}
 \item $N\longrightarrow N/N_+\longrightarrow T$, $$n\longmapsto \begin{cases} 0, & \operatorname{if} n\in N^\times, \\ \infty, & \operatorname{if} n\in N_+. \end{cases}$$
 \item $N\longrightarrow \diff N\longrightarrow T$, $$n\longmapsto \begin{cases} 0, & \operatorname{if} n\in N^\bullet, \\ \infty, & \operatorname{if} n=\infty. \end{cases}$$
\end{enumerate}
So there are two different $N$-sets $T$. In fact every prime ideal $\mathfrak{p}$ of $N$ gives an $N$-operation of $T$.
\exe
\df
Let $N$ be a binoid, $(S,p)$ be an $N$-set. We say that $T\subseteq S$ is an $N$-\emph{subset} of $S$, if $p\in T$ and for all $t\in T$ and $a\in N$ we have $a+t\in T$. $S$ and $\{p\}$ are called the trivial $N$-subsets.
\dfe
\lem
\label{union and intersection}
Let $N$ be a binoid, $(S,p)$ be an $N$-set. If $T_i\subseteq S, i\in I$, are $N$-subsets then $\bigcup_{i\in I}T_i$ and $\bigcap_{i\in I}T_i$ are also $N$-subsets.
\leme
\pr
$p\in T_i$ for every $i\in I$ so $p\in \bigcup_{i\in I}T_i$ and $p\in \bigcap_{i\in I}T_i$. Also for all $a\in N$ and $t\in \bigcup_{i\in I}T_i$($t\in \bigcap_{i\in I}T_i$) we have $a+t\in \bigcup_{i\in I}T_i$ ($a+t\in \bigcap_{i\in I}T_i$).
\pre
\rem
An ideal in $N$ is the same as an $N$-subset of $N$. ($N$ considered as an $N$-set)
\reme
\prop
\label{N/I Nset}
Let $N$ be a binoid and $I$ be an ideal of $N$, then $N/I$ is an $N$-set .
\prope
\pr
Let $$ N\longrightarrow N/I, n\longmapsto [n],$$ where $n\in N,[n]\in N/I$. Then this is a binoid homomorphism. So by Lemma \ref{homomorphism gives N set} we have the result.
\pre
\df
Let $N$ be a binoid and $(S,p)$ be an $N$-set. If there exist finitely many elements  $s_1,\dots,s_k\in S$ such that for all $s\in S$ we can write $s=n+s_j$ where $n\in N$, then we say that $S$ is a \emph{finitely generated $N$-set}. We call the elements $s_j$ $N$-\emph{generators}.
\dfe
\df
Let $N$ be a binoid and $(S,p)$ be an $N$-set. We call an element $t\in S$ \emph{reducible}, if we can write $t=n+s,n\in N,t\neq s\in S$, otherwise we call it \emph{irreducible}.
\dfe
\prop
\label{minimal generating set of N-set}
Let $N$ be a positive binoid, $(S,p)$ be a finitely generated $N$-set and let $\{s_1,\dots,s_k\}$ be a minimal generating set of $S$. Then the $s_i$ are irreducible and the set $\{s_1,\dots,s_k\}$ is uniquely determined and we call it $\gen_N S$.
\prope
\pr
If $s_i$ is reducible then there exist $n\in N$ and $s_i\neq s_j\in S$ such that $s_i=n+s_j$. So we do not need $s_i$ in the generating set. This contradicts that $\{s_1,\dots,s_k\}$ is the minimal generating set. Also if $s$ is an irreducible element then $s\in\{s_1,\dots,s_k\}$. So it is the set of all irreducible elements of $S$.
\pre
\prop
\label{Subset of f.g.}
Let $N$ be a finitely generated binoid and $S$ be a finitely generated $N$-set. If $T$ is an $N$-subset of $S$ then it is finitely generated.
\prope
\pr
Let $s_1,\dots,s_k$ be generators of $S$. If we set $$(s_i)=\{n\in N\mid n+s_i\in T\},i=1,\dots,k,$$ then $(s_i)$ is an ideal of $N$. By \cite[Proposition 2.1.14]{Simone} the ideals of a finitely generated binoid are finitely generated. So if $n_{i1},\dots,n_{im_i}$ are the generators of $(s_i)$ then $\{n_{ij}+s_i\mid i=1,\dots,k,j=1,\dots,m_i\}$ generate $T$.
\pre
\df
Let $N$ be a binoid. For a finite $N$-set $S$ we set $\# S=|S|-1=|S\setminus\{p\}|$, so we do not count the distinguished point.
\dfe
\prop
\label{birational to m}
Let $N\subseteq M$ be integral binoids. Let $M$ be birational over $N$ and finitely generated as an $N$-set. Then there exist $m\in N^\bullet$ such that $m+M\subseteq N$.
\prope
\pr
Let $m_1,\dots,m_r$ be the $N$-generators of $M$. Then by birationality we can write $m_i=a_i-b_i$, where $a_i,b_i\in N,b_i\neq \infty$. Hence if $m:=\sum_{i=1}^r b_i$ then $m+M\subseteq N$.
\pre

\section{Partition of N-sets}
In this section arbitrary binoids are assumed to be commutative. 
\df
Let $N$ be a binoid and let $(S_i,p_i),i\in I$, be a family of $N$-sets. Then we define the \emph{pointed union} of $S_i,i\in I,$ by $$\bigcupdot_{i\in I} S_i=(\biguplus_{i\in I}S_i)/\sim,$$ where $\biguplus$ is a disjoint union and $a\sim b$ if and only if $a=b$ or $a=p_j,~b=p_k$ for some $j,k$. So the pointed union just contracts the points $p_j$ to one point.
We write $S^{\cupdot r}=\bigcupdot_{i=1}^r S$ and in particular $N^{\cupdot r}$ for the $r$-folded pointed union of $N$ with itself. If $n$ is considered in the $j$-th component then we denote it by $(j,n)$.
\dfe
\df
Let $N$ be a binoid and $(S,p)\neq \{p\}$ be an $N$-set. $S$ is \emph{reducible}, if it can be written as a pointed union of non-trivial $N$-subsets. So \emph{irreducible}, if it is not the pointed union of non-trivial $N$-subsets.
\dfe

\df
Let $N$ be a binoid and $(S,p)$ be a $N$-set. $T\subseteq S$ is called an \emph{irreducible component} of $S$ if it is an irreducible $N$-subset and every $N$-set $T'$ with $T\subsetneq T'\subseteq S$ is reducible.
\dfe
 Let $N$ be a binoid and $(S,p)$ be an $N$-set. We can define an equivalence relation between two elements $s,t\in S^\bullet:=S\setminus \{p\}$ by the following: $t\sim_N s$ if and only if there exist elements $s=s_0,s_1,\dots,s_{k-1},t=s_k\in S^\bullet$ and $n_1,\dots,n_k\in N^\bullet$ such that $n_j+s_{j-1}=s_j$ or $s_{j-1}=n_j+s_j$ for all $j=1,\dots,k$. 
 It is easy to check that this is an equivalence relation. So from here we have $$S^\bullet=\biguplus_{i\in I} T_i$$ where $T_i$ are the equivalence classes. Equivalence classes are either equal or disjoint. So if we denote $S_i:=T_i\cup \{p\}$ then it will give us a pointed union $$S=\bigcupdot_{i\in I} S_i.$$ 
 If the index set $I$ is finite then we say that $S$ has a \emph{finite partition}.
 \prop
 Let $N$ be a binoid and $(S,p)$ be an $N$-set. If $S^\bullet=\biguplus_{i\in I} T_i$, where $T_i$ are its equivalence classes, then $S_i=T_i\cup \{p\}$ is an $N$-subset of $S$.
 \prope
 \pr
 Let $s\in S_i$ and $n\in N$, then either $n+s=p$ or $n+s\sim_N s$ so $n+s\in S_i$ in both cases.
 \pre
 \prop
 \label{irr set has one component}
 Let $N$ be a binoid and $(S,p)\neq \{p\}$ be an $N$-set. Then $S$ is irreducible if and only if it has only one equivalence class.
 \prope
 \pr
 Let $S$ be an irreducible $N$-set, then we have $S^\bullet=\biguplus_{i\in I} T_i$ and if we have more than one equivalence class then $S$ is the pointed union of non-trivial $N$-subsets, hence it gives us a contradiction.

 Now suppose that $S$ has only one equivalence class, which means that every two elements of $S^\bullet$ are equivalent to each other. If $S$ is reducible then $S=R\cupdot T$ and $R,T\neq\{p\}$. Hence there exist $s\in R^\bullet,t\in T^\bullet$ such that  $s\sim_N t$. 
 So there exist $s=s_0,s_1,\dots,s_{k-1},t=s_k\in S^\bullet=R^\bullet\uplus T^\bullet$ and $n_1,\dots,n_k\in N^\bullet$ such that $n_j+s_{j-1}=s_j$ or $s_{j-1}=n_j+s_j$ for all $j=1,\dots,k$. Hence there exists $j$ such that $s_j\in R$ and $s_{j+1}\in T$, so it gives us a contradiction.
 \pre
 \ex
Let $N$ be an integral binoid, then it is irreducible as an $N$-set. This is because $n,m\in N^\bullet$ are connected by $n+m\neq\infty$.
\exe
 \prop
 \label{equivalence class}
 Let $N$ be a binoid, $(S,p)$ be an $N$-set, then the equivalence classes of $S$ are exactly the irreducible components of $S$.
\prope
\pr
Let $S=\bigcupdot S_i$ where $S_i\setminus \{p\}$ are the equivalence classes of $S^\bullet$. By Proposition \ref{irr set has one component} $S_i$ is an irreducible $N$-subset. Now let $T$ be $S_i \subsetneq T\subseteq S$, then there exists $t\in T$ such that $t\notin S_i$. Hence $t\in S'=\bigcup_{j\neq i}S_j$ and from Lemma \ref{union and intersection}, $T\cap S'$ is $N$-set. So $T=S_i\cupdot (T\cap S')$ which means it is reducible.

Let $T\subseteq S$ be an irreducible component of $S$. By Proposition \ref{irr set has one component} we know that $T$ has only one equivalence class, so there exists $S_i$ such that $T\subseteq S_i$. But by Proposition \ref{irr set has one component} we also know that $S_i$ is an irreducible $N$-set. Hence $T=S_i$.
\pre
\prop
Let $N$ be a binoid and $(S,p)$ be an $N$-set. If $T$ is an irreducible $N$-subset of $S$ then $T\subseteq S_i$, where $S=\bigcupdot_{i\in I}S_i$ is the partition with the irreducible components.
\prope
\pr
If $t\in T\subseteq S$ then there exists $S_i$ such that $t\in S_i$. Now take any element $s\in T$. By Proposition \ref{irr set has one component} $t\sim_N s$, so this means $s\in S_i$.
\pre
\prop
Let $S$ be finitely generated over $N$, then it has a finite partition into irreducible $N$-subsets.
\prope
\pr
Let $s_1,\dots,s_k$ be the generators of $S$. If $s\in S_j$ where $S_j$ is one of the irreducible components of S, then $N+s\subseteq S_i$. Also $N+s_1,\dots,N+s_k$ will cover $S$. So the number of irreducible components is less or equal to $k$.
\pre

\section{Normalization of a binoid}
In this section arbitrary binoids are assumed to be commutative. 
\df
Let $M\subseteq N$ be binoids. Then we say that $N$ is \emph{pure integral} over $M$, if for every $n\in N$ there exist $k\in\n_+$ such that $kn\in M$.
\dfe
\df
Let $M\subseteq N$ be binoids. Then we say that $N$ is \emph{finite} over $M$, if $N$ is a finitely generated $M$-set.
\dfe
\lem
\label{finite set}
Let $M\subseteq N$ be binoids and suppose that $N$ is cancellative. Then $N$ is finite over $M$ if and only if $N$ is pure integral over $M$ and is a finitely generated $M$-binoid.
\leme
\pr
Let $N=M[x_1,\dots,x_t]=\{\sum_{i=1}^t n_ix_i+m\mid m\in M,n_i\in\n\}$, which means finitely generated as $M$-binoid. Also by pure integrality we have $k_ix_i\in M$, where $k_i\in\n_+$. Hence the elements of the form $\sum_{i=1}^t r_ix_i$ with $0\leqslant r_i<k_i$ are $M$-set generators of $N$, since we may write $n_i=a_ik_i+r_i,0\leqslant r_i<k_i$, and so $$\sum_{i=1}^t n_ix_i+m=\sum_{i=1}^t r_ix_i+\sum_{i=1}^t a_ik_ix_i+m.$$

Now let $y_1,\dots,y_s$ be generators of $N$ as an $M$-set. These are in particular generators of $N$ as an $M$-binoid. Let's take any element $y\in N$. So $y+y_i=y_{j(i)}+m_i,1\leqslant i\leqslant s$ and $m_i\in M$. If $y=y_t+m$ then $2y=y+y_t+m=y_{j(t)}+m_t+m$ and also it is not difficult to see that for every $k$ $$(k+1)y=y_{j^k(t)}+m_{j^{k-1}(t)}+\cdots+m_{j(t)}+m_t+m.$$
But we only have $s$ different $y_i$'s, so there exist $n,l\in\n_+$ such that $y_{j^n(t)}=y_{j^{n+l}(t)}$. Hence
\begin{align*}
(n+l+1)y&=y_{j^{n+l}(t)}+m_{j^{n+l-1}(t)}+\cdots+m_{j^{n+1}(t)}+m_{j^n(t)}+\cdots+m_{j(t)}+m_t+m\\
 &=(n+1)y+m_{j^{n+l-1}(t)}+\cdots+m_{j^{n+1}(t)}.
 \end{align*}
 So by cancellation we have $ly=m_{j^{n+l-1}(t)}+\cdots+m_{j^{n+1}(t)}\in M$.
\pre
Note: If our binoid is not cancellative then we only have the if part.
\ex
Let $M=\n^\infty=(x),N=M[y]/(2y=x+y)$. Then $\{y\}$ is an $M$-set generating system of $N$ but for any $k\in\n$, $ky\notin M$.
\exe
\df
Let $N$ be an integral and cancellative binoid and let $\diff N=N_{N^\bullet}$ be the difference group binoid of $N$. The binoid $$\hat{N}:=N^{normal}=\{x\in \diff N\mid  nx\in N,\operatorname{\;for\;some\;} n\in \n_+\}$$ is called the \emph{normalization} of $N$. We say that $N$ is \emph{normal} if $N=\hat{N}$.
\dfe

Let $N$ be a finitely generated, cancellative, integral and torsion-free binoid. We can assume that $(\diff N)^\bullet=\z^d$, where $d=\rk(N)$. Let $C=\r_+N^\bullet\subseteq \r^d$ then $C$ is the rational cone generated by $N$. This cone has by \cite[Theorem 1.15]{BrunsGubeladze}, a representation $$C=H^+_{\sigma_1}\cap\cdots\cap H^+_{\sigma_s}$$ as an irredundant intersection of half-spaces defined by rational linear forms $\sigma_i$ on $\r^d$.
Each of the hyperplanes $H_{\sigma_i}=\{a\in \r^d\mid \sigma_i(a)=0\}$ is generated as a vector space by integral vectors and we call them \emph{supporting hyperplanes}. Therefore we can assume that $\sigma_i$ is a linear form on $\z^d$ with nonnegative values on $C$ and we call it \emph{supporting forms}.

\prop
\label{normalization and k}
Let $N$ be a finitely generated, cancellative, integral, torsion-free binoid and $\hat{N}$ be the normalization of $N$. Then there exists $k\in \n,x_1,\dots,x_s\in \hat{N}$ such that $kx\in N$ for every $x\in \hat{N}$ and $\hat{N}=N\cup(x_1+N)\cup\cdots\cup(x_s+N)$.
\prope
\pr
Since $N$ is torsion-free, its difference group is $\diff N\cong\z^r$. If $n_1,\dots,n_l$ are the generators of $N$, then $\r_{\geqslant 0}N=\r_{\geqslant 0} \langle n_1,\dots,n_l\rangle$ is a finitely generated rational cone in $\r^r$ and so by \cite[Corollary 2.24]{BrunsGubeladze}, $\hat{N}^\bullet=\r_+N^\bullet\cap \z^r$ and it is finitely generated monoid by Gordan's lemma.
In particular, $\hat{N}$ is a finitely generated binoid. Since it is also pure integral and cancellative, it is finite as an $N$-set by Lemma \ref{finite set}.
Hence there exists $s\in\n$ such that $\hat{N}=N\cup(x_1+N)\cup\cdots\cup(x_s+N)$. Now $x_i\in \hat{N}$ implies that there exists $k_i\in \n$ such that $k_ix_i\in N$ for every $1\leqslant i \leqslant s$ and if we take $k=k_1k_2\cdots k_s$ then $kx\in N$ for every $x\in \hat{N}$.
\pre
\rem
Let $N$ be a finitely generated, positive, cancellative, integral and torsion-free binoid, then $\hat{N}$ is a finitely generated $N$-set and a finitely generated binoid.
\reme
\ex
\label{double N}
Let $a=(2,1),b=(3,0)\in (\n \times \z/2)^\infty$ be the generators of a binoid $N$. Then $(\diff N)^\times=\z \times \z/2$, since $2(2,1)-(3,0)=(1,0)$ and $3(2,1)-2(3,0)=(0,1)$. This binoid has torsion elements, since for example $(6,1)\neq(6,0)$ but $2(6,1)=2(6,0)$. The binoid is positive. Its normalization contains $(1,0)$ and $(0,1)$, hence $\hat{N}=(\n\times \z/2)^\infty$, which is not positive, only semipositive.
\exe

\prop
\label{normalization type}
Let $N$ be a finitely generated, positive, integral, cancellative and torsion-free binoid and $\hat{N}$ be the normalization of $N$. If $\hat{N}=N\cup(x_1+N)\cup\cdots\cup(x_s+N)$ then there exists $\alpha,p_1,\dots,p_s\in\n_+$ such that for all $\n\ni m\geqslant\alpha$ we have $p_imx_i\in N$ and if $p_i\neq 1$ then $(p_im+1)x_i,\dots,(p_im+p_i-1)x_i\notin N$ for $1\leqslant i \leqslant s$.
\prope
\pr
Let $A_i=\{a\in\n\mid ax_i\in N\},1\leqslant i \leqslant s$ and we know that these are the numerical semigroups. If we take $p_i:=\gcd(A_i,1\leqslant i \leqslant s)$ then it is known that there exists $\alpha_i\in\n$ such that for all $m\geqslant\alpha_i$ we have $p_imx_i\in N$, here $\gcd(A_i,1\leqslant i \leqslant s)$ is the greatest common divisor of $A_i$'s.
If $p_i\neq 1$ then $(p_im+1)x_i,\dots,(p_im+p_i-1)x_i\notin N$, because $\gcd(p_i,p_im+j)<p_i$, for $1\leqslant j\leqslant p_i-1$.
So, if we take $\alpha=\max \alpha_i$ then we have our result.
\pre

\prop
\label{m}
Let $N$ be a finitely generated, positive, integral, cancellative and torsion-free binoid and $\hat{N}$ be the normalization of $N$. Then there exists $m\in N$ such that $m+\hat{N}\subseteq N$.
\prope
\pr
By Proposition \ref{normalization and k} we know that $\hat{N}$ is a finitely generated $N$-set and by definition $\hat{N}$ is birational over $N$. Hence our result follows from Proposition \ref{birational to m}.
\pre
Let $\hat{N}$ be the normalization of $N$, then the set of `gaps' $\hat{N}\setminus N$ is, in a sense, small.
\prop
\label{c}
Let $N$ be a finitely generated, positive, integral, cancellative and torsion-free binoid. Then the ideal $$c(\hat{N}\setminus N):=\{x\in N\mid x+\hat{N}\subseteq N\}$$ is nonempty. Moreover, the set $\hat{N}\setminus N$ is contained in finitely many hyperplanes parallel to the facets of $\hat{N}$.
\prope
\pr
By Proposition \ref{m} we have that $m\in c(\hat{N}\setminus N)$, so it is nonempty. If $x\in c(\hat{N}\setminus N)$ then $x+n+\hat{N}\subseteq N$ for every $n\in N$, so $x+n\in c(\hat{N}\setminus N)$. This means that $c(\hat{N}\setminus N)$ is an ideal of $N$. If $z\in\hat{N}\setminus N$, then there exists at least one supporting form $\sigma_i$ of $N$ with $0\leqslant\sigma_i(z)<\sigma_i(m)$. So $z$ belongs to one of the hyperplanes $\{y\mid \sigma_i(y)=k\},k\in\z_{\geqslant 0},k<\sigma_i(m)$.
\pre
We will call $c(\hat{N}\setminus N)$ $conductor$. It is easy to check that the conductor is the largest ideal of $N$ that is also an ideal of $\hat{N}$.

\section{Homomorphisms of $N$-sets}
In this section arbitrary binoids are assumed to be commutative. 
\df
Let $N$ be a binoid and let $(S,p_S),(T,p_T)$ be $N$-sets. Then we call a map $\phi:S\rightarrow T$ an $N$-set \emph{homomorphism} if $\phi(n+s)=n+\phi(s)$, for every $s\in S,n\in N$.
\dfe
\ex
\label{pointed degree}
Let $N$ be a binoid, $S_i,T_i,1\leqslant i\leqslant k$ be $N$-sets and suppose we have homomorphisms $\phi_i:S_i\rightarrow T_i$ for every $1\leqslant i\leqslant k$. Then we can construct a homomorphism $\phi:\bigcupdot S_i\rightarrow \bigcupdot T_i$ such that $\phi(i,t)=(i,\phi_i(t))$.
\exe
\df
Let $N$ be a binoid, $(S,p_S),(T,p_T)$ be $N$-sets. If we have a homomorphism $f:S\rightarrow T$  then $$\im(f)=\{t\in T \mid t = f(s) \text{ for some } s \in S\}$$ $$\ker(f)=\{s\in S \mid  f(s)=p_T\}.$$
\dfe
\df
Let $N$ be a binoid, let $(T,p)$ be an $N$-set and $S\subseteq T$ an $N$-subset of $T$. The \emph{quotient} of $T$ by $S$ is the $N$-set $(T\setminus S, S)$ and it is denoted by $T/S$.
\dfe
\rem Here we think that every element of $S$ is collapsing into the point $p$.
\reme
Hence by definition we have a map $N\times T/S\longrightarrow T/S$, $$(n,t)\longmapsto\begin{cases} n+t, & \operatorname{if} t\notin S, \\ p, & \operatorname{if} t\in S, \end{cases}$$ which fulfills the conditions of an $N$-set. Also, we have a canonical $N$-set homomorphism $\varphi:T\rightarrow T/S$, given by $\varphi(t)=t$, if $t\notin S$ and $\varphi(t)=p$ otherwise.
\rem
Let $N$ be a binoid and $I$ be an ideal of $N$. If we consider $N$ as an $N$-set then we have the  quotient $N/(I+N)$. But $I+N=I$ so in this case we will just denote it with $N/I$.
\reme

\lem
\label{hom}
Let $N$ be a binoid, $(S,p_S),(T,p_T)$ $N$-sets and $S'\subseteq S$ an $N$-subset. If we have an $N$-set homomorphism $\phi:S\rightarrow T$ with $\phi(S')=p_T$ then there exists a unique  homomorphism $\tilde{\phi}:S/S'\rightarrow T$ such that the following diagram commutes.

\[\begin{tikzpicture}[node distance=1.8cm, auto]
  \node (S) {$S$};
  \node (T) [right of=S] {$T$};
  \node (A) [below of=S] {$S/S'$};
  \draw[->] (S) to node {$\phi$} (T);
  \draw[->, dashed] (A) to node  [swap]{$\tilde{\phi}$} (T);
  \draw[->] (S) to node [swap] {$\varphi$} (A);
\end{tikzpicture}
\]
If $\varphi$ is surjective, then $\tilde{\phi}$ is surjective.
\leme
\pr
Since $\phi$ is surjective, there exists at most one $\tilde{\phi}$. So we construct $\tilde{\phi}$ in the following way. $$\tilde{\phi}(\overline{s})= \begin{cases} \phi(s), & \textrm{ if } s\notin S', \\ p_T, & \textrm{ else}. \end{cases}$$ Then it is easy to see that this map is a well defined $N$-set homomorphism and that the diagram commutes.
Let $\phi$ be a surjective homomorphism, then for any $t\in T$ there exists $s\in S$ such that $\phi(s)=t$. So by construction $\tilde{\phi}(\varphi(s))=t$.
\pre
\prop
\label{quotient to union}
Let $N$ be a finitely generated binoid, $J$ be an ideal of $N$ and $S\subseteq T$ be $N$-sets. Then $(T/S)/(J+(T/S))= T/(S\cup (J+T))$.
\prope
\pr
We have an $N$-set surjective homomorphism
\begin{align*} T/S &\xrightarrow{\;\phi\;} T/(S\cup (J+T)),\\ t &\longmapsto [t].
\end{align*}
If   $t\in T/S$ and $\phi(t)=\infty$ then $t\in S\cup (J+T)$ so $t\notin S$ implies that $t\in J+T$. Hence it follows $t\in J+(T/S)$. Now let $t\in J+(T/S)$, then $\phi(t)=\infty$, so $\ker \phi= J+(T/S)$. So Lemma \ref{hom} gives surjective homomorphism $$(T/S)/(J+(T/S))\longrightarrow T/(S\cup (J+T))$$ and it is injective by construction.
\pre
\lem
\label{canonical isomorphism}
Let $N$ be a binoid, $I\subseteq N$ an ideal of $N$ and $(T,p)$ be an $N$-set. If $r$ is some positive integer then we have a canonical isomorphism $$T^{\cupdot r}/(I+T^{\cupdot r})\cong (T/(I+T))^{\cupdot r}.$$
\leme
\pr
We know that $I+T^{\cupdot r}=\{(j,n+t)\mid n\in I, t\in T, 1\leqslant j \leqslant r\}$ and by Example \ref{pointed degree} there exists a surjective homomorphism $\phi:T^{\cupdot r}\rightarrow (T/(I+T))^{\cupdot r}$. Since $\phi(n+(j,t))=\phi(j,n+t)=p$ we have $\phi(I+T^{\cupdot r})=p$, for $n\in I, (j,t)\in T^{\cupdot r}$.
Hence by Lemma \ref{hom}, there exists a surjective homomorphism $$\tilde{\phi}:T^{\cupdot r}/(I+T^{\cupdot r})\longrightarrow (T/(I+T))^{\cupdot r}.$$ If $\tilde{\phi}(j,t)=\tilde{\phi}(j',t')$, then $j=j',t=t'$ or $\tilde{\phi}(j,t)=\tilde{\phi}(j',t')=p$. In the last case $t=n+s,t'=n'+s'$, where $n,n'\in I,s,s'\in T$ and so $(j,t)=(j',t')=p$. Hence, $\tilde{\phi}$ is injective.
\pre
\lem
\label{surj hom}
Let $N$ be a binoid, $(S,p_S),(T,p_T)$ $N$-sets and $I\subseteq N$ be an ideal of $N$. If we have a surjective $N$-set homomorphism $\phi:S\rightarrow T$ then there exists a canonical surjective $N$-set homomorphism $\tilde{\phi}:S/(I+S)\longrightarrow T/(I+T)$.
\leme
\pr
We have a canonical homomorphism $\varphi:T\rightarrow T/(I+T)$, so $$\varphi\circ\phi:S\longrightarrow T/(I+T)$$ will give us a homomorphism with the property that $(\varphi\circ\phi)(I+S)=p_T$.
Hence by Lemma \ref{hom} we have a homomorphism $\tilde{\phi}:S/(I+S)\longrightarrow T/(I+T)$. (If $t\notin I+T$, then by surjectivity of $\phi$, there exists $s\in S$ such that $\phi(s)=t$. Then it is easy to see $\tilde{\phi}(\varphi(s))=t$.)
\pre

\lem
\label{surj map for f.g. set}
Let $N$ be a binoid and $T$ be an $N$-set. Then
\begin{enumerate}
\item For $t_1,\dots,t_r\in T$ we can define an $N$-set homomorphism $\phi:N^{\cupdot r}\rightarrow T$.
\item $t_1,\dots,t_r\in T$ is a generating system of $T$ over $N$ if and only if $\phi$ is a surjective homomorphism.
\item $T$ is finitely generated over $N$ if and only if there exists a surjective $N$-set homomorphism $N^{\cupdot r}\rightarrow T$.
\end{enumerate}
\leme
\pr
\begin{enumerate}
\item Let's take $\phi(i,0)=t_i$ then $\phi(i,n)=n+t_i$ and  $\phi((i,n)+m)=\phi(i,n+m)=(n+m)+t_i=m+(n+t_i)$. So $\phi$ is an $N$-set homomorphism.
\item If $t_1,\dots,t_r\in T$ is a generating system of $T$ over $N$ then we can write every element $t\in T$ as $t=n+t_j$ for some $j$ and $n\in N$. Hence $t=\phi(j,n)$.
If $t_1,\dots,t_r\in T$ and $\phi$ is a surjective homomorphism then for every element $t\in T$ there exist $n\in N$ and $1\leqslant j\leqslant r$ such that $t=\phi(j,n)=n+t_j$.
\item If $T$ is a finitely generated over $N$ then by 1. and 2. we have a surjective $N$-set homomorphism $N^{\cupdot r}\rightarrow T$. If we have a surjective homomorphism $N^{\cupdot r}\rightarrow T$ then $\{\phi(i,0)\mid 1\leqslant i \leqslant r\}$ will give us a generating system of $T$.
\end{enumerate}
\pre

\section{Smash product of $N$-sets}
In this section arbitrary binoids are assumed to be commutative. 
\df
Let $(S_i,p_i)_{i\in I}$ be a finite family of pointed sets and $\sim_\wedge$ the relation on $\prod_{i\in I} S_i$ given by $$(s_i)_{i\in I}\sim_\wedge (t_i)_{i\in I}:\Leftrightarrow s_i=t_i,\forall i\in I, \;\operatorname{or}\; s_j=p_j,t_k=p_k \operatorname{ \;for \; some\;} j,k\in I.$$
Then the pointed set $$\bigwedge_{i\in I}S_i:=(\prod_{i\in I} S_i)/\sim_\wedge$$ with distinguished point $[p_\wedge:=(p_i)_{i\in I}]$ is called the \emph{smash product} of the family $S_i,i\in I$.
\dfe
\df
Let $N$ be a binoid, $(S_i,p_i)_{i\in I}$ be a finite family of pointed $N$-sets and $\sim_{\wedge_N}$ the equivalence relation on
$\bigwedge_{i\in I} S_i$ generated by
$$\cdots\wedge n+s_i\wedge \cdots \wedge s_j \wedge\cdots \sim_{\wedge_N}\cdots\wedge s_i\wedge\cdots\wedge n+s_j \wedge\cdots,$$ for all $i,j\in I$ and $n\in N$.
Then $$\bigwedge_{i\in I} {\!_N} S_i:=(\bigwedge_{i\in I} S_i)/\sim_{\wedge_N}$$ is called the
\textbf{smash product} of the family $(S_i)_{i\in I}$ over $N$.
\dfe
\ex
Let $n,m\in \n$ and $N=\n^\infty$, $S_1=\{0,1,\dots,n,p_1=\infty\}$, $S_2=\{0,1,\dots,m,p_2=\infty\}$, with operations given by: $$k+a=(k+a)\;\operatorname{mod} \;n,$$ $$k+b=(k+b)\;\operatorname{mod} \;m$$ for $k\in N^\bullet, a\in S_1,b\in S_2$. Then we try to compute $S_1\wedge_N S_2$. By definition we have $(a+k)\wedge b=a\wedge (b+k)$ for every $k\in \n, a\in S_1,b\in S_2$. By the Euclidean algorithm there exist $u,v\in \z$ such that $un+vm=\operatorname{gcd}(n,m)$ and so
\begin{equation}
0\wedge \operatorname{gcd}(n,m)=0\wedge (un+vm)=0\wedge un=un\wedge 0=0\wedge 0
  \end{equation}
 and also
 \begin{equation}
a\wedge b=(0+a)\wedge b=0\wedge (b+a).
 \end{equation}
  From $(1.1),(1.2)$ we have $S_1\wedge_N S_2\cong\{0,1,2,\dots,\operatorname{gcd}(n,m)-1,\infty\}$
\exe
\ex
Let $N=\n^\infty$, $S_1=\{0,1,2,3,4,5,6,7,p_1=\infty\},$ $S_2=\r$ with $p_2=0$ and operations given by: $$n+a=(n+a)\;\operatorname{mod} \;8$$ for $n\in N^\bullet, a\in S_1$, $$n+x=x^{2^n}$$ for $n\in N^\bullet, x\in S_2$. Then we try to compute $S_1\wedge_N S_2$. If we take the elements $a\in S_1,x\in S_2,k\in \z$ then we have
\begin{equation}
a\wedge x=a+8\wedge x=a\wedge x^{2^8}=a\wedge x^{2^{8k}}=0\wedge x^{2^{a+8k}}.
\end{equation}

If $x<0$ then $a\wedge x=a-1\wedge x^2=a\wedge |x|$. So the smash product for negative $x$'s is just another copy of positive $x$'s.
Now we look at the function $\operatorname{ln ln}: \r_{>1} \rightarrow \r$. Then we have
$$\operatorname{ln ln}(x^{2^{a+8k}})=\operatorname{ln}(2^{a+8k}\operatorname{ln}(x))=
(a+8k)\operatorname{ln}2+\operatorname{lnln}(x)$$ and so by $(1.3)$ we can conclude that $\{a\wedge x: a\in S_1, 1<x\in S_2\}$
is one circle and the action on $S_1$ is a $2\pi/8$ degree rotation in this circle. Also similarly for $1>x>0$ is also another circle.
\exe
\ex
\label{sm}
Let $N=\n^\infty$, $T=\{0,\infty\}$ the trivial binoid, $(S,p)$ some $N$-set with an $N$-action given by $f:S\rightarrow S,n+s=f^n(s)$. Then we can define two different smash products $T\wedge_N S$, depending on how $N$ acts on $T$. \\
(a) Suppose that the operation of $N$ on $T$ is given by $n+0=\infty$ for $n\geqslant 1$ and $n+\infty=\infty$, where $n\in \n$ and $0,\infty \in T$. Then $0\wedge f(s)=0\wedge (1+s)=1+0\wedge s=\infty \wedge s=\infty$ in $T\wedge_N S$. It means that $$T\wedge_N S\cong(S/ \im(f),p)$$
(b) Suppose that the operation of $N$ on $T$ is given by $n+0=0$ and $n+\infty=\infty$, for $n\in\n$ and $0,\infty \in T$. Then we have $0\wedge s=(1+0)\wedge s=0\wedge(1+s)=0\wedge f(s)$. Hence we have $$T\wedge_N S\cong(S/(f(s)\sim s),p).$$
\exe
\ex
Let $N=(\n^2)^\infty$, $T=\{0,\infty\}$, $S=\c$. The operation on $T$ is given by $(n_1,n_2)+0=0$, where $\infty\neq(n_1,n_2)\in N,0\in T$. This operation factors through the difference group of $N$. Let the action of $(\n^2)^\infty$ on $\c$ be given by $(1,0)+z=u+z,(0,1)+z=v+z$ where $u,v$ are non-collinear complex numbers. From here we can check that $T\wedge_N \c\cong \c/\Gamma=torus$, where $\Gamma$ is the lattice generated by $u,v$.
\exe
\prop
\label{pointed union and smash}
Let $N$ be a binoid and $S$, $T_i,1\leqslant i\leqslant k$ be $N$-sets. Then $$S\wedge_N (\bigcupdot_{i=1}^k T_i)=\bigcupdot_{i=1}^k (S\wedge_N T_i).$$
\prope
\pr
Let $s\wedge t\in S\wedge_N (\bigcupdot_{i=1}^k T_i)$, then there exists $j$ such that $t\in T_j$ or $t=p$ the point of $\bigcupdot_{i=1}^k T_i$. But in both cases $s\wedge t\in \bigcupdot_{i=1}^k (S\wedge_N T_i).$ Now take $p\neq s\wedge t\in \bigcupdot_{i=1}^k (S\wedge_N T_i)$, then there exists $j$ such that $s\wedge t\in S\wedge_N T_j$, so $s\wedge t\in S\wedge_N (\bigcupdot_{i=1}^k T_i)$.
\pre
\prop
\label{quotient}
Let $N$ be a commutative binoid and $(S,p)$ an $N$-set. If $I$ is an ideal of $N$ then $$(N/I)\wedge_N S\cong S/(I+S).$$
\prope
\pr
Take the homomorphism $$\phi:(N/I)\wedge_N S\rightarrow S/(I+S),\;\phi([n]\wedge s)=\overline{n+s},$$ where $[n]\in N/I,s\in S,\overline{n+s}\in S/(I+S)$. This is well defined because if $[n_1]\wedge s_1=[n_2]\wedge s_2$, then, by the universal property of smash product, we have $[0]\wedge n_1+s_1=[0]\wedge n_2+s_2$, so $n_1+s_1=n_2+s_2$.
Also from $\phi(0\wedge s)=\overline{s}$ we deduce that $\phi$ is surjective.
If $$\phi([n_1]\wedge s_1)=\phi([n_2]\wedge s_2)$$ for $[n_1],[n_2]\in N/I$ and $s_1,s_2\in S$, then we have $\overline{n_1+s_1}=\overline{n_2+s_2}$ in $S/(I+S)$. Hence we get $n_1+s_1=n_2+s_2$ or $n_1+s_1,n_2+s_2\in I+S$. If the first case holds then we have also $[n_1]\wedge s_1=[n_2]\wedge s_2$. If the second case holds then $n_1+s_1=f_1+t_1,n_2+s_2=f_2+t_2$, where $f_1,f_2\in I,t_1,t_2\in S$. So $$[n_i]\wedge s_i=[0]\wedge n_i+s_i=[0]\wedge f_i+t_i=[f_i]\wedge t_i=[\infty]\wedge t_i=p$$ for $i=1,2$. Hence again $[n_1]\wedge s_1=[n_2]\wedge s_2$ and we have injectivity which means $\phi$ is an isomorphism.
\pre
\prop
\label{surj map to smash}
Let $N$ be a commutative binoid and $J$ an ideal of $N$. If $S$ is a finitely generated $N$-set then there is a surjective homomorphism $(N/J)^{\cupdot r}\rightarrow S\wedge_N N/J$, where $r$ is the number of generators of $S$.
\prope
\pr
From Proposition \ref{quotient} we have an isomorphism $(N/J)^{\cupdot r}\rightarrow (N\wedge_N N/J)^{\cupdot r}$ and by Proposition \ref{pointed union and smash} we have $(N\wedge_N N/J)^{\cupdot r}=N^{\cupdot r}\wedge_N N/J$. Also by Proposition \ref{surj map for f.g. set} we have a surjective homomorphism $\phi:N^{\cupdot r}\rightarrow S$. Hence $$\varphi:(N)^{\cupdot r}\wedge_N N/J\longrightarrow S\wedge_N N/J,\; \varphi(a\wedge [n])=\phi(a)\wedge [n]$$ is also surjective homomorphism.
\pre

\cor
\label{quotient over primary ideal and smash}
Let $N$ be a commutative, finitely generated, semipositive binoid and $J$ an $N_+$-primary ideal of $N$. If $S$ is a finitely generated $N$-set then $S\wedge_N N/J$ is finite.
\core
\pr
We know that by Proposition \ref{primary ideal quotient}, $N/J$ is a finite set, so by Proposition \ref{surj map to smash} we have our result.
\pre

We close this section with three results about the smash product over the trivial base.
\prop
\label{inverse of smash products}
Let $N_1,\dots,N_k$ be non-zero binoids, then $$(\bigwedge_{i=1}^k N_i)^\times=\prod_{i=1}^k N_i^\times.$$
\prope
\pr
Let $a_1\wedge\cdots\wedge a_k\in (\bigwedge_{i=1}^k N_i)^\times$, then there exist $b_i\in N_i, 1\leqslant i \leqslant k$ such that $$(a_1\wedge\dots\wedge a_k)+(b_1\wedge\cdots\wedge b_k)=0\wedge \cdots\wedge 0.$$ So we can conclude that $a_i+b_i=0\in N_i$ and $a_i\in N_i^\times$ for $1\leqslant i \leqslant k$. Which means $(\bigwedge_{i=1}^k N_i)^\times\subseteq\prod_{i=1}^k N_i^\times.$ The other inclusion is clear.
\pre
\prop
\label{number of elements of smash}
Let $N,M$ be finite binoids, then $$\#(N\wedge M)=\#N\cdot\#M.$$
\prope
\pr
Let $a,b\in N^\bullet$ and $c,d\in M^\bullet$, then $a\wedge c=b\wedge d$ if and only if $a=b$ and $c=d$. From here we have the result.
\pre
\lem
\label{smash of quotients}
Let $M,N$ be binoids and $I\subseteq M,J\subseteq N$ be ideals. Then $(I\wedge N)\cup (M\wedge J)$ is an ideal of $M\wedge N$ and $$(M\wedge N)/((I\wedge N)\cup (M\wedge J))\cong M/I\wedge N/J.$$
\leme
\pr
We have a surjective homomorphism
\begin{align*}
M\wedge N &\xrightarrow{\;\phi\;} M/I\wedge N/J\\ m\wedge n &\longmapsto [m]\wedge [n].
\end{align*}
If $m\wedge n\in (I\wedge N)\cup (M\wedge J)$ then $m\wedge n\in I\wedge N$ or $m\wedge n\in M\wedge J$ and in each case we get $[m]\wedge [n]=\infty$. Hence $(I\wedge N)\cup (M\wedge J)\subseteq \ker \phi$. Now let $m_1\wedge n_1\in \ker \phi$ then $[m_1]\wedge [n_1]=\infty$ so $[m_1]=\infty$ or $[n_1]=\infty$. This means $m_1\in I$ or $n_1\in J$ and $\ker \phi\subseteq (I\wedge N)\cup (M\wedge J)$. Applying Lemma \ref{hom} we have a surjective homomorphism $$(M\wedge N)/((I\wedge N)\cup (M\wedge J))\longrightarrow M/I\wedge N/J$$ and it is also injective so
$$(M\wedge N)/((I\wedge N)\cup (M\wedge J))\cong M/I\wedge N/J.$$
\pre

\section{Exact sequences}
In this section arbitrary binoids are assumed to be commutative. 
\df
Let $N$ be a binoid. A sequence $$S_0 \xrightarrow{\;\phi_1} S_1 \xrightarrow{\;\phi_2} S_2 \xrightarrow{\;\phi_3}\cdots\xrightarrow{\;\phi_n} S_n$$ of $N$-sets and $N$-set homomorphisms is called \emph{exact} if the image of each homomorphism is equal to the kernel of the next:
$$\im\phi_k = \ker\phi_{k+1}.$$
\dfe
\df
Let $N$ be a binoid. The exact sequence $$S_0 \xrightarrow{\;\phi_1} S_1 \xrightarrow{\;\phi_2} S_2 \xrightarrow{\;\phi_3}\cdots\xrightarrow{\;\phi_n} S_n$$ is called \emph{strongly exact} if $\phi_k$ is injective on $S_{k-1}\setminus \ker \phi_k$ for every $k$.
\dfe
\prop
\label{exact sequence}
Let $N$ be a binoid, $S\subseteq T$ and $U$ be $N$-sets. Then we have an exact sequence of $N$-sets $$\infty\xrightarrow{\;\phi_1} \{s\wedge u\mid s\wedge u=\infty \operatorname{ in } T\wedge_N U\} \xrightarrow{\;\phi_2} S\wedge_N U\xrightarrow{\;\phi_3} T\wedge_N U\xrightarrow{\;\phi_4} (T/S)\wedge_N U\xrightarrow{\;\phi_5}\infty.$$
\prope
\pr
We have an exact sequence $\infty \rightarrow S\hookrightarrow T\rightarrow T/S\rightarrow \infty$ and we can smash this sequence with $U$. Then we obtain a sequence $$\infty\xrightarrow{\;\phi_1} \{s\wedge u\mid s\wedge u=\infty \operatorname{ in } T\wedge_N U\} \xrightarrow{\;\phi_2} S\wedge_N U\xrightarrow{\;\phi_3} T\wedge_N U\xrightarrow{\;\phi_4} (T/S)\wedge_N U\xrightarrow{\;\phi_5}\infty,$$ where $\phi_1(\infty)=\infty,$ $\phi_2$ is the inclusion, $\phi_3(s\wedge u)=s\wedge u\in T\wedge_N U$, $\phi_4(t\wedge u)=[t]\wedge u$ and $\phi_5([t]\wedge u)=\infty$.
We know by definition that $\im\phi_1 =\{\infty\}= \ker\phi_{2}$, $\im\phi_2 =\{s\wedge u\mid s\wedge u=\infty \operatorname{ in } T\wedge_N U\}= \ker\phi_{3}$, $\im\phi_3 =S\wedge_N U= \ker\phi_{4}$, $\im\phi_4 =T/S\wedge_N U= \ker\phi_{5}$. So it is an exact sequence.
\pre
\ex
\label{example of not strongly exact}
Let $N=\n^\infty$ and $S$ be an $N$-set with an operation given by $n+s=f^n(s)$, like in Example \ref{n set}. Then we have an exact sequence $$\infty\longrightarrow N_+\longrightarrow N\longrightarrow N/N_+\longrightarrow\infty$$ and we can smash this sequence with $S$ over $N$.
Then by Proposition \ref{exact sequence}, we get an exact sequence $$\infty\rightarrow \{n\wedge s\mid n\geqslant 1, n\wedge s=\infty \;\operatorname{in}\; N\wedge_N S\}\rightarrow N_+\wedge_N S\rightarrow N\wedge_N S\rightarrow N/N_+\wedge_N S\rightarrow\infty.$$
We have also isomorphisms $$N\wedge_N S\longrightarrow S,\;n\wedge s\longmapsto n+s=f^n(s),$$ $$N_+\wedge_N S\longrightarrow S,\;n\wedge s\longmapsto (n-1)+s=f^{n-1}(s),$$ and from here we get $\{n\wedge s\mid n\geqslant 1, n\wedge s=\infty \;\operatorname{in}\; N\wedge_N S\}\cong\{t\in S\mid f(t)=p\}=\ker f$ and $N/N_+\wedge_N S\cong S/\im f$(as was also shown in Example \ref{sm}(a)). So we have an exact sequence $$\infty\longrightarrow \ker f\longrightarrow S\xrightarrow{\;f\;} S \longrightarrow S/\im f\longrightarrow\infty.$$ If $S=\{a,b,p\}$, $f(a)=f(b)=b$ and $f(p)=p$ then this sequence is not strongly exact.
\exe
\prop
\label{strongly exact sequence}
Let $N$ be a binoid, $J\subseteq N$ an ideal and $S\subseteq T$  $N$-sets. Then we have a strongly exact sequence of $N$-sets $$\infty\xrightarrow{\phi_1} \{s\wedge [a]\mid s\wedge [a]=\infty \operatorname{ in } T\wedge_N N/J\} \xrightarrow{\phi_2} S\wedge_N N/J\xrightarrow{\phi_3} T\wedge_N N/J\xrightarrow{\phi_4} T/S\wedge_N N/J\xrightarrow{\phi_5}\infty$$ which is the same as $$\infty\xrightarrow{\phi_1}S\cap(J+T))/(J+S) \xrightarrow{\phi_2}
S/(J+S)\xrightarrow{\phi_3}T/(J+T)\xrightarrow{\phi_4}(T/S)/(J+T/S))\xrightarrow{\phi_5}\infty.$$
If $S=I$ is an ideal of $N$ and $T=N$, then we have the strongly exact sequence
$$\infty\longrightarrow (I\cap J)/(I+J) \longrightarrow I/(I+J) \longrightarrow N/J\longrightarrow (N/I)/(J+N/I)\longrightarrow\infty.$$
\prope
\pr
From Proposition \ref{exact sequence}, when $U=N/J$, we have an exact sequence $$\infty\xrightarrow{\phi_1} \{s\wedge [a]\mid s\wedge [a]=\infty \operatorname{ in } T\wedge_N N/J\} \xrightarrow{\phi_2} S\wedge_N N/J\xrightarrow{\phi_3} T\wedge_N N/J\xrightarrow{\phi_4} T/S\wedge_N N/J\xrightarrow{\phi_5}\infty.$$
By Proposition \ref{quotient} we know that $S\wedge_N N/J\cong S/(J+S)$ and $T\wedge_N N/J\cong T/(J+T)$.
Let $$s\wedge [a]\in\{s\wedge [a]\mid s\wedge [a]=\infty \operatorname{ in } T\wedge_N N/J\}.$$ Then $S/(J+S)\ni[a+s]=\infty\in T/(J+T)$. So we have $a+s\in J+T$, which means that $a+s\in S\cap (J+T)$. Hence we have $$\{s\wedge [a]\mid s\wedge [a]=\infty \operatorname{ in } T\wedge_N N/J\}\cong(S\cap(J+T))/(J+S).$$ Also, by Proposition \ref{quotient} and Proposition \ref{quotient to union}, we get $$T/S\wedge_N N/J\cong (T/S)/(J+T/S)=T/(S\cup(J+T)).$$ So we can rewrite the previous exact sequence as
$$\infty\xrightarrow{\phi_1}(S\cap(J+T))/(J+S) \xrightarrow{
\phi_2}S/(J+S)\xrightarrow{\phi_3}T/(J+T)\xrightarrow{\phi_4}(T/S)/(J+T/S))\xrightarrow{\phi_5}\infty$$ or, equivalently
$$\infty\xrightarrow{\phi_1}(S\cap(J+T))/(J+S) \xrightarrow{
\phi_2}S/(J+S)\xrightarrow{\phi_3}T/(J+T)\xrightarrow{\phi_4}T/(S\cup(J+T))\xrightarrow{\phi_5}\infty.$$
Here $\phi_1(\infty)=\infty$, $\phi_2$ is the inclusion (so it is injective), $\phi_3$ is an inclusion (injective) on $S\setminus (J+T)$, $\phi_4$ is surjective and outside of the kernel it is a bijection, $\phi_5([s])=\infty$. So it means that our sequence is a strongly exact sequence.

If $S=I$ and $T=N$ then we get
\[\infty\xrightarrow{\;\phi_1}I\cap J/(J+I) \xrightarrow{\;\phi_2}I/(J+I)\xrightarrow{\;\phi_3}N/J\xrightarrow{\;\phi_4}(N/I)/(J+N/I)
\xrightarrow{\;\phi_5}\infty. \qedhere \]
\pre

\prop
\label{general equation of exact seq}
Let $N$ be a binoid and $\infty \rightarrow S_1 \rightarrow S_2 \rightarrow\cdots\rightarrow S_n\rightarrow \infty$ a strongly exact sequence of finite $N$-sets. Then $$\sum_{i=1}^n (-1)^i\# S_i=0.$$
\prope
\pr
Write $S_i=K_i\uplus R_i\uplus \{p_i\}$, with maps
$$\begin{aligned}
  \phi_i:S_{i-1}&\longrightarrow S_i,\\
      R_{i-1}&\xrightarrow{bij} K_i,\\
      K_{i-1}&\longrightarrow p_i,\\
      p_{i-1}&\longrightarrow p_i,
  \end{aligned}
$$
where $1 \leqslant i \leqslant n+1,$ $$S_0=S_{n+1}=\{\infty\},\;p_0=p_{n+1}=\infty,$$ and $$R_0=K_0=K_1=K_{n+1}=R_{n+1}=\varnothing.$$ Then we can conclude that $$\sum_{i=1}^n (-1)^i\# S_i=\sum_{i=1}^n (-1)^i(|K_i|+|R_i|)=\sum_{i=1}^n (-1)^i(|K_i|+|K_{i+1}|)=-| K_1 |+(-1)^n |K_{n+1} |=0.$$
\pre
\cor
\label{equality of e.s. corollary}
Let $N$ be a commutative, finitely generated, semipositive binoid and let $I$ be an ideal of $N$. If $J$ is an $N_+$-primary ideal of $N$ then $$\# N/J+\# I\cap J/(I+J)=\# I/(I+J)+\# (N/I)/(J+N/I).$$
\core
\pr
We want to apply Proposition \ref{general equation of exact seq} to the strongly exact sequence $$\infty\longrightarrow I\cap J/(I+J) \longrightarrow I/(I+J) \longrightarrow N/J\longrightarrow (N/I)/(J+N/I)\longrightarrow\infty$$ from Proposition \ref{strongly exact sequence}. To do this we have to show that the involved $N$-sets are finite. From Proposition \ref{primary ideal quotient} we know that $N/J$ is a finite set.
Also we know that $I$ is a finitely generated $N$-set, so by Proposition \ref{surj map to smash} we have a surjective homomorphism $(N/J)^{\cupdot r}\rightarrow I\wedge_N N/J$. Hence $|I\wedge_N N/J |=|I/(J+I)|\leqslant |N/J|^r$, which is a finite set. So we can apply Proposition \ref{general equation of exact seq} and get the result.
\pre

\section{Dimension}

\df
The \emph{combinatorial dimension} of a binoid $N$, denoted by $\dim N$, is the supremum of the lengths of strictly increasing chains of prime ideals of $N$.
\dfe
We will assume that the zero binoid $\{\infty\}$ has dimension $-1$. The trivial binoid $\{0,\infty\}$ has dimension 0.

In general the Krull dimension of a binoid algebra $K[N]$ and $\dim N$ are not the same.
\ex
Let $N=(\n^2)^{\infty}/(x+y=2y)$. It is separated and we have $\dim_{\krull} K[N]=1$, as $K[N]=K[X,Y]/Y(Y-X)$. But we have the following chain of prime ideals $(\infty)\subset(y)\subset(x,y)$ so $\dim N=2$.
\exe
\prop
\label{dimension decrease}
Let $N$ be a finitely generated, integral binoid and $I\neq \{\infty\}$ be an ideal of $N$. Then $\dim N/I< \dim N$.
\prope
\pr
By \cite[Corollary 2.2.10]{Simone} we know that if $k=\dim N/I$ then there exist $\mathfrak{p}_0,\dots,\mathfrak{p}_k\in \spec N$ such that $\infty\neq I\subseteq \mathfrak{p}_0\subset\cdots\subset\mathfrak{p}_k$. By integrality, $\{\infty\}$ is a prime ideal of $N$ and so $\dim N> \dim N/I$.
\pre
\prop
\label{dimension decreasing of union of minimal primes}
Let $N$ be a finitely generated commutative binoid and $\mathfrak{p},\mathfrak{q}$ different minimal prime ideals of $N$. Then $\dim N/(\mathfrak{p}\cup \mathfrak{q})<\min\{\dim N/\mathfrak{p},\dim N/\mathfrak{q}\}$.
\prope
\pr
By minimality $\mathfrak{p}\nsubseteq \mathfrak{q}$ and $\mathfrak{q}\nsubseteq \mathfrak{p}$ and it is easy to check that $\mathfrak{p}\cup \mathfrak{q}$ is also a prime ideal. Hence $\mathfrak{p},\mathfrak{q}\subset \mathfrak{p}\cup \mathfrak{q}$, which means $\dim N/(\mathfrak{p}\cup \mathfrak{q})<\dim N/\mathfrak{p}$ and $\dim N/(\mathfrak{p}\cup \mathfrak{q})<\dim N/\mathfrak{q}$.
\pre
\tr
\label{dimension of smash product}
Let $M,N$ be non zero commutative binoids of finite dimension. Then $$\dim M \wedge N=\dim M+ \dim N.$$
\tre
\pr
Let $\mathfrak{p}\subset M,\mathfrak{q}\subset N$ be prime ideals. Then we know that there are binoid homomorphisms $\chi_\mathfrak{p}:M\rightarrow \{0,\infty\}$ and $\chi_\mathfrak{q}:N\rightarrow \{0,\infty\}$. Then by the universal property of $M\wedge N$ we have a binoid homomorphism
\begin{align*}
 \chi_\mathfrak{p}+\chi_\mathfrak{q}:M\wedge N &\longrightarrow\{0,\infty\}\\
 m\wedge n&\longmapsto \chi_\mathfrak{p}(m)+\chi_\mathfrak{q}(n),
\end{align*}
and it will give us the prime ideal $\ker(\chi_\mathfrak{p}+\chi_\mathfrak{q})=\mathfrak{p}\wedge N\cup M\wedge\mathfrak{q}$ of $M\wedge N$, and every prime ideal of $M\wedge N$ is of this form.
\begin{claim}
 Let $\mathfrak{p},\mathfrak{p'}\in\spec M,\mathfrak{q},\mathfrak{q'}\in\spec N$. Then $\mathfrak{p}\wedge N\cup M\wedge\mathfrak{q}\subseteq \mathfrak{p'}\wedge N\cup M\wedge\mathfrak{q'}$ if and only if $\mathfrak{p}\subseteq \mathfrak{p'}$ and $\mathfrak{q}\subseteq \mathfrak{q'}$.
\end{claim}
\pr
The if part is obvious. So let us assume $\mathfrak{p}\wedge N\cup M\wedge\mathfrak{q}\subseteq \mathfrak{p'}\wedge N\cup M\wedge\mathfrak{q'}$ and $\mathfrak{p}\nsubseteq \mathfrak{p'}$ or $\mathfrak{q}\nsubseteq \mathfrak{q'}$. Since $\mathfrak{p}\nsubseteq \mathfrak{p'}$, there exists $f\in\mathfrak{p}$ such that $f\notin \mathfrak{p'}$ and $(\chi_\mathfrak{p}+\chi_\mathfrak{q})(f\wedge 0)=\infty+0=\infty$ but $(\chi_\mathfrak{p'}+\chi_\mathfrak{q'})(f\wedge 0)=0+0=0$. So this contradicts the inclusion $\mathfrak{p}\wedge N\cup M\wedge\mathfrak{q}\subseteq \mathfrak{p'}\wedge N\cup M\wedge\mathfrak{q'}$.
\pre
Now let $\dim M=m$ and $\dim N=n$. Then we have strictly increasing chains of prime ideals $$\mathfrak{p}_0\subset\mathfrak{p}_1\subset\cdots\subset \mathfrak{p}_m\subseteq M$$ and $$\mathfrak{q}_0\subset\mathfrak{q}_1\subset\cdots\subset \mathfrak{q}_n\subseteq N.$$
Hence from Claim 1, we get a strictly increasing chain of prime ideals of $M\wedge N$
$$(\mathfrak{p}_0\wedge N\cup M\wedge\mathfrak{q}_0)\subset\cdots\subset(\mathfrak{p}_0\wedge N\cup M\wedge\mathfrak{q}_n)\subset(\mathfrak{p}_1\wedge N\cup M\wedge\mathfrak{q}_n)\subset\cdots\subset(\mathfrak{p}_m\wedge N\cup M\wedge\mathfrak{q}_n),$$ which means $\dim M\wedge N\geqslant m+n.$

Now if we have a strictly increasing chain $$(M\wedge  \mathfrak{p}_0\cup \mathfrak{q}_0\wedge N)\subset \cdots\subset (M\wedge  \mathfrak{p}_k\cup \mathfrak{q}_k\wedge N)\in \spec M\wedge N,$$ where $k>m+n$. Then we have increasing chains $\mathfrak{p}_0\subseteq\cdots\subseteq\mathfrak{p}_k$ and $\mathfrak{q}_0\subseteq\cdots\subseteq\mathfrak{q}_k$ in $M$ and $N$ respectively. If $\mathfrak{p}_i=\mathfrak{p}_{i+1}$ then $\mathfrak{q}_i\subset\mathfrak{q}_{i+1}$. So, if only $m+1$ different prime ideals of $M$ occur in the sequence then $$\mathfrak{p}_0=\dots=\mathfrak{p}_{i_1},\mathfrak{p}_{i_1+1}=\dots=\mathfrak{p}_{i_2},\dots,\mathfrak{p}_{i_m+1}=\dots=\mathfrak{p}_{i_{m+1}},$$ where $i_{m+1}=k$. So we have $$\mathfrak{q}_0\subset\dots\subset\mathfrak{q}_{i_1},\mathfrak{q}_{i_1+1}\subset\dots\subset\mathfrak{q}_{i_2},\dots,\mathfrak{q}_{i_m+1}\subset\dots\subset\mathfrak{q}_{i_{m+1}}$$ and here we have at least $n+2$ different prime ideals of $N$. Hence it is a contradiction.
\pre

\section{Algebras and Modules}
In this section we will assume that $K$ is a commutative ring and arbitrary binoids are assumed to be commutative.

\df
The \emph{binoid algebra} of a binoid $N$ is the quotient algebra $$KN/\langle X^\infty\rangle=:K[N],$$
where $KN$ is the monoid algebra of $N$ and $\langle X^\infty\rangle$ is the ideal in $KN$ generated by the element $X^\infty$.
\dfe
So we can consider $K[N]$ as the set of all formal sums $\sum_{m\in M}r_m X^m,$ where $M\subseteq N^\bullet$ is finite, $r_m\in K$ and the multiplication is given by

$$r_nX^n\cdot r_m X^m=\left\{
    \begin{array}{ccc}
     (r_nr_m)X^{n+m},&\operatorname {\;if\;} n+m\in N^\bullet \\
     0,&\;\operatorname {if\;} n+m=\infty. \\
      \end{array}
       \right.$$

Let $N$ be a binoid and $(S,p)$ an $N$-set. Then we can define the $K[N]$-module $K[S]$ as the set of all formal sums $\sum_{s\in U}r_s X^s,$ where $U\subseteq S^\bullet$ is finite, $r_s\in K$ and the multiplication is given by
$$(\sum_{m\in M}r_m X^m)\cdot (\sum_{s\in U}r_s X^s)=\left\{
    \begin{array}{ccc}
     \sum_{m\in M,s\in U}r_mr_s X^{m+s},&\operatorname {\;if\;} m+s\in S^\bullet \\
     0,&\operatorname {if\;} m+s=p. \\
      \end{array}
       \right.$$
 Here $M$ is a finite subset of $N$ and $U$ is a finite subset of $S$.
\prop
Let $N$ be a binoid and $I$ an ideal of $N$. Then $$K[I]:=\{\sum_{a\in J}r_a X^a\mid  J\subseteq I \operatorname{finite\;subset\;}\}$$ is an ideal of $K[N]$.
\prope
\pr
This follows from $$(\sum_{m\in M}r_m X^m)\cdot (\sum_{a\in J}r_a X^a)=\sum_{m\in M,a\in J}r_mr_a X^{m+a}\in K[I].$$
\pre
\prop
\label{K algebra}
Let $N$ be a binoid. Then the following hold.
\begin{enumerate}
 \item If $M$ is a subbinoid of $N$ then the corresponding binoid algebra $K[M]$ is a $K$-subalgebra of $K[N]$.
 \item If $\mathfrak{a}$ is an ideal of $K$ then $(K/\mathfrak{a})[N]\cong K[N]/\mathfrak{a}K[N]$.
 \item If $I$ is an ideal of $N$ then $K[N/I]\cong K[N]/K[I]\cong KN/KI$.
 \item If $S$ is a subbinoid of $M$ then $\{X^a\mid a\in S\}$ defines a multiplicative system in $K[N]$, which will also be denoted by $S$, and there is an isomorphism $K[N]_S:=S^{-1}K[N]\cong K[N_S]$.
 \item If $B$ is a $K$-algebra then $B\otimes_K K[N]\cong B[N]$.
 \item Is $S$ is a multiplicative system in $K$ then $K[N]_S\cong K_S[N]$.
\end{enumerate}
\prope
\pr
See the proof of Corollary 3.2.8 in \cite{Simone}.
\pre
\prop
\label{finite set and dimension}
Let $N$ be a binoid and $S$ a finite $N$-set. Then $\# S=\dim_K K[S]$.
\prope
\pr
Let $S=\{p,s_1,\dots,s_m\}$ then by definition of $K[S]$, \{$X^{s_1}$, \dots,$X^{s_m}$\} is a $K$-basis of $K[S]$.
\pre
An $N$-set homomorphism $\phi:S\rightarrow T$ induces a $K[N]$-module homomorphism $K[\phi]:K[S]\rightarrow K[T]$.
\prop
Let $N$ be a binoid and $S,T$ be $N$-sets. Then the following are equivalent.
\begin{enumerate}
 \item $\phi:S\rightarrow T$ is a surjective (injective) homomorphism,
 \item $K[\phi]:K[S]\rightarrow K[T]$ is a surjective (injective) homomorphism for all commutative rings,
 \item $K[\phi]:K[S]\rightarrow K[T]$ is a surjective (injective) homomorphism for some commutative ring $K\neq 0$.
\end{enumerate}
\prope
\pr
(1) to (2). Let $\phi:S\rightarrow T$ be a surjective homomorphism of $N$-sets and $\sum_{t\in J} r_t X^t\in K[T]$, where $J$ is some finite subset of $T$. Then we can write $\sum_{t\in J} r_t X^t=\sum_{\phi(s)\in J} r_{\phi(s)} X^{\phi(s)}$, so $K[\phi](\sum_{\phi(s)\in J} r_{\phi(s)} X^s)=\sum_{t\in J} r_t X^t$.

(2) to (3) is obvious.

(3) to (1). Let $K$ be some nonzero commutative ring and suppose that $K[\phi]:K[S]\rightarrow K[T]$ is a  surjective ring homomorphism. So for every $t\in T$ we know that there is $\sum_s r_s X^s\in K[S]$ such that $X^t=K[\phi](\sum_s r_s X^s)=\sum_s r_s X^{\phi(s)}$. Hence $t=\phi(s)$ for some $s\in S$.
\pre

\prop
\label{exact sequence of K algebra}
Let $N$ be a binoid, $$\infty \longrightarrow S_1 \longrightarrow S_2 \longrightarrow\cdots\longrightarrow S_n\longrightarrow \infty$$ a strongly exact sequence of finite $N$-sets and $K$ a commutative ring. Then we have an exact sequence of $K[N]$-modules $$0 \longrightarrow K[S_1] \longrightarrow K[S_2] \longrightarrow\cdots\longrightarrow K[S_n]\longrightarrow 0.$$
\prope
\pr
For every $i$, let $\phi_i:S_i\rightarrow S_{i+1}$ be the $N$-set homomorphism of the sequence. Then we can define homomorphisms $K[\phi_i]: K[S_i]\rightarrow K[S_{i+1}]$ by setting $K[\phi_i](\sum_{s\in U}r_s X^s)=\sum_{s\in U}r_s X^{\phi_i(s)}$.
We have a strongly exact sequence of $N$-sets, so $$(K[\phi_{i+1}]\circ K[\phi_i])=K[\phi_{i+1}\circ\phi_i]=K[p_{i+1}]=0,$$ which means $\im K[\phi_i]\subseteq\ker K[\phi_{i+1}]$. We have the other inclusion because $\phi_i$ is an injection outside the kernel, which implies that $K[\phi_i]:K[S_i/\ker \phi_i]\rightarrow K[S_{i+1}]$ is injective.
\pre
\ex
Let $S=\{a,b,p\}$ be as in Example \ref{example of not strongly exact}. Then we have an exact sequence of $N$-sets $$\infty\longrightarrow\ker f=\infty\xrightarrow{\;i\;\,} S\xrightarrow{\;f\;} S\longrightarrow S/\im f\longrightarrow\infty.$$ We have $$K[f](X^a-X^b)=X^{f(a)}-X^{f(b)}=X^b-X^b=0,$$ but $X^a-X^b\notin \im K[i]=\{\infty\}$. So strong exactness is a necessary condition for Proposition \ref{exact sequence of K algebra}.
\exe

\prop
\label{smash to tensor}
Let $N$ be a binoid and $S,T$ the $N$-sets. Then we have $$K[S\wedge_N T]\cong K[S]\otimes_{K[N]} K[T]$$ and $$K[S\cupdot T]\cong K[S]\oplus K[T].$$
\prope
\pr
The first isomorphism follows from Corollary 3.5.2 in \cite{Simone}.

From the inclusions $S\hookrightarrow S\cupdot T$ and $T\hookrightarrow S\cupdot T$ we get $K[N]$-homomorphisms $K[S]\rightarrow K[S\cupdot T]$ and $K[T]\rightarrow K[S\cupdot T]$. So we have a $K[N]$-homomorphism $$\phi:K[S]\oplus K[T]\rightarrow K[S\cupdot T].$$
Also from the homomorphisms $S\cupdot T \xrightarrow{T\rightarrow p} S$, $S\cupdot T \xrightarrow{S\rightarrow p} T$ we get the module homomorphisms $K[S\cupdot T]\rightarrow K[S]$ and $K[S\cupdot T]\rightarrow K[T]$. Hence we have a homomorphism $$\varphi:K[S\cupdot T]\rightarrow K[S]\oplus K[T].$$
It is easy to check that $\phi\circ \varphi$ and $\varphi\circ \phi$ are the identity maps.
\pre

\prop
\label{quotient for tensor}
Let $R$ be a commutative ring and $V$ be an $R$-module. If $I$ is an ideal of $R$ then $(R/I)\otimes_R V\cong V/IV$. Here $IV$ is the submodule generated by all products $fv$, where $f\in I$ and $v\in V$.
\prope
\pr
It is easy to check that $V/IV$ is an $R$-module. If $f\in I$ then for every $v\in V$ we have $fv=0\in V/IV$. So $I\subseteq \ann_R(V/IV)$ and $V/IV$ is an $R/I$-module.
So we have an $R/I$-module homomorphism
$$\begin{aligned}
  (R/I)\otimes_R V&\longrightarrow V/IV,\\
  \sum_i a_i\otimes v_i&\longmapsto [\sum_i a_iv_i].
  \end{aligned}
$$
This map is a surjective homomorphism by $1\otimes v\mapsto [v]$. If $\sum a_i\otimes v_i\mapsto [\sum a_iv_i]=0\in V/IV$ then $\sum_i a_iv_i=\sum_j b_jw_j$, where $b_j\in I$ and $w_j \in V$. Hence $$1\otimes\sum_i a_iv_i=1\otimes \sum_j b_jw_j=\sum_j b_j\otimes w_j=\sum_j 0\otimes w_j=0.$$ So it is injective.
\pre
\prop
\label{K algebra of quotient}
Let $N$ be a binoid, $S$ be an $N$-set and $I$ be an ideal of $N$. Then we have $$K[S/(I+S)]\cong K[S]/(K[I]K[S]).$$
\prope
\pr
By Proposition \ref{quotient} we have $S/(I+S)\cong (N/I)\wedge_N S$, so $K[S/(I+S)]\cong K[(N/I)\wedge_N S]$. Hence from Proposition \ref{K algebra} (3), Proposition \ref{smash to tensor} and Proposition \ref{quotient for tensor} we get
\[
\begin{aligned}
K[S/(I+S)]&\cong K[(N/I)\wedge_N S] \\
&\cong K[N/I]\otimes_{K[N]} K[S]\\
&\cong (K[N]/K[I])\otimes_{K[N]} K[S] \\
&\cong K[S]/(K[I]K[S]).
\end{aligned}\qedhere
\]
\pre
\lem
\label{semipositive binoid algebra}
Let $N$ be a semipositive binoid and $K$ some field with characteristic $p$ which does not divide $|N^\times|$. Then $K[N_+]$ is the intersection of finitely many maximal ideals.

In particular, if $N$ is a positive binoid then $K[N_+]$ is a maximal ideal of $K[N]$.
\leme
\pr
Because of the fundamental theorem for finite abelian groups we can write $$N^\times=\z/(\alpha_1)\times\cdots\times \z/(\alpha_r).$$ By assumption $p$ does not divide $\alpha_1,\dots,\alpha_r$. So from these conditions we can deduce that
$$\begin{aligned}
\overline{K}N^\times &\cong \overline{K}[X_1,\dots,X_r]/(X_1^{\alpha_1}-1,\dots,X_r^{\alpha_r}-1)\\
&=\overline{K}[X_1,\dots,X_r]/\big((X_1-\xi_{11})\cdots(X_1-\xi_{1\alpha_1}),\dots,
(X_r-\xi_{r1})\cdots(X_r-\xi_{r\alpha_r})\big)\\
&\cong\overline{K}[X_1]/\big((X_1-\xi_{11})\cdots(X_1-\xi_{1\alpha_1})\big)\otimes\cdots
\otimes\overline{K}[X_r]/\big((X_r-\xi_{r1})\cdots(X_r-\xi_{r\alpha_r})\big)\\
&\cong\overline{K}^{|N^\times|},
  \end{aligned}$$
where $\overline{K}$ is the algebraic closure of $K$ and $\xi_{ij}$ are the $\alpha_i$-th roots of unity.
Hence the maximal ideals of $\overline{K}N^\times$ have the form $$\mathfrak{m}_{i_1,\dots,i_r}=(X_1-\xi_{1i_1},\dots,X_r-\xi_{ri_r}).$$ So we have finitely many maximal ideals in $\overline{K} N^\times$ with this form. 
We also know that the intersection of all maximal ideals of $\overline{K} N^\times$ is equal to $\nil (\overline{K} N^\times)$ and this is $0$. 
Under the homomorphism $KN^\times\hookrightarrow\overline{K}N^\times$ the preimage of a maximal ideal is maximal and therefore the intersection of the maximal ideals of $KN^\times$ is 0 as well.

Let $K[\pi]:K[N]\rightarrow K[N/N_+]$ be the homomorphism induced by $\pi:N\rightarrow N/N_+\cong (N^\times)^\infty$. Then 
$K[\pi]^{-1}(\mathfrak{m_i})$ is a maximal ideal of $K[N]$, where $\mathfrak{m_i}$ is a maximal ideal of $K[N^\times]$. So 
\[ K[N_+]=K[\pi]^{-1}(0)=\bigcap_i K[\pi]^{-1}(\mathfrak{m_i}).\qedhere\]
\pre
\ex
Let $p$ be a prime number. Then
$$(\z/p\z)[X]/(X^p-1)=(\z/p\z)[X]/(X-1)^p$$ is not reduced and has only one maximal ideal $(X-1)$ so it is a local ring.
\exe

\chapter{Hilbert-Kunz theory of binoids}

In this chapter we will introduce the Hilbert-Kunz function and the Hilbert-Kunz multiplicity of a binoid and compute it in some cases. Also we will prove that the Hilbert-Kunz multiplicity of a commutative, finitely generated, semipositive, cancellative and reduced binoid is a rational number.

\section{Hilbert-Kunz function and Hilbert-Kunz multiplicity of a ring}

Let $(R,\mathfrak{m})$ be a commutative Noetherian local ring of dimension $d$ and containing a field $K$ of positive characteristic $p$. For an ideal $I$ and a prime power $q=p^e$ we define the ideal $I^{[q]}=\langle a^q\mid a\in I\rangle$ as the ideal generated by the $q$th powers of the elements of $I$.
Let $I$ be an $\mathfrak{m}$-primary ideal of $R$, $M$ a finite $R$-module. Then the $R$-modules $M/I^{[q]}M$ has finite length. So the \emph{Hilbert-Kunz function} of $M$ with respect to $I$ is $$\hkf(I,M)(q)=\len(M/I^{[q]}M).$$ If $M=R,I=\mathfrak{m}$ then we have the classical Hilbert-Kunz function $\hkf(\mathfrak{m},R)(q)=HK_R(q)$, introduced by Kunz \cite{Kunz}. Also he showed that $R$ is regular if and only if $HK_R(q)=q^d$ for all $q$.
In \cite{Monsky}, P. Monsky proved the following theorem.
\tr[Monsky]
There are a real constant $c(M)$ and a natural number $d$ (the dimension of $M$) such that $$\len(M/I^{[q]}M)=c(M)q^d+O(q^{d-1}).$$
\tre
So we have the next definition.
\df
The \emph{Hilbert-Kunz multiplicity} $e_{HK}(I,M)$ of $M$ with respect to $I$ is defined as:
$$e_{HK}(I,M):=\lim_{q\rightarrow\infty}\dfrac{\len(M/I^{[q]}M)}{q^d}.$$
\dfe

From here there are many questions related to Hilbert-Kunz function and multiplicity.
\prob
 Is the Hilbert-Kunz multiplicity always a rational number?
\probe
\prob
 Is there any interpretation in characteristic 0?
\probe
To formulate the next three problems we fix a relative situation. Let $R$ be a $\z$-algebra of finite type with a ring homomorphism $R\xrightarrow{\phi} \z$. Let $R_p$ be the localization of $R\otimes_\z \z/p\z$ at the maximal ideal given by $\phi$. Let $I\subseteq R$ be an ideal such that the extension $I_p$ in $R_p$ is primary to the maximal ideal for almost all $p$. In this situation we can look at $e_{HK}^{R_p}(I_p)$ and its dependence on the prime number $p$. In particular we have the following problems.
\prob
When and how does $e_{HK}(I_p,R_p)$ depend on the characteristic $p$?
\probe
\prob
Does the limit $$\lim_{p\rightarrow\infty} e_{HK}(I_p,R_p)$$ exist?
\probe
\prob[C.Miller]
 Does the limit $$\lim_{p\rightarrow \infty} \dfrac{\len(R_p/I_p^{[p]})}{p^d}$$ exist?
\probe
In most known cases the Hilbert-Kunz multiplicity is a rational number, for example for monoid rings (\cite{Eto}), toric rings (\cite{Watanabe},\cite{Bruns}), monomial ideals and binomial hypersurfaces (\cite{Conca}), rings of finite Cohen-Macaulay type (\cite{Seibert}), two-dimensional graded rings (\cite{Bre06},\cite{Trivedi2}) and Fermat ring $K[X_1,\dots,X_n]/(X_1^{d_1}+\cdots+X_n^{d_n})$ (\cite{MonHan}).

Let $G$ be a finite group acting linearly on a polynomial ring $B=K[x_1,\ldots ,x_n]$ with invariant ring $A=B^G$ and let $I = A_+ B$ be the Hilbert ideal in $B$. Let $R$ be the localization of $A$ at the irrelevant ideal. If $K$ has positive characteristic, then
$$e_{HK}(R) = \frac{ \operatorname{dim}_K \, (B/I ) }{\#( G )}.$$ This is a rational number which depends only on the invariant ring, not on its representation.
But H. Brenner proved in \cite{Brenner} that there exist three-dimensional quartic hypersurface domains and Artinian modules with irrational Hilbert-Kunz multiplicity.

There are many situations where the Hilbert-Kunz multiplicity is independent on the characteristic $p$. For example the Hilbert-Kunz multiplicity of toric rings (\cite{Watanabe}) and of invariant rings are independent from the characteristic of the base field as long as $p$ does not divide the group order. But there are also examples where the Hilbert-Kunz multiplicity is depending on the characteristic. Let $$R=\z[X,Y,Z]/(X^4 +Y^4+Z^4)$$ be the Fermat-quartic and $R_p=R\otimes \z/p\z$. The Hilbert-Kunz multiplicity\footnote{By the Hilbert-Kunz multiplicity of a standard-graded ring we mean the Hilbert-Kunz multiplicity of the localization at the irrelevant ideal.} is
$$e_{HK}(R_p)=\begin{cases} 3 \text{ for } p = \pm 1 \,\operatorname{mod} \, 8, \\ 3 + \frac{1}{p^2}\text{ for } p= \pm 3 \, \operatorname{mod} \, 8 \, .\end{cases}$$
In this example we see also that $\lim_{p\rightarrow\infty}e_{HK}(R_p)=3$.

We can ask when the limit of the Hilbert-Kunz multiplicity exists for $p\rightarrow\infty$. If it exists then this limit is a candidate for the Hilbert-Kunz multiplicity in characteristic zero. This  leads us to the question of whether a characteristic zero Hilbert-Kunz multiplicity could be defined directly. In all known cases this limit exists.

H. Brenner, J. Li and C. Miller (\cite{BreLiMil}) observed that in all the known cases where $\lim_{p\rightarrow\infty}e_{HK}(R_p)$ exists then this double limit can be replaced by the limit of Problem 5.

The Hilbert-Kunz function has also another description.
Let $R$ be a commutative ring, $V$ be an $R$-module and $\phi:R\rightarrow R$ a ring homomorphism. Then $V$ is also an $R$-module by $r\cdot a:=\phi(r)a$ and we denote it by $^\phi V$.
So if $\phi$ is an $e$th iteration of the Frobenius homomorphism $F:R\rightarrow R$, $f\mapsto f^p$, then we have an $R$ module $^eR$ by $r\cdot a:=F^{e*}(r)a$.
\prop
Let $R$ be a local ring with maximal ideal $\mathfrak{m}$, $I$ an $\mathfrak{m}-$primary ideal of $R$ and $q=p^e$. Then we have $R/I^{[q]}=$ $^eR\otimes_R R/I$.
\prope
\pr
If we take $V=$ $^eR$ then by Proposition \ref{quotient for tensor} $^eR\otimes_R R/I=$ $^eR/I$ $ ^eR$. So $I$ $^eR\cong I^{[q]}$ implies $^eR/I$ $ ^eR\cong R/I^{[q]}$.
\pre

\section{The Hilbert-Kunz function of a binoid}

In this section we introduce the Hilbert-Kunz function of a binoid.

\lem
Let $N$ be a finitely generated, semipositive binoid, $T$ a finitely generated $N$-set and $\mathfrak{n}$ an $N_+$-primary ideal of $N$. Then for every positive integer $q$ we have $$\# T/([q]\mathfrak{n}+T)\leqslant r|N^\times|+D q^s,$$ where $r$ is the number of generators of $T$, $s$ is the number of generators of $N_+$ and $D$ is some constant.
\leme
\pr
Let $t_1,\dots,t_r$ be generators of $T$ and $n_1,\dots,n_s$ be generators of $N_+$. If $t\in T$ then either $$t=u+t_i,$$ where $u\in N^\times$, or $$t=a_1n_1+\cdots+a_sn_s+t_i,$$ where $a_j\in \n,1\leqslant i \leqslant r$.
There are at most $r|N^\times|$ elements of the first type. In the second case we know that there exist $d_i\in\n$ such that 
$d_in_i\in \mathfrak{n}$ for $1\leqslant i \leqslant s$. So if $a_j\geqslant qd_j$, for some $j$, then $t\in [q]\mathfrak{n}+T$, 
which means 
\[\# T/([q]\mathfrak{n}+T)\leqslant r|N^\times|+r\cdot q^s\cdot\prod_{i=1}^s d_i.\qedhere\]
\pre
\df
Let $N$ be a finitely generated, semipositive binoid, $T$ a finitely generated $N$-set and $\mathfrak{n}$ an $N_+$-primary ideal of $N$. Then we call the number $$\hkf^N(\mathfrak{n},T,q)=\hkf(\mathfrak{n},T,q):=\# T/([q]\mathfrak{n}+T)=|T/([q]\mathfrak{n}+T)|-1$$ the \emph{Hilbert-Kunz function} of $\mathfrak{n}$ on the $N$-set $T$ at $q$.
\dfe
In particular, for $T=N$ and $\mathfrak{n}=N_+$ we have $$\hkf(N,q):=\hkf(N_+,N,q)=\#N/[q]N_+.$$
\ex
If $N=\z^\infty$ then $N_+=\{\infty\}$ and $N/[q]N_+=N$ is not a finite set.
\exe
\ex
\label{N^n}
Let $N=(\n^n)^\infty$ then $\hkf(N,q)=\# N/[q]N_+=q^n$. Since $N_+=(\n^n)^\infty\setminus \{0\}$ we have $[q]N_+=\bigcup_{i=1}^n qe_i+N$, where $e_i$ is the $i$-th standard vector of $\n^n$. Hence $$N/[q]N_+=\{(a_1,\dots,a_n)\in \n^n\mid 0\leqslant a_i\leqslant q-1\}\cup\{\infty\},$$ which means $\# N/[q]N_+=q^n$.
\exe
\lem
\label{HKF of annihilator}
Let $N$ be a finitely generated and semipositive binoid, $(T,p)$ a finitely generated $N$-set and $\mathfrak{n}$ an $N_+$-primary ideal of $N$. Let $\mathfrak{a}\subseteq \ann T$ be an ideal of $N$, where $\ann T$ is the annihilator of $T$. Then $T$ is also an $N/\mathfrak{a}$-set and we have  $$\hkf^N(\mathfrak{n},T,q)=\hkf^{N/\mathfrak{a}}((\mathfrak{n}\cup\mathfrak{a})/\mathfrak{a},T,q).$$
\leme
\pr
Let us first show that $(\mathfrak{n}\cup\mathfrak{a})/\mathfrak{a}$ is an $(N/\mathfrak{a})_+$-primary ideal of $N/\mathfrak{a}$. If there is an element $[m]\in (N/\mathfrak{a})_+$, where $m\in N$ then there exists $l\in\n$ such that $lm\in\mathfrak{n}$ and so $l[m]=[lm]\in(\mathfrak{n}\cup\mathfrak{a})/\mathfrak{a}$.

Let us define an action $(N/\mathfrak{a})\times T\longrightarrow T$ by $(\overline{n},t)\longrightarrow \overline{n}+t=n+t$. Then it is easy to see that this is well defined, that 
$$\begin{tikzpicture}[node distance=2cm, auto]
  \node (N) {$N\times T$};
  \node (T) [right of=N]  {$T$};
  \node (N1) [below of=N]   {$(N/\mathfrak{a})\times T$};
  \node (T1) [below of=T]   {$T$};
  \draw[->] (N) to node {}(T);
  \draw[->] (N) to node {}(N1);
  \draw[->] (N1) to node {}(T1);
  \draw[->] (T) to node {}(T1);
  \end{tikzpicture}
$$
commutes and so $T$ is an $N/\mathfrak{a}$-set. For every $q\in \n$ we have 
$[q]\mathfrak{n}+T=[q](\mathfrak{n}\cup\mathfrak{a})/\mathfrak{a}+T$, as this is the image of $[q]\mathfrak{n}\times T$. Therefore 
\[\#T/([q]\mathfrak{n}+T)=\#T/([q](\mathfrak{n}\cup\mathfrak{a})/\mathfrak{a}+T).\qedhere\]
\pre
\lem
\label{HKF ineq}
Let $N$ be a finitely generated, semipositive binoid, $S$ and $T$ are finitely generated $N$-sets and $\mathfrak{n}$ an $N_+$-primary ideal of $N$. Suppose that we have a surjective $N$-set homomorphism $\phi:S\rightarrow T$. Then for all $q$ $$\hkf(\mathfrak{n},S,q)\geqslant \hkf(\mathfrak{n},T,q).$$
\leme
\pr
By definition $[q]\mathfrak{n}$ is an ideal of $N$. So by Lemma \ref{surj hom} we have a surjective homomorphism 
$S/([q]\mathfrak{n}+S)\longrightarrow T/([q]\mathfrak{n}+T)$. Hence 
\[\hkf(\mathfrak{n},S,q)= |S/([q]\mathfrak{n}+S) |-1\geqslant | T/([q]\mathfrak{n}+T) |-1=\hkf(\mathfrak{n},T,q).\qedhere\]
\pre
\df
Let $N$ be a finitely generated, semipositive binoid, $T$ a finitely generated $N$-set and $\mathfrak{n}$ an $N_+$-primary ideal of $N$. Then the \emph{Hilbert-Kunz multiplicity} of $\mathfrak{n}$ on the $N$-set $T$ is defined by $$e_{HK}(\mathfrak{n},T):=\lim_{q\rightarrow \infty} \frac{\hkf^N(\mathfrak{n},T,q)}{q^{\dim N}},$$ if this limit exists.

In particular for $T=N$ and $\mathfrak{n}=N_+$ we set $e_{HK}(N):=e_{HK}(N_+,N)$ and denote it the \emph{Hilbert-Kunz multiplicity} of $M$.
\dfe
Note that here we work with the combinatorial dimension of $N$, which is not the appropriate definition in the non-cancellative case.
\tr
\label{bound of HKF}
Let $N$ be a finitely generated, semipositive binoid and $\mathfrak{n}$ an $N_+$-primary ideal of $N$. If $T$ is a finitely generated $N$-set and $\hkf(\mathfrak{n},N,q)\leqslant Cq^{\dim N}$ for every $q\in \n$ and some $C$ constant (in particular, if $e_{HK}(\mathfrak{n},N)$ exists) then there exists $\alpha$ such that $\hkf(\mathfrak{n},T,q)\leqslant \alpha q^{\dim N}$.
\tre
\pr
By Lemma \ref{canonical isomorphism}, we have a canonical isomorphism $N^{\cupdot r}/([q]\mathfrak{n}+N^{\cupdot r})\longrightarrow (N/[q]\mathfrak{n})^{\cupdot r}$, which means that $\# N^{\cupdot r}/([q]\mathfrak{n}+N^{\cupdot r})=\#(N/[q]\mathfrak{n})^{\cupdot r}= r\cdot \#(N/[q]\mathfrak{n})\leqslant r\cdot Cq^{\dim N}$.

From Lemma \ref{surj map for f.g. set}, we have a surjective map $\phi:N^{\cupdot r}\longrightarrow T$, where $r$ is the number of generators of $T$. Also by Lemma \ref{surj hom}, we know that there exists a surjective homomorphism $$\tilde{\phi}:N^{\cupdot r}/([q]\mathfrak{n}+N^{\cupdot r}) \longrightarrow T/([q]\mathfrak{n}+T).$$  So by Lemma \ref{HKF ineq} we have $\hkf(\mathfrak{n},T,q)\leqslant \alpha q^{\dim N}$, where $\alpha=r\cdot C$.
\pre

\lem
\label{multiplicity of smash product}
Let $M$ and $N$ be finitely generated, semipositive binoids. Then we have $$\hkf(M\wedge N,q)=\hkf(M,q)\cdot \hkf(N,q).$$
\leme
\pr
We know that (see Theorem \ref{dimension of smash product}) $$(M\wedge N)_+=(M\wedge N)\setminus (M^\times\wedge N^\times)=(M_+\wedge N)\cup (M\wedge N_+).$$
Take an element $q(m\wedge n)+m_1\wedge n_1\in [q](M\wedge N)_+$, where $m\wedge n\in (M\wedge N)_+$. Then $m\wedge n\in M_+\wedge N$ or $m\wedge n\in M\wedge N_+$, so $qm+m_1\in [q]M_+$ or $qn+n_1\in [q]N_+$, which means $$[q](M\wedge N)_+\subseteq ([q]M_+\wedge N)\cup (M\wedge [q]N_+).$$
We can also easily check the other inclusion, so we have
$$[q](M\wedge N)_+= ([q]M_+\wedge N)\cup (M\wedge [q]N_+).$$
Hence, from this result and Lemma \ref{smash of quotients}, we get
$$(M\wedge N)/[q](M\wedge N)_+\cong (M/[q]M_+)\wedge (N/[q]N_+).$$
By assumption $M,N$ are semipositive, so we know $M/[q]M_+$ and $N/[q]N_+$ are finite biniods. Hence, by Proposition
\ref{number of elements of smash}, we can conclude that
\[\# \big((M\wedge N)/[q](M\wedge N)_+\big)=\# (M/[q]M_+) \cdot \# (N/[q]N_+). \qedhere \]
\pre
\tr
\label{multiplicity of smash}
Let $M$ and $N$ be binoids such that $e_{HK}(M)$ and $e_{HK}(N)$ exist. Then $e_{HK}(M\wedge N)$ exists and $$e_{HK}(M\wedge N)=e_{HK}(M)\cdot e_{HK}(N).$$
\tre
\pr
By definition
\begin{align*}
e_{HK}(M)\cdot e_{HK}(N)&=\lim_{q\rightarrow \infty} \frac{\hkf(M,q)}{q^{\dim M}}\cdot \lim_{q\rightarrow \infty} \frac{\hkf(N,q)}{q^{\dim N}}\\
&=\lim_{q\rightarrow \infty} \frac{\hkf(M,q)\cdot \hkf(N,q)}{q^{\dim M+\dim N}}\\
&\overset{\text{Lemma } \ref{multiplicity of smash product}}{=}\lim_{q\rightarrow \infty} \frac{\hkf(M\wedge N,q)}{q^{\dim M+\dim N}}\\
&\overset{\text{Theorem } \ref{dimension of smash product}}{=}\lim_{q\rightarrow \infty} \frac{\hkf(M\wedge N,q)}{q^{\dim M\wedge N}}\\
&=e_{HK}(M\wedge N).
\end{align*}
(here by assumption $e_{HK}(M)$ and $e_{HK}(N)$ exist, so $M,N$ are semipositive binoids) 
\pre
\rem
We have the same results for primary ideals in Lemma \ref{multiplicity of smash product} and Theorem \ref{multiplicity of smash} and the proofs are the same.
\reme

\section{Hilbert-Kunz function of finite binoids}
In this section we will compute the Hilbert-Kunz function of finite binoids and separation is finite.
\prop
\label{HKfiniteness}
Let $N$ be a finitely generated, semipositive, separated binoid. Then the following are equivalent.
\begin{enumerate}
\item The Hilbert-Kunz function $\hkf(N,q)$ is eventually constant.
\item $N$ is finite.
\end{enumerate}
\prope
\pr
1) to 2). Since $\hkf(N,q)$ is eventually constant so $[q]N_+$ is also eventually constant and by separateness we have $\bigcap_q [q]N_+=\{\infty\}$. So there exists $q_0\in\n$ such that $[q]N_+=\{\infty\}$ for every $q\geqslant q_0$. Let $N$ be generated by $e_1,\dots,e_n$. Then for every $u\in N$ we can write $u=\sum_{i=1}^n a_ie_i$ and if $\sum_{i=1}^n a_i>nq_0$ then $u=\infty$. So $N$ is finite.

2) to 1). Let $a\in N$. If $na=ma\neq \infty,n\neq m\in \n_+$ then $a\in\bigcap_q [q]N_+$, but this would contradict separateness. Since $N$ is finite, we have $m_aa=\infty$ for some $m_a\in \n$. So $[q]N_+=\{\infty\}$ for $q\geqslant \max_{a\in N} m_a$. This means that $N/[q]N_+=N$ and so $\hkf(N,q)=\# N$ for $q$ large enough.
\pre
\cor
Let $N$ be a finite and separated binoid. Then $$\# N=\lim_{q\rightarrow \infty} \hkf(N,q).$$
\core
\pr
Our binoid is finite, so it is finitely generated and semipositive. Hence this equality was proved in Proposition \ref{HKfiniteness}.
\pre
\ex
\label{finite group's HK}
Let $T$ be a finite group. Then it is easy to check that $HKF(T^\infty,q)=|T|$ and in particular $e_{HK}(T^\infty)=|T|$.
\exe
\prop
\label{HK for separated}
If $N$ is a finitely generated, semipositive binoid then we have $$\hkf(N^{\sep},q)=\hkf(N,q).$$
\prope
\pr
We have the following commutative diagram
$$\begin{tikzpicture}[node distance=2cm, auto]
  \node (N) {$N$};
  \node (Nq) [right of=N]  {$N/{[q]N_+}$};
  \node (N1) [below of=N]   {$~~~N/{\bigcap_n nN_+}$};
  \draw[->] (N) to node {$f$} (Nq);
  \draw[->] (N) to node [swap] {$g$} (N1);
  \draw[->, dashed] (N1) to node {$h$} (Nq);
\end{tikzpicture}
$$
where $f$ and $g$ are the natural surjective homomorphisms and $\ker f=[q]N_+$. By Proposition \ref{ordinaryforbinuos} we get $\bigcap_n nN_+\subseteq [q]N_+=\ker f$. So by \cite[Lemma 1.6.3]{Simone}, there exists a unique surjective homomorphism $$h:N/{\bigcap nN_+}\rightarrow N/{[q]N_+}$$ such that $f=h\circ g$.
If we denote $M:=N^{\sep}=N/{\bigcap_n nN_+}$ and take $x\in [q]M_+$ then $x=\infty\in M$ or $x=y+qn$ where $y\in M^\bullet, n\in (M_+)^\bullet$. In the first case we have $x\in \ker h$. In the second case we have an equation $\widetilde{x}=\widetilde{y}+q\widetilde{n}$ in $N$.
This means that $\widetilde{x}\in [q]N_+=\ker f$. So $\infty=f(\widetilde{x})=h(g(\widetilde{x}))=h(x)$ and we have $[q]N_+\subseteq \ker h$.
Let us take an element $x\in \ker h$. Then there exists $\widetilde{x}\in N$ such that $g(\widetilde{x})=x$. Hence $f(\widetilde{x})=h(g(\widetilde{x}))=h(x)=\infty$ so $\widetilde{x}\in [q]N_+$. From here we can write $\widetilde{x}=z+qm$, $z\in N,m\in N_+$ and $x=g(\widetilde{x})=g(z+qm)=g(z)+q g(m)\in [q]M_+$. So this means that
$\ker h=[q]M_+$ and from here we can conclude 
\[M/[q]M_+\cong N/[q]N_+.\qedhere\]
\pre
\prop
\label{finitenessequivalence}
Let $N$ be a finitely generated, semipositive binoid. Then the following are equivalent.
\begin{enumerate}
\item $\# N^{\sep}=\#(N/\bigcap_n nN_+)$ is finite.
\item $\hkf(N,q)$ is eventually constant.
\end{enumerate}
\prope
\pr
From Proposition \ref{HK for separated} it follows that $\hkf(N^{\sep},q)=\hkf(N,q)$ and if we use Proposition \ref{HKfiniteness} then we are done.
\pre
\ex
Let $N=(\n^2)^{\infty}/(3x=5x+2y,3y=2x+5y)$ and $q\geqslant 3$. Every element $ax+by$, where $a\geqslant 3$ or $b\geqslant 3$ is equal to $\infty$ in $N/[q]N_+$. Because if $a\geqslant 3$ then $ax+by=(a-3)x+3x+by=(a-3)x+5x+2y+by=(a+2)x+(b+2)y=\dots=(a+2k)x+(b+2k)y$ for every $k\in \n$.
From here there exists $k$ such that $a+2k>q$ so $ax+by=\infty$. By this result we can easily see that $N^{\sep}=N/(3x,3y)=(\n^2)^\infty/(3x,3y)$ and so $\# (N/[q]N_+)=9$ for every $q\geqslant 3$. In particular $N^{\sep}$ is finite. The binoid homomorphism $N\rightarrow\z$ given by $x\mapsto 1$ $y\mapsto -1$, shows that for every $c\neq d\in \n$ we have $cx\neq dx$ in $N$ and so $N$ has infinitely many elements. So in  Proposition \ref{finitenessequivalence} we can not change $N^{\sep}$ to $N$.
\exe
\ex
Let $N=(\q_{\geqslant 0},+,\infty)$. For every $q\in \n \setminus \{0\}$ we have $qN_+=N_+$ and $[q]N_+=N_+$. So $\hkf(N^{\sep},q)=\hkf(N,q)=1$.
\exe
\ex
Let $N=\langle e_n \mid n\in \n\rangle$ be endowed with the following infinitely many relations for generators
\begin{align*}
e_1&=e_2+e_3,\\
e_2&=e_4+e_5+e_6\\
e_3&=e_7+e_8+e_9\\
e_4&=e_{10}+e_{11}+e_{12}+e_{13},
\end{align*}
and etc.
\begin{align*}
e_{\sum_{k=1}^{n-1}k!+1}&=e_{\sum_{k=1}^nk!+1}+\cdots+e_{\sum_{k=1}^n k!+n+1}\\
e_{\sum_{k=1}^{n-1}k!+2}&=e_{\sum_{k=1}^nk!+n+2}+\cdots+e_{\sum_{k=1}^n k!+2(n+1)}\\
& \hspace*{2cm}\vdots \\
e_{\sum_{k=1}^{n-1}k!+n!}&=e_{\sum_{k=1}^nk!+(n+1)(n!-1)+1}+\dots+e_{\sum_{k=1}^n k!+(n+1)n!}
\end{align*}
for every $n\geqslant 2$.\\
We can easily check that $\bigcap_n nN_+=N_+$ and $\bigcap_n [n]N_+=\{\infty\}$. So $\hkf(N^{\sep},q)=1$ and $\hkf(N,q)$ is not even defined.
\exe

\section{Some Hilbert-Kunz functions of 2 and 3 generated binoids}

 Let $N$ be a commutative, finitely generated binoid with minimal generating family $x_1,\dots,x_r$. Then $N$ is isomorphic to $(\n^r)^\infty/\sim_\varepsilon$, where
 $$\varepsilon:(\n^r)^{\infty}\longrightarrow N$$
is the canonical binoid surjective homomorphism with $\varepsilon(e_i)=x_i,i=1,\dots,r$. Now, by R\'{e}dei's theorem (\cite[Theorem 1.6.9]{Simone}), the congruence $\sim_\varepsilon$ can be generated by finitely many relations $\mathcal{R}:\varepsilon(a)=\varepsilon(b),a,b\in(\n^r)^{\infty}$.
 Hence $N$ is given by finitely many generators and relations, which we denote by $$(\n^r)^{\infty}/(\mathcal{R}_1,\dots,\mathcal{R}_n).$$
 So we can think of the elements in the binoid as points in $\n^r$ modulo these relations.
 \df
A relation of the form $\mathcal{R}:\varepsilon(x)=\infty$ is called a \emph{monomial} relation. A \emph{binomial} relation is a relation of the form $\mathcal{R}:\varepsilon(a)=\varepsilon(b)\neq \infty$.
\dfe
\lem
Let $N$ be a commutative, finitely generated and positive binoid. Then it has the maximal ideal $N_+=N\setminus \{0\}=\langle x_1,\dots,x_r\rangle$, where $\{x_1,\dots,x_r\}$ is the minimal generating system of $N$.
\leme
\pr
Since $x_i$ is a minimal generator, we have $x_i\neq 0$ and so $x_i$ is a non-unit, hence $x_i\in N_+$ and so $\langle x_1,\dots,x_r\rangle\subseteq N_+$. If $x\in N_+$, then $x=\sum_{i=1}^r n_ix_i$ and at least one $n_i\geqslant 1$. Then $x\in \langle x_1,\dots,x_r\rangle$.
\pre
From this Lemma and the isomorphism $(\n^r)^\infty/[q]\n_+^r\rightarrow N/[q]N_+$ we can conclude that $\# N/[q]N_+$ is equal to the number of all equivalence classes in $[0,q-1] \times\cdots\times [0,q-1]\subseteq \n^r$ except the class of $\infty$, where the relations are $qx_i\sim \infty,i=1,\dots,r$ and the induced binoid relations.
With this direct approach we can compute the Hilbert-Kunz function, if the relations are easy. So now we will compute the Hilbert-Kunz function of some 2 or 3 generated positive binoids.

Let $N=(\n^2)^{\infty}/(nx_1+mx_2=kx_1+lx_2)$ be a binoid where $n,m,k,l\in \n$. If $n=m=0$ or $k=l=0$ then there is a relation where one side is 0. Hence $N$ is not positive and $N_+\neq \langle x_1,x_2\rangle$. So we have $n+m>0$ and $k+l>0$.
The equivalence relation $\sim$ is given by:
\begin{equation}
(i+n,j+m)\sim(i+k,j+l)
\end{equation}
 for every $i,j \in \n$ and $(a,b)\sim \infty$ in $N/[q]N_+$, where $q\gg 0$. This means that $(a,b)\sim (i,j)$ for some $i,j\in \n$ where $i\geqslant q$ or $j \geqslant q$.

\df
Let $N=(\n^r)^{\infty}/\sim$ be a binoid and suppose that one of the relation is a binomial relation
$$n_1x_1+\dots+n_rx_r= m_1x_1+\dots+m_rx_r$$ with $\sum n_i\geqslant 1,\sum m_i\geqslant 1$ and $(n_1,\dots,n_r)\neq (m_1,\dots,m_r)$. If  $n_1\geqslant m_1,\dots,n_r\geqslant m_r$ or $n_1\leqslant m_1,\dots,n_r\leqslant m_r$ then we call it a \emph{positive} relation otherwise a \emph{negative} relation.
Also we call the vector $(n_1-m_1,\dots,n_r-m_r)$ the \emph{relation vector} or just \emph{vector}.
\dfe

\prop
\label{monomial to binomial}
Let $N=(\n^2)^{\infty}/(nx_1+mx_2=kx_1+lx_2)$, $M=(\n^2)^{\infty}/(nx_1+mx_2=\infty)$ be two binoids where $n,m,k,l,q\in \n,q\gg 0$ and $k+l>0$. If $n<k$ and $m<l$ then $$\hkf(N,q)=\hkf(M,q)=(n+m)q-nm.$$
\prope
\pr
We have $(i+n,j+m)\sim(i+k,j+l)$ for all $i,j\in \n$, so with $i=k-n,j=l-m$ we get $$(n,m)\sim (k,l)\sim(2k-n,2l-m)\sim (n+s(k-n),m+s(l-m)),$$ for all $s\in\n$. Hence, if $q$ is large enough then $(n,m)\sim \infty$ in $N/[q]N_+$, which means $\hkf(N,q)=\hkf(M,q)=(n+m)q-nm$.
\pre
\lem
\label{integral iff binomial}
Let $N=(\n^r)^{\infty}/\sim$ be a positive binoid. Then all the relations are binomial if and only if $N$ is integral.
\leme
\pr
All the relations are binomial or monomial.
If all the relations are binomial then $\sum_{i=1}^r n_ix_i\neq\infty$ for all $n_i\in\n$. If a binoid is integral then we do not have monomial relations.
\pre
\prop
\label{make integral}
Let $N$ be a commutative, finitely generated and positive binoid. Then there exists an integral binoid which has the same Hilbert-Kunz function as $N$ for $q\gg 0$.
\prope
\pr
We can assume that our binoid is $(\n^r)^{\infty}/\sim$, with finitely many relations. If the binoid is not integral then by Lemma \ref{integral iff binomial} there are monomial relations.
Let $M$ be a binoid with the same generators $x_1,\dots,x_r$ and same binomial relations but every monomial relation $\sum_{i=1}^r n_ix_i=\infty$ is replaced by $\sum_{i=1}^r n_ix_i=\sum_{i=1}^r (n_i+1)x_i$. Then we can show that $\hkf(N,q)=\hkf(M,q)$ by the same argument as in proof of Proposition \ref{monomial to binomial}.
\pre
\rem
Here a binoid became integral but it looses cancellativity. So we do not know the right dimension of the new binoid.
\reme
\prop
  Let $N=(\n^2)^{\infty}/(nx_1+mx_2=kx_1+lx_2)$ be a binoid where $n,m,k,l,q\in \n,q\gg 0$ and $n>k,m<l$. Then we have
$\# N/[q]N_+=(k+m)q-km+\# N'$, where $$N'=(\n^2)^{\infty}/((n-k)x_1=(l-m)x_2,(q-k)x_1=(q-m)x_2=\infty).$$
\prope
\pr
If $x<k,y<q$ or $x<q,y<m$ then all the points $(x,y)$ have a different equivalence class. Hence $\# N/[q]N_+=(k+m)q-km+\# N'$.
\pre
\prop
\label{HKF of dim 2}
Let $N'=(\n^2)^{\infty}/((n-k)x_1=(l-m)x_2,(q-k)x_1=(q-m)x_2=\infty)$ be a binoid and $q\gg 0$. Let
$q-m=(l-m)s+t, s,t \in \n,0\leqslant t<l-m$ (without loss of generality we can assume that $m\geqslant k$). Then we have the following results.
\begin{enumerate}
   \item If $l-m>n-k$ or ($l-m=n-k$ and $n-k\leqslant m+t$) then $\# N'=(n-k)(q-m)$,
   \item If $l-m=n-k$ and $n-k>m+t$ then $\# N'=(n-k)(q-m)-t(n-k-m-t)$,
   \item If $l-m<n-k$ then $\# N'=(l-m)(q-k)$.
\end{enumerate}
\prope
\pr
If $l-m> n-k$ then the points $(x,y)$, where $x<n-k,y<q-m$, are all in different equivalence classes and are not equivalent to $\infty$, since $x+(n-k)s<q$. It is also clear that all other points are equivalent to these points or $\infty$. From here $\# N'$ is equal to $(n-k)(q-m)$. If $l-m=n-k$ then $(x,y)\sim\infty$, where $x\geqslant m+t$ and $y\geqslant q-m-t$ so $\# N'$ is equal to $(n-k)(q-m)-t(n-k-m-t)$ if $n-k>m+t$, or $(n-k)(q-m)$ otherwise.

If $l-m<n-k$ then the computation is the same.
\pre

Let $N=(\n^2)^{\infty}/(nx_1+mx_2=kx_1+lx_2,n_1x_1+m_1x_2=k_1x_1+l_1x_2)$ be a binoid, where $n,m,k,l,n_1,m_1,k_1,l_1,q\in \n$ and $q\gg 0$.
The equivalence relations are given by:
\begin{equation}
(i+n,j+m)\sim(i+k,j+l)
\end{equation}
\begin{equation}
(i+n_1,j+m_1)\sim(i+k_1,j+l_1)
\end{equation}
for every $i,j \in \n$. $(x,y)\sim\infty$ in $N/[q]N_+$ means that$(x,y)\sim (i,j)$
for some $i,j\in \n$, where $i\geqslant q$ or $j \geqslant q$.

\prop
Let $N=(\n^2)^{\infty}/(nx_1+mx_2=kx_1+lx_2,n_1x_1+m_1x_2=k_1x_1+l_1x_2)$ be a binoid, where $n,m,k,l,n_1,m_1,k_1,l_1,q\in \n$ and $n+m>0,k+l>0,n_1+m_1>0,k_1+l_1>0$. If $q$ is a large enough integer and $(n-k)(m_1-l_1)\neq (n_1-k_1)(m-l)$ then $$\# N/[q]N_+=(a+b)q+C,$$ where $C$ is a constant and $a=\min(n,k,n_1,k_1),b=\min(m,l,m_1,l_1)$.
\prope
\pr
All the points $(x,y)$ have a different equivalence class in $N/[q]N_+$, where $x<a$ or $y<b$. Hence we have $$\# N/[q]N_+=(a+b)q-ab+\# N',$$ where $N'=(\n^2)^{\infty}/((n-a)x_1+(m-b)x_2=(k-a)x_1+(l-b)x_2,(n_1-a)x_1+(m_1-b)x_2=(k_1-a)x_1+(l_1-b)x_2,(q-a)x_1=\infty,(q-b)x_2=\infty)$. So for the relations in $N'$ we have the following cases:
\begin{enumerate}
 \item $Ax_1=Bx_1+Cx_2$, $Dx_1+Ex_2=Fx_2$ and $(A-B)(F-E)\neq CD$
 \item $Ax_1=Bx_2$, $Cx_1+Dx_2=Ex_1+Fx_2$ and $A(D-F)\neq B(E-C).$
\end{enumerate}
In the first case we have
$$\begin{aligned}
(AF+BE)x_1&=F(Bx_1+Cx_2)+BEx_1\\
&=FBx_1+C(Dx_1+Ex_2)+BEx_1\\
&=FBx_1+CDx_1+E(Cx_2+Bx_1)\\
&=(FB+CD+EA)x_1
\end{aligned}
$$
and similarly we can show $$(AF+BE)x_2=(FB+CD+EA)x_2.$$ By assumption $AF+BE\neq FB+CD+EA$, hence we can conclude that
$$\begin{aligned}
(AF+BE)x_1&=(FB+CD+EA)x_1=(AF+BE)x_2\\
&=(FB+CD+EA)x_2=\infty\in N'/[q]N'_+.
\end{aligned}$$
So $\# N'/[q]N'_+ \leqslant (\min(AF+BE,FB+CD+EA))^2$.\\

In the second case we can show $$\begin{aligned}
                                  (AD+BC)x_1=(BE+AF)x_1,\\
                                  (AD+BC)x_2=(BE+AF)x_2,
                                 \end{aligned}$$
and by assumption $AD+BC\neq BE+AF$, so $\# N'/[q]N'_+ \leqslant (\min(AD+BC,BE+AF))^2$.
\pre

\prop
Let $N=(\n^2)^{\infty}/(nx_1+mx_2=kx_1+lx_2,n_1x_1+m_1x_2=k_1x_1+l_1x_2)$ be a binoid, where $n,m,k,l,n_1,m_1,k_1,l_1,q\in \n$ and $n+m>0,k+l>0,n_1+m_1>0,k_1+l_1>0$. Without loss of generality  we can assume that $n>k,m<l,n_1>k_1,m_1<l_1$ and let
\begin{align}
(n-k,m-l)&=i(A,-B),\\
(n_1-k_1,m_1-l_1)&=j(A,-B)
\end{align}
where $i,j,A,B\in \n$ and $\operatorname{gcd}(i,j)=1$. If $q$ is a large enough integer then we have $$\# N/[q]N_+=(a+b+\operatorname{min}(A,B))q+\operatorname{Const},$$ where $a=\min(k,k_1),b=\min(m,m_1)$.
\prope
\pr
All the points $(x,y)$ have a different equivalence class in $N/[q]N_+$, where $x<a$ or $y<b$. Hence we have $$\# N/[q]N_+=(a+b)q-ab+\# N',$$ where $N'=(\n^2)^{\infty}/((n-a)x_1+(m-b)x_2=(k-a)x_1+(l-b)x_2,(n_1-a)x_1+(m_1-b)x_2=(k_1-a)x_1+(l_1-b)x_2,(q-a)x_1=\infty,(q-b)x_2=\infty)$.
From $\operatorname{gcd}(i,j)=1$ there exist $u,v \in \n$ such that $iu=jv\pm 1$.
If we have
\begin{align}
iu=jv+1
\end{align}
 then $u(n-k,m-l)-v(n_1-k_1,m_1-l_1)=(iu-jv)(A,-B)=(A,-B)$.
So it is easy to check that
$$\begin{aligned}
(x,y)\sim (x+A,y-B),
\end{aligned}$$
where $x\geqslant v(n_1-k_1)+k_1,y\geqslant B+b$ or $x\geqslant a,y\geqslant u(l-m)+m.$ From $(2.4)$ and $(2.5)$ we get $$n-k=iA,l-m=iB,n_1-k_1=jA,l_1-m_1=jB,$$ and we also know $u(l-m)=uiB$, $v(n_1-k_1)=vjA \stackrel{(2.6)}{=}(ui-1)A$. So we can find some integer $C$ such that
$$(x+a,y+CuiB+b)\sim (x+A+a,y+(Cui-1)B+b)\sim \ldots $$$$\sim (x+(Cui-1)A+a,y+B+b)\sim (x+CuiA+a,y+b),$$ for every $x,y\in \n$, where for every $k\in \n$ we should have $$x+kA+a\geqslant v(n_1-k_1)+k_1$$ or $$y+(Cui-k)B+b\geqslant u(l-m)+m.$$
So now it is easy see that $\# N/[q]N_+=(a+b+\operatorname{min}(A,B))q+\text{Const}$.

For the case $nu+1=mv$, the proof is exactly the same as in the case $nu=mv+1$.
\pre

If we have more than 2 relations (vectors) then there are two cases: there exist two non-collinear vectors or all vectors are collinear. Then we can deduce similar results as above.\\

Let $N=(\n^3)^{\infty}/(nx_1=mx_2+kx_3)$ be a binoid, where $n,m,k,q\in \n_+$. The equivalence
 relation is given by:
 \begin{equation}(a,b,c)\sim(a-n,b+m,c+k)
 \end{equation}
where $a\geqslant n$ and $a,b,c\in \n$.

Also $(a,b,c)=\infty$ in $N/[q]N_+$ means that $(a,b,c)\sim (i,j,k)$ for some $i,j,k\in \n$, where $i\geqslant q$ or $j \geqslant q$ or
 $k \geqslant q$.
\prop
\label{HKF of 3 generated binoid}
Let $N=(\n^3)^{\infty}/(nx_1=mx_2+kx_3)$ be a binoid, where $n,m,k,q\in \n_+$. Let $q\gg 0$ and write $q=ns+t,0\leqslant t<n$. Then we have the following results.
\begin{enumerate}
\item If $n<\max(m,k)$ then $\# N/[q]N_+=nq^2$,
\item If $n>\max(m,k)$ then $\# N/[q]N_+=(m+k-\frac{mk}{n})q^2-\frac{mk}{n}t(n-t)$,
\item If $n=m$ and $n\geqslant k$ then $$\# N/[q]N_+=nq^2-\dfrac{t(n-t)(n-k)}{n}q-kt^2+\frac{kt^3}{n},$$
\end{enumerate}
\prope
\pr
\begin{enumerate}
\item If $n<\max(m,k)$ then the points $(x,y,z)\in [0,n-1] \times [0,q-1]\times [0,q-1]$ are neither equivalent to $\infty$ nor one to each other. So $\# N/[q]N_+=nq^2$.
\item If $n>\max(m,k)$ then all the points where $x<n$ are not equivalent to each other, so we only have to count the points that are not infinity. Hence all the points in $\{ (x,y,z): x<n,y<ms,z<q\}$ and $\{ (x,y,z): x<n,y<q,z<ks\}$ are not equivalent to $\infty$. But we have double counted the points in $\{ (x,y,z): x<n,y<ms,z<ks\}$.
It is also not difficult to see that all the points $\{ (x,y,z): x<t,ms\leqslant y<m(s+1),ks\leqslant z<q\}$ and $\{ (x,y,z): x<t,m(s+1)\leqslant y<q,ks\leqslant z<k(s+1)\}$ are not equivalent to $\infty$. Hence we have
 $$
 \begin{aligned}
\# N/[q]N_+&=|\{ (x,y,z): x<n,y<ms,z<q\}|\\
&\enspace+|\{ (x,y,z): x<n,y<q,z<ks\}|\\
&\enspace-|\{ (x,y,z): x<n,y<ms,z<ks\}|\\
&\enspace+|\{ (x,y,z): x<t,ms\leqslant y<m(s+1),ks\leqslant z<q\}|\\
&\enspace+|\{ (x,y,z): x<t,m(s+1)\leqslant y<q,ks\leqslant z<k(s+1)\}|\\
&=nmsq+nksq-nmsks+mt(q-ks)+kt(q-m(s+1))\\
&=(m+k-\frac{mk}{n})q^2-\frac{mk}{n}t(n-t).
 \end{aligned}
$$
\item From the above counting we can conclude that
\[\begin{aligned}
\# N/[q]N_+&=|\{ (x,y,z): x<n,y<ns,z<q\}|\\
&\enspace+|\{ (x,y,z): x<n,y<q,z<ks\}|\\
&\enspace-|\{ (x,y,z): x<n,y<ns,z<ks\}|\\
&\enspace+|\{ (x,y,z): x<t,ns\leqslant y<q,ks\leqslant z<q\}|\\
&=n^2sq+nksq-n^2ks^2+t^2(q-ks)\\
&=nq^2-\frac{t(n-t)(n-k)}{n}q-kt^2+\frac{kt^3}{n}.
\end{aligned}\qedhere
\]
\end{enumerate}
\pre

\section[Hilbert-Kunz multiplicity of toric binoids]{Hilbert-Kunz multiplicity of integral, cancellative, torsion-free and normal binoids}

Let $N$ be a finitely generated, positive, integral, cancellative and torsion-free binoid. Then we can assume that $(\diff N)^\times=\z^d$, where $d=\rk(N)$, and $C=\r_+N^\bullet\subseteq \r^d$ is the rational cone generated by $N$. From Chapter 1, Section 1.4, we know that $$C=H^+_{\sigma_1}\cap\cdots\cap H^+_{\sigma_s}$$ where $\sigma_i$ are the supporting linear forms.

Let $N$ be a finitely generated, positive, integral, cancellative and torsion-free binoid. Then the  Hilbert-Kunz multiplicity of the binoid algebra $K[N]$ is a rational number, which was already shown in the papers \cite{Bruns}, \cite{Watanabe}, \cite{Eto}. So now we will prove a similar theorem. The proof is similar to the proof of K. Eto.
\tr
\label{HKF for normal binoid}
Let $N$ be a finitely generated, positive, integral, cancellative, torsion-free and normal binoid. If $\mathfrak{n}$ is an $N_+$-primary ideal of $N$ then $$e_{HK}(\mathfrak{n},N)=\lim_{q\rightarrow \infty}\dfrac{\# N/[q]\mathfrak{n}}{q^{\dim N}}\in \q.$$
In particular it exists.
\tre
\pr
We can assume that $(\diff N)^\times=\z^d$, where $d=\dim N$. Let $\{v_1,\dots,v_k\}\subseteq N$ be the minimal set of generators of $N$. We know that $\mathfrak{n}$ is finitely generated so let $\{a_1,\dots,a_t\}$ be the set of generators of $\mathfrak{n}$. 
The ideal $\mathfrak{n}$ is $N_+$-primary so there exist $\alpha_i\in \n$ such that $\alpha_iv_i\in \mathfrak{n}$, for $1\leqslant i\leqslant k$, which means $\alpha_iv_i\in \langle a_1,\dots,a_t\rangle$. For every $0\neq u\in C=\r_+ N^\bullet$ we have $u=\sum_{i=1}^k c_iv_i$, $c_i\geqslant 0$, and there exists $j$ such that $c_j>0$. So we have  $$\alpha_j u=\sum_{i=1}^k \alpha_j c_iv_i=\alpha_jv_j+n\in a_l+C,$$ for some $l$, where $n\in N$. Hence the shifted cones $a_i+C$, $i=1,\dots,t$ define inside $C$ a bounded region which we denote by $B$, i.e. $$B=C\setminus \bigcup_{i=1}^t a_i+C.$$ The volume of $B$ is a rational number as it is given by rational cones.

Let $q$ be some positive integer and let us denote the sets $$\frac{1}{q}M=\Big\{\Big(\frac{n_1}{q},\dots,\frac{n_d}{q}\Big)\suchthat{\frac{1}{n}} (n_1,\dots,n_d)\in M\Big\}\subseteq \z^d\otimes_{\z}\q$$
for every binoid $M\subseteq \z^d$. There is a binoid isomorphism between $M$ and $\frac{1}{q}M$, and $\diff \frac{1}{q}M=\frac{1}{q}\z^d$. So it is easy to see that $$\# N/[q]\mathfrak{n}=\# N/\bigcup^t_{i=1}(qa_i+N)=\#\Big(\frac{1}{q}N/\bigcup^t_{i=1}(a_i+\frac{1}{q}N)\Big).$$
By normality of $N$ we have $\frac{1}{q}N^\bullet=\r_+N^\bullet\cap \frac{1}{q}\z^d=C\cap \frac{1}{q}\z^d$. 
So the elements of $\Big(\frac{1}{q}N/\bigcup^t_{i=1}(a_i+\frac{1}{q}N)\Big)^\bullet$ are the elements of $B\cap\frac{1}{q}\z^d$.

Hence the number of elements belonging to the set $\frac{1}{q}N/\bigcup^t_{i=1}(a_i+\frac{1}{q}N)$ divided by $q^d$ tends to the volume of 
$B$ as $q$ tends to infinity. Thus 
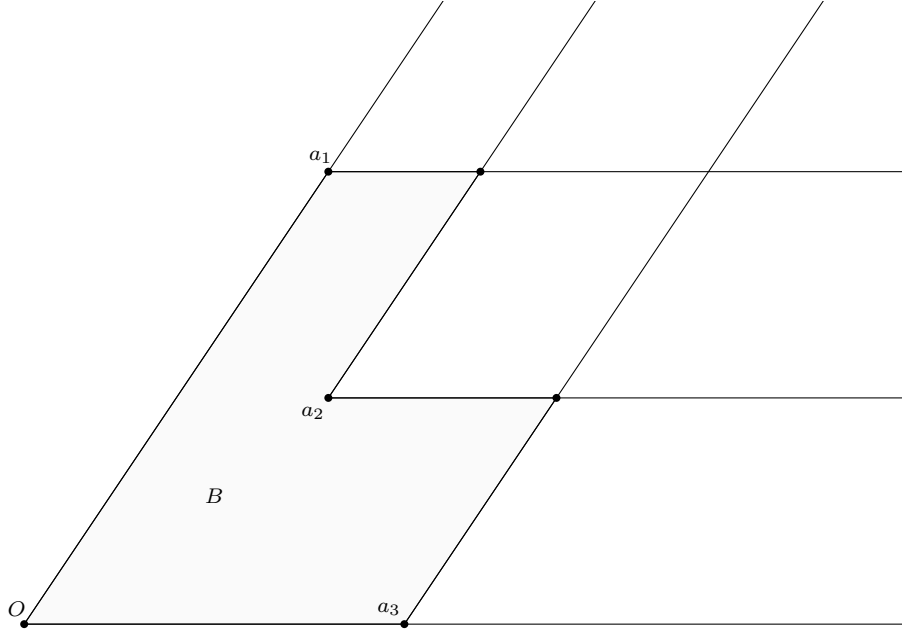
\[\lim_{q\rightarrow \infty}\dfrac{\# N/[q]\mathfrak{n}}{q^d}=\vol(B)\in \q.\qedhere\]
\pre
\definecolor{cccccc}{rgb}{0.8,0.8,0.8}

\begin{tikzpicture}[line cap=round,line join=round,>=triangle 45,x=1.0cm,y=1.0cm]

\clip(-0.95,-0.38) rectangle (11.54,8.25);
\fill[color=cccccc,fill=cccccc,fill opacity=0.1] (0,0) -- (4,6) -- (6,6) -- (4,3) -- (7,3) -- (5,0) -- cycle;
\draw [domain=0.0:11.544444] plot(\x,{(-0--6*\x)/4});
\draw [domain=0.0:11.544444] plot(\x,{(-0-0*\x)/5});
\draw [color=black] (0,0)-- (4,6);
\draw [color=black] (4,6)-- (6,6);
\draw [color=black] (6,6)-- (4,3);
\draw [color=black] (4,3)-- (7,3);
\draw [color=black] (7,3)-- (5,0);
\draw [color=black] (5,0)-- (0,0);
\draw [domain=4.0:11.544444] plot(\x,{(-6--3*\x)/2});
\draw [domain=4.0:11.544444] plot(\x,{(--9-0*\x)/3});
\draw [domain=4.0:11.544444] plot(\x,{(--12-0*\x)/2});
\draw [domain=5.0:11.544444] plot(\x,{(-15--3*\x)/2});
\begin{scriptsize}
\fill [color=black] (0,0) circle (1.5pt);
\draw[color=black] (-0.1,0.2) node {$O$};
\fill [color=black] (4,6) circle (1.5pt);
\draw[color=black] (3.9,6.2) node {$a_1$};
\fill [color=black] (5,0) circle (1.5pt);
\draw[color=black] (4.8,0.2) node {$a_3$};
\fill [color=black] (6,6) circle (1.5pt);
\fill [color=black] (4,3) circle (1.5pt);
\draw[color=black] (3.8,2.8) node {$a_2$};
\draw[color=black] (2.5,1.7) node {$B$};
\fill [color=black] (7,3) circle (1.5pt);
\end{scriptsize}
\end{tikzpicture}

\section{Reduction to the integral and reduced case}

In this section we show how to reduce the existence and rationality of the Hilbert-Kunz multiplicity to the integral and the reduced case.
\lem
\label{reduction of torsion free}
Suppose that for all finitely generated, semipositive, cancellative, integral and torsion-free binoids of dimension less or equal to $d$ the Hilbert-Kunz multiplicity exists. Then
\begin{enumerate}
 \item Let $N$ be a finitely generated, semipositive, cancellative, torsion-free, reduced  binoid of dimension $d$ and let $\{\mathfrak{p}_1,\dots,\mathfrak{p}_s\}$ be all minimal prime ideals of  dimension $d$. Then the Hilbert-Kunz multiplicity of $N$ exists, and $$e_{HK}(N)=\sum_{i=1}^s e_{HK}(N/\mathfrak{p}_i).$$
 \item Let $N$ be a finitely generated, semipositive, cancellative binoid of dimension $d$ and torsion-free except for nilpotent elements. Then we have $$\hkf(N,q)\leqslant Cq^{\dim N},$$ for all $q$, where $C$ is some constant.
\end{enumerate}
\leme
\pr
\begin{enumerate}
 \item We have finitely many minimal prime ideals $\{\mathfrak{p}_1,\dots,\mathfrak{p}_n\}$ and $n\geqslant s$. From reducedness we have $\bigcap_{i=1}^n \mathfrak{p}_i=\nil(N)=\{\infty\}$. If $f\notin [q]N_+$ and $f\in [q]N_+\cup \mathfrak{p}_i$ for all $i$, then $f\in \bigcap_{i=1}^n \mathfrak{p}_i=\{\infty\}$, which means $N\setminus[q]N_+\subseteq \bigcup_{i=1}^n N\setminus([q]N_+\cup\mathfrak{p}_i)$. Also by Proposition \ref{quotient to union} and $$N\setminus[q]N_+\cong (N/[q]N_+)\setminus \{\infty\},$$   $$N\setminus([q]N_+\cup\mathfrak{p}_i)\cong (N/([q]N_+\cup\mathfrak{p}_i))\setminus \{\infty\}$$ we can conclude that $$N\setminus[q]N_+\subseteq \bigcup_{i=1}^n (N/\mathfrak{p}_i)\setminus[q](N/\mathfrak{p}_i)_+.$$ But we know that $N/\mathfrak{p}_i$ is finitely generated, semipositive, cancellative, torsion-free and integral so by assumption we know that $$\hkf(N/\mathfrak{p}_i,q)=e_{HK}(N/\mathfrak{p}_i)q^{d_i}+O(q^{d_i-1}),$$ where $d_i:=\dim N/\mathfrak{p}_i$. Hence only the minimal prime ideals 
with dimension $d_i=d$ are important and we can write 
 \begin{equation}
\# N/[q]N_+\leqslant \sum_{i=1}^s e_{HK}(N/\mathfrak{p}_i)q^d+O(q^{d-1}).
 \end{equation}

Let us denote $H:=N/[q]N_+$. Then $$N/([q]N_+\cup\mathfrak{p}_i)=(H\setminus \mathfrak{p}_i)\cup\{\infty\}$$ and $$N/([q]N_+\cup\mathfrak{p}_i\cup\mathfrak{p}_j)=\big((H\setminus \mathfrak{p}_i)\cap (H\setminus \mathfrak{p}_j)\big)\cup\{\infty\}.$$ Hence from set theory we know that $$|H\setminus \{\infty\}|\geqslant |\bigcup_{i=1}^n H\setminus\mathfrak{p}_i|\geqslant \sum_{i=1}^n |H\setminus\mathfrak{p}_i|-\sum_{1\leqslant i<j\leqslant n}|(H\setminus \mathfrak{p}_i)\cap (H\setminus \mathfrak{p}_j)|.$$ So from here we have $$\#  N/[q]N_+\geqslant \sum_{i=1}^n \# N/([q]N_+\cup\mathfrak{p}_i)-\sum_{1\leqslant i<j\leqslant n}\# N/([q]N_+\cup\mathfrak{p}_i\cup\mathfrak{p}_j).$$
But we know that $N/(\mathfrak{p}_i\cup\mathfrak{p}_j)$ is finitely generated, semipositive, cancellative, torsion-free and integral so by Proposition \ref{quotient to union} and the assumption we get 
\begin{align*}
\sum_{i\neq j}\# N/([q]N_+\cup\mathfrak{p}_i\cup\mathfrak{p}_j)&=\sum_{i\neq j}\hkf(N/(\mathfrak{p}_i\cup\mathfrak{p}_j),q)\\
&=\sum_{i\neq j}\Big(e_{HK}(N/(\mathfrak{p}_i\cup\mathfrak{p}_j))q^{d_{ij}}+O(q^{d_{ij}-1})\Big), 
\end{align*}
 where $d_{ij}$ is the dimension of $N/(\mathfrak{p}_i\cup\mathfrak{p}_j)$. By Proposition \ref{dimension decreasing of union of minimal primes} we also know that $d_{ij}<\min\{d_i,d_j\}$. So we have 
\begin{equation}
\#  N/[q]N_+\geqslant \sum_{i=1}^s e_{HK}(N/\mathfrak{p}_i)q^d+O(q^{d-1}). 
\end{equation}
Hence from (2.8) and (2.9) we will get $$e_{HK}(N)=\lim_{q\rightarrow \infty} \dfrac{\hkf(N,q)}{q^d}=\sum_{i=1}^s e_{HK}(N/\mathfrak{p}_i).$$

\item We know that $\nil(N)$ is an ideal of $N$ and $N_{\red}:=N/\nil(N)$ is reduced. Let $a_1,\dots,a_s$ be the generators of $\nil(N)$ and $k_ia_i=\infty$ but $(k_i-1)a_i\neq \infty$, where $k_i\in \n,1\leqslant i\leqslant s$.
We have a finite decreasing sequence $$N=M_0\longrightarrow M_1\longrightarrow \cdots \longrightarrow M_{t-1}\longrightarrow M_t=N_{\red},$$ where $M_{k_0+\cdots+k_i+j}=M_{k_0+\cdots+k_i+j-1}/((k_{i+1}-j)a_{i+1}),1\leqslant j\leqslant k_{i+1}-1,0\leqslant i\leqslant s-1$ and $k_0=0$. Here we have $2(k_{i+1}-j)a_{i+1}=\infty$ in $M_{k_0+\cdots+k_i+j}$.
So in particular there exists a sequence such that $$M_{i+1}=M_i/(f_i),2f_i=\infty \operatorname{\;in\;}M_i,0\leqslant i\leqslant t-1$$ and $M_0=N,M_t=N_{\red}$. Hence there exists also such a sequence with this property of minimal length $l$.
We will use induction on $l$ to prove this Lemma.

For $l=0$ we are in the reduced situation and the statement follows from the case above. So suppose that $l$ is arbitrary and that for smaller $l$ the statement is already proven.

Now suppose we have a sequence with minimal length $l$ such that $$M_{i+1}=M_i/(f_i),2f_i=\infty\operatorname{\;in\;} M_i,0\leqslant i\leqslant l-1$$ and $M_0=N,M_l=N_{\red}$. So this means that we have a strongly exact sequence $$\infty \longrightarrow (f)+N\longhookrightarrow N\longtwoheadrightarrow N/(f)\longrightarrow\infty,$$
where $f=f_0$.  Hence $M:=N/(f)$ has a length of a decreasing sequence with the described property smaller than $l$, so by the induction hypothesis we know that $\# M/[q]M_+\leqslant Dq^d$ for some constant $D$. By Corollary \ref{equality of e.s. corollary} we also have $$\# N/[q]N_+ +\# \big((f)\cap [q]N_+\big)/((f)+[q]N_+)=\# ((f)+N)/((f)+[q]N_+)+\# M/([q]N_++M).$$
Since $2f=\infty$, we know that $f$ annihilates the $N$-set $(f)+N$. So $(f)+N$ is a finitely generated $M$-set ($f$ is the only $M$-set generator of $(f)+N).$
Hence by Lemma \ref{surj map for f.g. set} we have a surjective homomorphism $M\longrightarrow (f)+N$, which implies 
$$\# ((f)+N)/((f)+[q]N_+)\leqslant\# M/[q]M_+.$$ But we also know that $[q]M_+\subseteq [q]N_++M$, so we have that 
$$\# M/([q]N_++M)\leqslant\# M/[q]M_+.$$ From here we can conclude that 
\[\# N/[q]N_+\leqslant\# M/[q]M_++\# M/([q]N_++M)\leqslant 2Dq^d.\qedhere\]
\end{enumerate}
\pre
The same reduction steps hold when we allow torsion.
\lem
\label{reduction}
Suppose that for all finitely generated, semipositive, cancellative and integral binoids of dimension less or equal to $d$ the Hilbert-Kunz multiplicity exists. Then
\begin{enumerate}
 \item Let $N$ be a finitely generated, semipositive, cancellative and reduced  binoid of dimension $d$ and let $\{\mathfrak{p}_1,\dots,\mathfrak{p}_s\}$ be all minimal prime ideals of dimension $d$. Then the Hilbert-Kunz multiplicity of $N$ exists, and $$e_{HK}(N)=\sum_{i=1}^s e_{HK}(N/\mathfrak{p}_i).$$
 \item Let $N$ be a finitely generated, cancellative and semipositive binoid of dimension $d$. Then $$\hkf(N,q)\leqslant Cq^d,$$ for all $q$, where $C$ is some constant.
\end{enumerate}
\leme
\pr
The proof is similar to the one of Lemma \ref{reduction of torsion free}.
\pre

\lem
\label{finite birational extantion}
Suppose that for all finitely generated, cancellative, semipositive, integral binoids of dimension less than $d$, the Hilbert-Kunz function is bounded by $Cq^{d-1}$ for some constant $C$. Let $N$ be a finitely generated, semipositive, cancellative, integral binoid of dimension $d$ and let $\mathfrak{n}$ be an $N_+$-primary ideal of $N$. Let $M$ be a finite $N$-binoid which is birational over $N$. Suppose that $e_{HK}(\mathfrak{n}+M,M)$ exists. Then $e_{HK}(\mathfrak{n},N)$ exists and it is equal to $e_{HK}(\mathfrak{n}+M,M)$.
\leme
\pr
By Proposition \ref{strongly exact sequence} applied to the $N$-sets $S=N,T=M$ and the ideal $J=[q]\mathfrak{n}$ we have the strongly exact sequence $$\infty\longrightarrow (N\cap([q]\mathfrak{n}+M))/[q]\mathfrak{n} \longrightarrow N/[q]\mathfrak{n}\longrightarrow M/([q]\mathfrak{n}+M)\longrightarrow  (M/N)/([q]\mathfrak{n}+M/N)\longrightarrow\infty.$$
If $a\in M_+$ and $m_1,\dots,m_r$ are $N$-generators of $M$ then $a=m_i+n$, where $n\in N$. Hence from the property that $\mathfrak{n}$ is an $N_+$-primary ideal of $N$, we can conclude that $ka=km_i+kn\in \mathfrak{n}+M$, for some $k\in \n$. This means that $\mathfrak{n}+M$ is an $M_+$-primary ideal. So, by Proposition \ref{general equation of exact seq}, 
we know that $$\# (N\cap([q]\mathfrak{n}+M))/[q]\mathfrak{n} + \# M/([q]\mathfrak{n}+M) = \# N/[q]\mathfrak{n} +\#(M/N)/([q]\mathfrak{n}+M/N)$$
and we have $$\# M/([q]\mathfrak{n}+M) \leqslant \# N/[q]\mathfrak{n} +\#(M/N)/([q]\mathfrak{n}+M/N).$$

By Proposition \ref{birational to m} we know that there exists $b\in N$ such that $b+M\subseteq N$, and it is easy to check that $I:=b+M\neq\{\infty\}$ is an ideal of $N$. We also know that $I$ annihilates $M/N$, so by Lemma \ref{N/ann-set}, $M/N$ is an $N/I$-set and $N/I$ is a finitely generated, semipositive, cancellative binoid and torsion-free up to nilpotence.
So by Lemma \ref{HKF of annihilator} we have that $$\# (M/N)/([q]\mathfrak{n}+M/N)=\hkf(\mathfrak{n},M/N,q)=\hkf^{N/I}((I\cup\mathfrak{n})/I,M/N,q)$$ and by Proposition \ref{dimension decrease} we know that $\dim N/I<d$, so from the assumptions and by Lemma \ref{reduction of torsion free} we get $$\hkf^{N/I}((I\cup\mathfrak{n})/I,M/N,q)\leqslant C q^{d-1}.$$
Hence we have
\begin{equation}
 \#  M/([q]\mathfrak{n}+M)-Cq^{d-1} \leqslant \# N/[q]\mathfrak{n}.
\end{equation}
We know that $I=b+M$ is isomorphic to $M$, because $N$ is an integral binoid. So, by Corollary \ref{equality of e.s. corollary}, we have that $$\# N/[q]\mathfrak{n}+\# (I\cap [q]\mathfrak{n})/(I+[q]\mathfrak{n})=\# I/(I+[q]\mathfrak{n})+\# (N/I)/([q]\mathfrak{n}+N/I)$$ and from here we have $$\# N/[q]\mathfrak{n}\leqslant\# I/(I+[q]\mathfrak{n})+\# (N/I)/([q]\mathfrak{n}+N/I).$$
But we know that $\# (N/I)/([q]\mathfrak{n}+N/I)=\hkf(\mathfrak{n},N/I,q)$ and $I$ annihilates $N/I$. Similarly we can show that $$\hkf^N(\mathfrak{n},N/I,q)=\hkf^{N/I}((I\cup\mathfrak{n})/I,N/I,q)\leqslant C'q^{d-1},$$ where $C'$ is some constant.
Since $I\cong M$, we have that $\# I/(I+[q]\mathfrak{n})=\# M/([q]\mathfrak{n}+M)$, so we can conclude that  $$\# N/[q]\mathfrak{n}\leqslant\# M/([q]\mathfrak{n}+M)+C'q^{d-1}.$$
From the last result and $(2.10)$ we have that $$e_{HK}(\mathfrak{n},N)=\lim_{q\to\infty}\dfrac{\# N/[q]\mathfrak{n}}{q^d}=\lim_{q\to\infty}\dfrac{\# M/([q]\mathfrak{n}+M)}{q^d}.$$
From Lemma \ref{extended ideal commutes with q} we have $\# M/([q]\mathfrak{n}+M)=\# M/[q](\mathfrak{n}+M)$, which means that $$e_{HK}(\mathfrak{n},N)=e_{HK}(\mathfrak{n}+M,M).$$
\pre

\section[Hilbert-Kunz multiplicity of not normal binoids]{Hilbert-Kunz multiplicity of integral, cancellative, torsion-free binoids}

In this Section we will prove existence and rationality of the Hilbert-Kunz multiplicity of integral, cancellative, torsion-free binoids. This result is for fixed positive characteristic essentially due to K. Eto, see \cite[Theorem 2.2]{Eto}.
\tr
\label{not normal}
Let $N$ be a finitely generated, positive, cancellative, torsion-free and integral binoid. Let $\mathfrak{n}$ be an $N_+$-primary ideal of $N$. Then $e_{HK}(\mathfrak{n},N)$ exists and this is a rational number.
\tre
\pr
We work with the normalization $\hat{N}$ of $N$, which is finite over $N$ and birational over $N$. We also can show that $\mathfrak{n}+\hat{N}$ is an $\hat{N}_+$-primary ideal. So by Theorem \ref{HKF for normal binoid} we know that $e_{HK}(\mathfrak{n}+\hat{N},\hat{N})$ exists and it is a rational number. 
So, by Lemma \ref{finite birational extantion}, we have 
\[e_{HK}(\mathfrak{n},N)=e_{HK}(\mathfrak{n}+\hat{N},\hat{N}).\qedhere \]
\pre
We also can prove this Theorem using the convex geometrical approach to commutative, finitely generated, positive, cancellative, torsion-free, integral binoid. We give a second proof for this in the case $\mathfrak{n}=N_+$ and we use some easy Lemmas from Section 3.3 below.

\pr
We will use induction on $d=\dim N$ to prove this theorem.

Let $d=1$. Then our binoid becomes a numerical semigroup so, by Proposition \ref{numerical semigroup} $\# N/[q]N_+=nq+1$, for $q\gg 0$, where $n$ is the smallest positive element of $N$. Hence $e_{HK}(N)$ exists and it is equal to $n\in\n$.

Now we will assume that our theorem is true for every binoid fulfilling the assumptions and having smaller dimension than $\dim N=d$. Let $\{v_1,\dots,v_k\}\subseteq N$ be a set of generators for $N$. We know that $\hat{N}^\bullet=\r_+N^\bullet\cap\z^d$ and
$C:=\r_+N^\bullet=H^+_{\sigma_1}\cap\cdots\cap H^+_{\sigma_s}$, where $\sigma_i$ are linear integer forms. By Proposition \ref{m} there exists $n\in N^\bullet$ such that $n+\hat{N}\subseteq N$ and $I:=n+\hat{N}$ is a finitely generated ideal of $N$ and of $\hat{N}$. 
We can write $$N=(N\cap I)\uplus (N\setminus I)$$ and from here we can write $$N/[q]N_+=(N\cap I)/([q]N_+\cap I)\cupdot \big((N\setminus I)\setminus([q]N_+\cap (N\setminus I))\big).$$
We know that $N\cap I=I$ and from Lemma, \ref{inner cone of f.g. set} we have $N\setminus I=\bigcup_{j\in J} T_j$, where $J=J_1\uplus\cdots\uplus J_s$ is a finite index set, $T_j=N\cap H_{i,\alpha}$ and $H_{i,\alpha}:=\{y\mid \sigma_i(y)=\alpha\}$ is a hyperplane parallel to the facet hyperplane $H_{\sigma_i}$.
By Lemma \ref{plane and binoid intersection}, we know that $N_i:=H_{\sigma_i}\cap N$ (we will think that $\infty$ is added) is a finitely generated binoid with a smaller dimension and also, by Lemma \ref{planes and binoid intersection}, that the $T_j$'s, $j\in J_i,$ are finitely generated $N_i$-sets. Hence we can write
$$\begin{aligned}
 \# (N\setminus I)\setminus([q]N_+\cap (N\setminus I))&=\#\bigcup_{j\in J}T_j/([q]N_+\cap\bigcup_{j\in J}T_j)\\
 &\leqslant \sum_{j\in J}\# T_j/([q]N_+\cap\bigcup_{j\in J}T_j) \\
 &\leqslant\sum_{j\in J}\# T_j/([q]N_+\cap T_j)\\
 &=\sum_{j\in J}\# T_j/(([q]N_++N)\cap H_{i,\alpha})\\
 &=\sum_{j\in J}\# T_j/([q]N_+\cap H_{i,\alpha}+T_j)\\
 &=\sum_{j\in J}\# T_j/([q](N_i)_++T_j).
\end{aligned}$$

By the induction hypothesis $$e_{HK}(N_i)=\lim_{q\rightarrow \infty} \frac{\# N_i/[q](N_i)_+}{q^{\dim N_i}}$$ exists and so, by Theorem \ref{bound of HKF}, we have that
\begin{equation*}\# (N\setminus I)\setminus([q]N_+\cap (N\setminus I))\leqslant\sum_{j\in J}\# T_j/([q]N_{i+}+T_j)= \sum_{j\in J}\hkf(T_j,q)\leqslant\beta q^{d-1}
\end{equation*}
for some constant $\beta$. Similarly we can prove that
\begin{equation}\# (\hat{N}\setminus I)\setminus(([q]N_++\hat{N})\cap (\hat{N}\setminus I))\leqslant \gamma q^{d-1},\end{equation}
where $\gamma$ is some constant.
Hence, we can also write $$\hat{N}/([q]N_++\hat{N})=(\hat{N}\cap I)/(([q]N_++\hat{N})\cap I)\cupdot (\hat{N}\setminus I)\setminus(([q]N_++\hat{N})\cap (\hat{N}\setminus I))$$ and so by Theorem \ref{HKF for normal binoid} and by $(2.11)$ we can conclude that $$\lim_{q\rightarrow \infty}\frac{\#(\hat{N}\cap I)/(([q]N_++\hat{N})\cap I)}{q^d}=\lim_{q\rightarrow \infty}\frac{\#\hat{N}/([q]N_++\hat{N})}{q^d}$$ is a rational number.

Let us check that $$\bigcup^k_{i=1} (qv_i+I)\subseteq [q]N_+\cap I\subseteq ([q]N_++\hat{N})\cap I.$$
If $qv_i+a$ with $a\in I$, then $qv_i+a\in I$ and $qv_i+a\in [q]N_+$, so the first inclusion is clear. The second inclusion is clear from $[q]N_+\subseteq [q]N_++\hat{N}$.
Hence we have inequalities $$\# (N\cap I)/\bigcup^k_{i=1} (qv_i+I)\geqslant \#(N\cap I)/(([q]N_++\hat{N})\cap I).$$
But we know that $N\cap I=I\cong \nn$, so $$\# (N\cap I)/\bigcup^k_{i=1} (qv_i+I)=\#\nn/([q]N_++\nn)$$ and $$\#(N\cap I)/(([q]N_++\hat{N})\cap I)=\#(\nn\cap I)/(([q]N_++\hat{N})\cap I).$$ Hence  $$\lim_{q\rightarrow\infty}\frac{\#(N\cap I)/([q]N_+\cap I)}{q^d}$$ exists, because of the following inequalities:
$$\begin{aligned}
\lim_{q\rightarrow \infty}\frac{\#\hat{N}/([q]N_++\hat{N})}{q^d}&=\lim_{q\rightarrow\infty}\frac{\# (N\cap I)/\bigcup^k_{i=1} (qv_i+I)}{q^d} \\
&\geqslant \lim_{q\rightarrow\infty}\frac{\#(N\cap I)/([q]N_+\cap I)}{q^d}\\
&\geqslant \lim_{q\rightarrow\infty}\frac{\#(N\cap I)/(([q]N_++\hat{N})\cap I)}{q^d}\\
&=\lim_{q\rightarrow \infty}\frac{\#\hat{N}/([q]N_++\hat{N})}{q^d}.
\end{aligned}$$
So from this and $(2.9)$ we can conclude that 
\[e_{HK}(N)=\lim_{q\rightarrow\infty}\frac{\# N/[q]N_+}{q^d}=\lim_{q\rightarrow \infty}\frac{\#\hat{N}/([q]N_++\hat{N})}{q^d}\in \q.\qedhere\]
\pre

\ex
Let $N=(\n^3)^\infty /(ae_1+be_2=ce_3)$ be a binoid and $\gcd(b,c)=1$. By condition $\gcd(b,c)=1$ there exists $u,v\in \n$ such that $-ub+vc=\pm1$.
Let us assume $-ub+vc=1$. Then we can think of our binoid as generated by vectors $x=(0,1),y=(c,au),z=(b,av)$ in $\z^2$. So we have $ax+by=cz$.
Because of $(0,1)\in N$ we have $(c,0),(b,0)\in N$ and by condition $\gcd(b,c)=1$ we get $(1,0)\in N$, hence $\diff N=\z^2$.

Let $K=(0,1),B=(b,av),D=(c,au),F=(b,0),G=(c,0),L=(0,av),H=(0,au),I=(c,au+1)$ and $A=\overleftrightarrow{KI}\cap l,C=\overleftrightarrow{OD}\cap l$, where $l$ is the line $x=b$. A point $E$ is on the line $x=c$ such that  $\overleftrightarrow{BE}\parallel \overleftrightarrow{OD}$.
\begin{enumerate}
   \item If $c>\max(a,b)$ then $c\geqslant b+1$ and $u=\frac{vc-1}{b}\geqslant\frac{vb+v-1}{b}=v+\frac{v-1}{b}\geqslant v$.

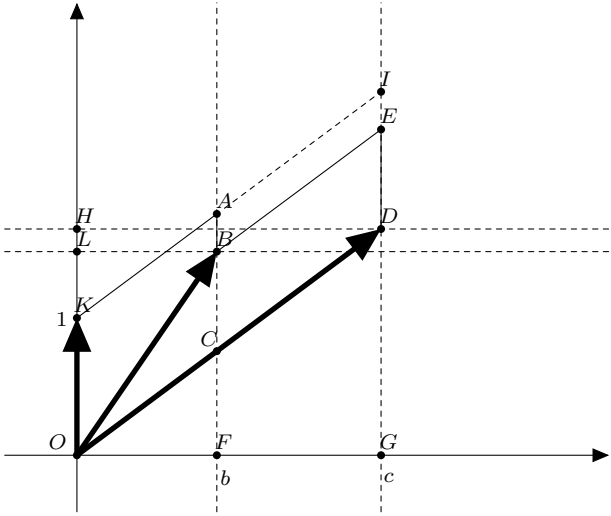
\begin{tikzpicture}[line cap=round,line join=round,>=triangle 45,x=1.0cm,y=1.0cm]
\draw[->,color=black] (-0.95,0) -- (7,0);
\draw[->,color=black] (0,-0.75) -- (0,6);
\draw [->,line width=2pt] (0,0) -- (4,3);
\draw [->,line width=2pt] (0,0) -- (0,1.82);
\draw [line width=0.4pt,dash pattern=on 2pt off 2pt] (4,-0.75) -- (4,6);
\draw [line width=0.4pt,dash pattern=on 2pt off 2pt,domain=-0.95:7] plot(\x,{(-3-0*\x)/-1});
\draw [line width=0.4pt,dash pattern=on 2pt off 2pt] (1.84,-0.75) -- (1.84,6);
\draw [line width=0.4pt,dash pattern=on 2pt off 2pt,domain=-0.95:7] plot(\x,{(-2.7-0*\x)/-1});
\draw [line width=0.4pt,dash pattern=on 2pt off 2pt,domain=1.94:4.06] plot(\x,{(--7.28--3*\x)/4});
\draw [->,line width=2pt] (0,0) -- (1.84,2.7);
\draw (0,1.82)-- (1.84,3.2);
\draw (1.84,3.2)-- (1.84,2.7);
\draw (1.84,2.7)-- (4,4.32);
\draw (4,4.32)-- (4,3);
\draw (0,0)-- (4,3);
\begin{scriptsize}
\fill [color=black] (0,0) circle (1.5pt);
\draw[color=black] (-0.25,0.18) node {$O$};
\fill [color=black] (4,3) circle (1.5pt);
\draw[color=black] (4.1,3.18) node {$D$};
\fill [color=black] (0,1.82) circle (1.5pt);
\draw[color=black] (0.09,1.99) node {$K$};
\draw[color=black] (-0.2,1.8) node {$1$};
\fill [color=black] (4,0) circle (1.5pt);
\draw[color=black] (4.1,0.18) node {$G$};
\draw[color=black] (4.1,-0.3) node {$c$};
\fill [color=black] (0,3) circle (1.5pt);
\draw[color=black] (0.1,3.18) node {$H$};
\fill [color=black] (1.84,0) circle (1.5pt);
\draw[color=black] (1.93,0.18) node {$F$};
\draw[color=black] (1.95,-0.3) node {$b$};
\fill [color=black] (0,2.7) circle (1.5pt);
\draw[color=black] (0.09,2.87) node {$L$};
\fill [color=black] (1.84,3.2) circle (1.5pt);
\draw[color=black] (1.94,3.38) node {$A$};
\fill [color=black] (1.84,2.7) circle (1.5pt);
\draw[color=black] (1.94,2.87) node {$B$};
\fill [color=black] (4,4.32) circle (1.5pt);
\draw[color=black] (4.1,4.5) node {$E$};
\fill [color=black] (4,4.82) circle (1.5pt);
\draw[color=black] (4.06,4.99) node {$I$};
\fill [color=black] (1.84,1.38) circle (1.5pt);
\draw[color=black] (1.74,1.55) node {$C$};
\end{scriptsize}
\end{tikzpicture}

By direct computation we can show that
$d(C,F)=\dfrac{abu}{c}$, so $d(B,C)=av-\dfrac{abu}{c}=a(v-\dfrac{bu}{c})=\dfrac{a}{c}<1$.
Which means, via Theorem \ref{HKF for normal binoid}, that the bounded region whose area we have to compute is $OKABEDO$. Hence $$e_{HK}(N)=c-(1-\dfrac{a}{c})(c-b)=a+b-\dfrac{ab}{c}.$$

   \item If $b<c\leqslant a$ then $d(B,C)=\dfrac{a}{c}\geqslant 1$. So $B$ is not in between of $A$ and $C$. Hence the bounded region is $OKID$ and $e_{HK}(N)=c$.
   \item Let $c\leqslant b$. Then also in this case the bounded region is $OKID$ and $e_{HK}(N)=c$.
  \end{enumerate}
\exe
\rem
We know the answer in this example also from the Proposition \ref{HKF of 3 generated binoid}.
\reme

\section{Cancellative, torsion-free binoids}

\tr
\label{not integral, reduced}
Let $N$ be a commutative, finitely generated, positive, cancellative, torsion-free and  reduced binoid and let $\{\mathfrak{p}_1,\dots,\mathfrak{p}_s\}$ be all minimal prime ideals with dimension $d=\dim N$.
Then $$e_{HK}(N)=\lim_{q\rightarrow \infty} \dfrac{\hkf(N,q)}{q^d}=\sum_{i=1}^s e_{HK}(N/\mathfrak{p}_i)$$ exists and this is a rational number.
\tre
\pr
This follows from Lemma \ref{reduction of torsion free} (1) and Theorem \ref{not normal}.
\pre

\df
Let $\trian$ be a simplicial complex on the vertex set $V$. The binoid associated to $\trian$ is given by $F(V)/I_{\trian}=:M_{\trian}$, where $I_{\trian}$ is the ideal $\{ f\in F(V)\mid  \supp(f)\nsubseteq \trian\}$ of the free binoid $F(V)\cong (\n^{|V|})^\infty$.
A binoid $M$ is called a \emph{simplicial binoid} if there exists a simplicial complex $\trian$ such that $M\cong M_{\trian}$. Then $\trian$ is unique and called the \emph{underlying simplicial complex} of $M$.
\dfe
\cor
Let $\trian$ be a simplicial complex on the vertex set $V$ and let $M_{\trian}$ be the corresponding simplicial binoid. Then $$e_{HK}(M_{\trian})=k,$$ where $k$ is the number of facets of maximal dimension.
\core
\pr
By \cite[Theorem 6.5.8]{Simone}, $M_{\trian}$ is a finitely generated, positive, cancellative, torsion-free, reduced binoid. The faces of $\trian$ correspond to the prime ideals of $M_{\trian}$ via $$F \longmapsto \langle V\setminus F\rangle=\mathfrak{p}_F$$ and the facets correspond to the minimal prime ideals.
Moreover, we have $$M_{\trian}/\mathfrak{p}_F=(\n^{|F|})^\infty$$ and $e_{HK}((\n^m)^\infty)=1$ by Example \ref{N^n}.
So if $F_1,\dots,F_k$ are the facets of maximal dimension then $\mathfrak{p}_i=\mathfrak{p}_{F_i},i=1,\dots,k$ 
are the minimal prime ideals of maximal dimension and hence by Theorem \ref{not integral, reduced} we have 
\[e_{HK}(M_{\trian})=\sum_{i=1}^k e_{HK}(M_{\trian}/\mathfrak{p_i})=k.\qedhere\]
\pre
\lem
\label{not integral and torsion free up to nilpotence}
Let $N$ be a commutative, finitely generated, positive, cancellative binoid which is torsion-free up to nilpotence. Then $\hkf(N,q)\leqslant Cq^d$,where $C$ is some constant and $d=\dim N$.
\leme
\pr
This follows from Theorem \ref{not normal} and Lemma \ref{reduction of torsion free} (2).
\pre
\cor
\label{bound for N/I}
Let $N$ be a finitely generated, integral, positive, cancellative and torsion-free binoid and let $I\neq \infty$ be an ideal of $N$. Then $\hkf(N/I,q)\leqslant Cq^d$,where $C$ is some constant and $d=\dim N/I$.
\core
\pr
We know by \cite[Lemma 2.1.20]{Simone}, that $N/I$ is torsion-free up to nilpotence and finitely generated, positive, cancellative. So by Lemma \ref{not integral and torsion free up to nilpotence} we have the result.
\pre
\lem
\label{ideal case}
Let $N$ be a commutative, finitely generated, positive, cancellative, integral and torsion-free binoid and let $I\neq \infty$ be an ideal of $N$. Then $e_{HK}(I)$ exists and is equal to $e_{HK}(N)$.
\leme
\pr

By Corollary \ref{equality of e.s. corollary} we have
\begin{equation*}
\# N/[q]N_+ +\# I\cap [q]N_+/(I+[q]N_+)=\# I/(I+[q]N_+)+\# (N/I)/[q](N/I)_+.
\end{equation*}

Also by Corollary \ref{bound for N/I} we know that $\hkf(N/I,q)\leqslant Dq^{d'}$, where $d'=\dim N/I<d$. Hence from $(2.11)$ we get
\begin{equation}
\# I/(I+[q]N_+)\geqslant \hkf(N,q)-\# (N/I)/[q](N/I)_+\geqslant \hkf(N,q)-Dq^{d'}.
\end{equation}

Let $\infty\neq f\in I$. Then we have an exact sequence $\infty\rightarrow (f)\rightarrow I\rightarrow I/(f)\rightarrow\infty$. So we have $$\# I/(I+[q]N_+)\leqslant \# (f)/((f)+[q]N_+)+\# (I/(f))/[q](I/(f))_+.$$ By Lemma \ref{HKF of annihilator}, the $N$-set  $I/(f)$ is an $N/(f)$-set.
By Corollary \ref{bound for N/I} we know that \[\hkf(N/(f),q)\leqslant Dq^{d''},\] where $d''=\dim N/(f)<d$ and by Theorem \ref{bound of HKF} we have $\hkf(I/(f),q)\leqslant \alpha q^{d''}$. Hence we can conclude that
\begin{equation}
 \# I/(I+[q]N_+)\leqslant \# (f)/((f)+[q]N_+)+\alpha q^{d''}.
\end{equation}

But we know that $(f)=f+N$ is isomorphic to $N$ as an $N$-set so $$\# (f)/((f)+[q]N_+)=\hkf(N,q).$$
Now by $(2.12)$ and $(2.13)$ we have $$\begin{aligned}e_{HK}(N)&=\lim_{q\rightarrow\infty} \frac{\hkf(N,q)-Dq^{d'}}{q^d}\\
&\leqslant \lim_{q\rightarrow\infty} \frac{\# I/(I+[q]N_+)}{q^d}=e_{HK}(I)\\
&\leqslant\lim_{q\rightarrow\infty} \frac{\hkf(N,q)+\alpha q^{d''}}{q^d}=e_{HK}(N)\end{aligned}$$
\pre

\section{Integral and cancellative binoids with torsion}

Let $N$ be a finitely generated, semipositive, cancellative and integral binoid. We know that $$N\subseteq \diff N\cong(\z^m\times T)^\infty=(\z^m)^\infty\wedge T^\infty,$$ where smashing is over the trivial binoid $\t=\{0,\infty\}$.
Here  $T$ is the torsion part of the difference group, which is a finite commutative group, hence $T=\z/k_1\times\cdots \times \z/k_l$. We will write elements $x\in N^\bullet$ as $x=f\wedge t$ with $f\in\z^m$ and $t\in T$. This representation is unique.

The relation $\sim_{\tf}$ on $N$ given by $a\sim_{\tf} b$ if $na=nb$ for some $n\geqslant 1$ is a congruence and $N_{\tf}:=N/\sim_{\tf}$ is a torsion-free binoid which we call the \emph{torsion-freefication}.
If $F=\{f\in (\z^m)^\infty\mid  \exists t\in T, f\wedge t\in N\}$ then $N_{\tf}\cong F\subseteq \z^m$. Hence we may assume $N\subseteq F\wedge T^\infty$ and $f\wedge t_1\sim_{\tf} g\wedge t_2$ if and only if there exists $n\in\n$ such that $n(f\wedge t_1)=n (g\wedge t_2)$, which means $f=g\in F$. We define the subsets $F\wedge t:=\{g\wedge t \mid g\in F\}$, $F_t:=\{f\in F \mid f\wedge t\in M\}$ and $M_t:=F_t\wedge t$.

\prop
Let $F$ be a finitely generated, positive, cancellative, integral and torsion-free binoid. If $T$ is a finite group then
$e_{HK}(F\wedge T^\infty)=e_{HK}(F)\cdot |T|$.
\prope
\pr
This follows from Theorem \ref{not normal}, Example \ref{finite group's HK} and Theorem \ref{multiplicity of smash}.
\pre
In the following we write $N\subseteq F\wedge T^\infty$, where $F\cong N_{\tf}$, $T$ a finite group and $\diff N=\diff F\wedge T^\infty$. The strategy is to relate the Hilbert-Kunz multiplicity of $N$ with that of $F\wedge T^\infty$.
\lem
\label{correspondence of ideals}
Let $F$ be a binoid, $T^\infty$ be a group binoid. Then we have a bijection between ideals of $F$ and ideals of $F\wedge T^\infty$. The $F_+$-primary ideals correspond to $(F\wedge T^\infty)_+$-primary ideals.
\leme
\pr
We have an inclusion $$F \xrightarrow{\;i\;} F\wedge T^\infty,~f \longmapsto f\wedge 0.$$ So we can consider for an ideal in $F$ its extended ideal in $F\wedge T^\infty$, in other words we use the map $\mathfrak{a}\mapsto \mathfrak{a}+F\wedge T^\infty$. For ideal generators $f_j\wedge t_j,j\in J$, in $F\wedge T^\infty$ we have $$\langle f_j\wedge t_j\mid j\in J\rangle=\langle f_j\wedge 0\mid j\in J\rangle\subseteq F\wedge T^\infty,$$ because $0\wedge t_j$ are units.
Hence we have the map $\mathfrak{b}\mapsto i^{-1}(\mathfrak{b})$, where $\mathfrak{b}$ is an ideal of $F\wedge T^\infty$, which is inverse to the extension map.

Let $\mathfrak{p}\subseteq F$ be an $F_+$-primary ideal and $f\wedge t\in (F\wedge T^\infty)_+$. Then there exists $l\in\n$ such that $lf \in \mathfrak{p}$, so $l(f\wedge t)=lf\wedge lt=(lf\wedge 0)+(0\wedge lt) \in \mathfrak{p}+F\wedge T^\infty$, because $0\wedge lt$ is a unit.

Now let $\mathfrak{q}$ be an $F\wedge T^\infty_+$-primary ideal and $f\in F_+$. Then $f\wedge 0\in F\wedge T^\infty$ is not a unit, so there exists $m\in \n$ such that $m(f\wedge 0)=mf\wedge 0\in \mathfrak{q}$. Hence $mf \in i^{-1}(\mathfrak{q})$.
\pre
\lem
\label{finite birational ext}
Let $N$ be a finitely generated, semipositive, cancellative, integral binoid. If $F=N_{\tf}$ and $T$ is the torsion subgroup of $\diff N$ then $F\wedge T^\infty$ is finite and birational over $N$.
\leme
\pr
We can assume that $\diff N=(\z^d)^\infty\wedge T^\infty$, where $d=\dim N$, $T=\{t_1,\dots,t_m\}$ is a finite abelian group. We know $F\wedge T^\infty \subseteq \diff N$, so it is clear that $F\wedge T^\infty$ is birational over $N$. If $f\wedge t\in (F\wedge T^\infty)^\bullet$ then there exists $t'\in T$ such that $f\wedge t'\in N$ and
\begin{equation}
m(f\wedge t)=mf\wedge mt=mf\wedge 0=mf\wedge mt'=m(f\wedge t')\in N,
\end{equation}
 which means $F\wedge T^\infty$ is pure integral over $N$. Let $f_i\wedge t_i,1\leqslant i \leqslant k$ be the generators of $N$, then $\{f_i\wedge t_j,\mid 1\leqslant i \leqslant k,1\leqslant j \leqslant m\}$ will give us an $N$-generating system of $F\wedge T^\infty$.
So it means that $F\wedge T^\infty$ is finite over $N$.
\pre
Note that the notion birational makes sense though the corresponding binoid algebras are not integral in general.
\lem
\label{F smash T}
Let $F$ be a finitely generated, positive, cancellative, torsion-free, integral binoid, $T^\infty$ be a torsion group binoid and $\mathfrak{p}$ be an $F_+$-primary ideal of $F$. Then $$\hkf^{F\wedge T^\infty}(\mathfrak{p}+F\wedge T^\infty,F\wedge T^\infty,q)=|T|\cdot \hkf^F(\mathfrak{p},F,q).$$
In particular, the Hilbert-Kunz multiplicity exists and $$e_{HK}^{F\wedge T^\infty}(\mathfrak{p}+F\wedge T^\infty,F\wedge T^\infty)=|T|\cdot e_{HK}^F(\mathfrak{p},F).$$
\leme
\pr
For every $q\in\n$ we have by Lemma \ref{smash of quotients} $$F\wedge T^\infty/([q]\mathfrak{p}\wedge T^\infty\cup F\wedge (\infty))=F\wedge T^\infty/[q]\mathfrak{p}\wedge T^\infty\cong F/[q]\mathfrak{p}\wedge T^\infty,$$
and $[q]\mathfrak{p}\wedge T^\infty=[q]\mathfrak{p}+F\wedge T^\infty$. Hence $$F\wedge T^\infty/([q]\mathfrak{p}+F\wedge T^\infty)\cong F/[q]\mathfrak{p}\wedge T^\infty,$$ which means $$\hkf^{F\wedge T^\infty}(\mathfrak{p}+F\wedge T^\infty,F\wedge T^\infty,q)=|T|\cdot \hkf^F(\mathfrak{p},F,q).$$
The existence of the Hilbert-Kunz multiplicity follows from Theorem \ref{not normal}.
\pre

Let $N$ be a binoid with $N\subseteq F\wedge T^\infty$, where $F\cong N_{\tf}$, $T$ a finite group. Then we have the following diagram.
$$\begin{tikzpicture}[node distance=2cm, auto]
  \node (F) {F};
  \node (FT) [below of=F]  {$F\wedge T^\infty$};
  \node (N) [left of=FT]   {$N$};
  \draw[->] (N) to node {} (FT);
  \draw[->] (F) to node [swap] {$i$} (FT);
\end{tikzpicture}
$$
where $i(f)=f\wedge 0$.

\tr
\label{HKM with torsion binoid}
Let $N$ be a finitely generated, semipositive, cancellative, integral binoid and $\mathfrak{n}$ be an $N_+$-primary ideal of $N$. Then $$e_{HK}^N(\mathfrak{n},N)=|T|\cdot e_{HK}^F(\mathfrak{m},F),$$ where $F=N_{\tf}$, $T$ is the torsion subgroup of $\diff N$ and $\mathfrak{m}=i^{-1}(\mathfrak{n}+F\wedge T^\infty)$ is an ideal of $F$.
\tre
\pr
 We know that $\mathfrak{n}+F\wedge T^\infty$ is a primary ideal and by Lemma \ref{correspondence of ideals} that $\mathfrak{m}$ is an $F_+$-primary ideal of $F$. So by Lemma \ref{F smash T} we know that $e_{HK}(\mathfrak{m}+F\wedge T^\infty,F\wedge T^\infty)$ exists and is equal to $|T|\cdot e_{HK}^F(\mathfrak{m},F)$.
 By Lemma \ref{correspondence of ideals} we know that $\mathfrak{m}+F\wedge T^\infty=\mathfrak{n}+F\wedge T^\infty$ so $$e_{HK}^{F\wedge T^\infty}(\mathfrak{n}+F\wedge T^\infty,F\wedge T^\infty)=e_{HK}^{F\wedge T^\infty}(\mathfrak{m}+F\wedge T^\infty,F\wedge T^\infty).$$
 Because of Lemma \ref{finite birational ext} and using induction over the dimension we can apply Lemma \ref{finite birational extantion}. Hence $e_{HK}^N(\mathfrak{n},N)$ exists and is equal to $$e_{HK}^{F\wedge T^\infty}(\mathfrak{n}+F\wedge T^\infty,F\wedge T^\infty)=|T|\cdot e_{HK}^F(\mathfrak{m},F).$$
\pre

\tr
\label{f.g,s.p,c,r binoid}
Let $N$ be a finitely generated, semipositive, cancellative, reduced binoid and $\mathfrak{n}$ be an $N_+$-primary ideal of $N$. Then $e_{HK}(\mathfrak{n},N)$ exists and is rational number.
\tre
\pr
This follows directly from Lemma \ref{reduction} and Theorem \ref{HKM with torsion binoid}.
\pre
\ex
For $a\in \n_+$ we have a binoid $\langle x,y\rangle/ax=ay$. Then this binoid can be realized as $\langle (1\wedge 0),(1\wedge 1)\rangle\subseteq \n\wedge (\z/a)^\infty$. In this case, the torsion-freefication is $\n$ and the Hilbert-Kunz multiplicity is $a$ by Theorem \ref{HKM with torsion binoid}, since $e_{HK}((\n)^\infty)=1$ by Example \ref{N^n} and $HKF(T^\infty,q)=|T|$ by Example \ref{finite group's HK}.
\exe
\ex
Let $a=(2,1),b=(3,0)\in (\n \times \z/2)^\infty$ be the generators of a binoid $N$. From Example \ref{double N} we know that this binoid is not torsion-free. It is not difficult to see that  $N_{\tf}\cong\n^\infty\setminus \{1\}$ and the torsion group is $T=\z/2$. So by Theorem \ref{HKM with torsion binoid}, we have $$e_{HK}(N)=|T|\cdot e_{HK}(N_{\tf})=2\cdot 2=4.$$
\exe

\lem
\label{3dim ex}
Let $N=\langle X,Y,Z\rangle/aX+bY=cZ$ be a binoid, $a,b,c\in\n_+$ and $d=\gcd(a,b,c)$.
\begin{enumerate}
 \item The map $\phi:N \rightarrow \big(\z^2\times \z/c\big)^\infty$ given by $$X \longmapsto (c,0,0),~ Y \longmapsto (0,c,0),~Z \longmapsto (a,b,1),~\infty\longmapsto\infty,$$ is a well defined injective binoid homomorphism and $N$ is a cancellative binoid.
 \item The difference group of $N$ is isomorphic to $\z^2\times \z/d$.
 \item For $d=1$ the homomorphism $\psi:N \longrightarrow (\z^2)^\infty$, given by $$X \longmapsto (c,0),~ Y \longmapsto (0,c),~ Z \longmapsto (a,b),\infty\longmapsto \infty$$
 is injective and our binoid is torsion-free.
 \item The torsion-freefication of $N$ is $$N_{\tf}\cong N'=\langle X',Y',Z'\rangle/\frac{a}{d}X'+\frac{b}{d}Y'=\frac{c}{d}Z',$$ where $\delta:N\rightarrow N',~X\mapsto X',~Y\mapsto Y',~Z\mapsto Z'$.
\end{enumerate}
 
\leme
\pr
\begin{enumerate}
 \item It is clear that this map is a well defined binoid homomorphism. Let $$\phi(eX+fY+gZ)=\phi(e'X+f'Y+g'Z).$$ This means that $$(ec+ga,fc+gb,g)=(e'c+g'a,f'c+g'b,g')\in\z^2\times\z/c.$$ So we have 
 \begin{align*}
  ec+ga &=e'c+g'a,\\
  fc+gb &=f'c+g'b,\\
  g &\equiv g'\mod c 
 \end{align*}
and we can write after some trivial manipulations
  $g =g'-ck,~ e =e'+ak,f =f'+bk,$ where $k\in \z$.
This means that $eX+fY+gZ=e'X+f'Y+g'Z$ in $N$, hence our map is injective. The cancellativity of $\z^2\times \z/c$ implies that $N$ is cancellative.

\item Let us consider the map $\varphi:\Gamma N\rightarrow \big(\z^2\times \z/c\big)^\infty $ with $\varphi|_N=\phi$. Let us denote $x=(c,0,0)$, $y=(0,c,0)$, $z=(a,b,1)$ in $\z^2\times \z/c$. If we write $a=da',b=db',c=dc'$ then $\gcd(a',b',c')=1$. So we have $(0,0,c')=c'z-a'x-b'y$ which means $(0,0,c')\in \im \varphi$. Hence $\Gamma N$ contains a torsion element of order $d$.

Now let  $ux+vy+wz$ be a torsion element. Then $ux+vy+wz=(uc+wa,vc+wb,w)$ and by torsion we have $uc+wa=0$, $vc+wb=0$. From here we have $w=-\frac{uc}{a}=-\frac{vc}{b}$ and $\frac{u}{a'}=\frac{v}{b'}$. If $\frac{u}{a'}=\frac{v}{b'}\notin \z$ then there exists a prime power $p^i$ such that $p^i\mid a'$, $p^i\mid b'$ and $p^i\nmid u$,$p^i\nmid v$. But $w=-\frac{uc'}{a'}\in\z$ implies that $p\mid c'$ and this contradicts $\gcd(a',b',c')=1$. So $\frac{u}{a'}=\frac{v}{b'}=m\in \z$ which means that $w=mc'$. Hence the torsion group of $\Gamma N$ is $\z^d$. Therefore $\Gamma N\cong \z^2\times \z/d$.

\item It is clear that $\psi$ is a well defined binoid homomorphism. Let $$\psi(eX+fY+gZ)=\psi(e'X+f'Y+g'Z).$$ From here we get $(ec+ga,fc+gb)=(e'c+g'a,f'c+g'b)$ and $$\frac{e-e'}{a}=\frac{f-f'}{b}=\frac{g'-g}{c}=\frac{p}{q}\in\q,$$ where $p\in \z$, $q\in\n_+,~\gcd(p,q)=1$. From $e-e'=a\frac{p}{q}\in \z$ we deduce that $q\mid a$ and similarly $q\mid b$ and $q\mid c$. Hence from $\gcd(a,b,c)=1$ we get $q=1$ and $e-e'=ap,~f-f'=bp,~g'-g=cp$. So this means that $eX+fY+gZ=e'X+f'Y+g'Z$ in $N$, hence our map is injective. From here we also get that our binoid is torsion-free.

\item We have the following commutative diagram.
$$\begin{tikzpicture}[node distance=2cm, auto]
  \node (N) {$N$};
  \node (gost) [right of=N]{};
  \node (Z2C) [right of=gost]  {$\big(\z^2\times \z/c\big)^\infty$};
  \node (N1) [below of=N]   {$N'$};
   \node (N2) [right of=N1]   {$N'$};
  \node (Z2) [right of=N2]{$(\z^2)^\infty$};
  \draw[->] (N) to node {$\phi$} (Z2C);
  \draw[->>] (N) to node {$\delta$} (N1);
  \draw[->] (Z2C) to node  {$\pi$} (Z2);
  \draw[->] (N1) to node {$[d]$} (N2);
  \draw[->] (N2) to node {$\psi$} (Z2);
\end{tikzpicture}$$
Since $\phi$ and $\psi\circ [d]$ are injective we can understand which elements are identified under $\delta$ by looking at $\pi$.

Let $n,n'\in N$. So $\phi(n)=(m,t),\phi(n')=(m',t')\in\z^2\times \z/c$. They are identified under the projection $\pi$ if and only if $m=m'$.

Also $n\equiv n'$ in the torsion-freefication if and only if there exists $k\in\n_+$ such that $kn=kn'$ holds in $N$. If this holds, then $\phi(kn)=(km,kt)=\phi(kn')=(km',kt')$ and $km=km'$ implies $m=m'$ in $\z^2$. If $\delta(n)=\delta(n')$, then $m=m'$ and therefore $cn=cn'$ in $N$. Hence the map $\delta:N\twoheadrightarrow N'$ gives the torsion-freefication of $N$.
\end{enumerate}
\pre
\ex
Let $N=\langle X,Y,Z\rangle/4X+12Y=16Z$ be a binoid. From Lemma \ref{3dim ex} we have an injective binoid homomorphism $\phi:N\rightarrow \big(\z^2\times \z/16\z\big)^\infty$. So $\phi(N)$ has generators $(16,0,0),(0,16,0),(4,12,1)$. If we choose the new generators $u=(4,-4,1),~v=(0,16,0),~w=(0,0,4)$, then $\phi(X)=4u+v-w$, $\phi(Y)=v$, $\phi(Z)=u+v$. Hence the difference group of our biniod is isomorphic to $\z^2\times \z/4$. The torsion-freefication of $N$ is $F=\langle X',Y',Z'\rangle/X'+3Y'=4Z'$ with $e_{HK}(F)=\frac{13}{4}$ by Theorem \ref{not normal} ($F$ is not normal), so $e_{HK}(N)=13$ by Theorem \ref{HKM with torsion binoid}.
\exe

We include some further observations about the relation between $N$ and $F\wedge T^\infty$, which give in part another proof of our result.
\lem
\label{element about inclusion}
Let $N\subseteq F\wedge T^\infty$ be a finitely generated, semipositive, cancellative and integral binoid. Then there exists $f\in F^\bullet$ such that $f\wedge 0+F\wedge T^\infty\subseteq N$.
\leme
\pr
Let $T=\z/{k_1}\times \cdots \times \z/{k_l}$. Then we can write $(0\wedge e_j)=(f_j\wedge t_j)-(f_j\wedge t'_j)$ where $e_j=(0,\dots,0,1,0,\dots,0)\in T$ and $(f_j\wedge t_j),(f_j\wedge t'_j)\in N$. We set $f\wedge 0:=\sum_{j=1}^l k(f_j\wedge t'_j)$, where $k=\sum_{j=1}^l k_j$. Then we have
$$\begin{aligned}
f\wedge(i_1,\dots,i_l)&=f\wedge 0+\sum_{j=1}^l i_j (0\wedge e_j)\\
&=\sum_{j=1}^l k(f_j\wedge t'_j)+\sum_{j=1}^l i_j (f_j\wedge t_j-f_j\wedge t'_j)\\
&=\sum_{j=1}^l i_j (f_j\wedge t_j)+\sum_{j=1}^l (k-i_j)(f_j\wedge t'_j)\in N,
  \end{aligned}$$
  where $0\leqslant i_j< k_j$.
For every $g\in F$ there exists $t'\in T$ such that $(g,t')\in N$. So from the above we can write  $(f+g)\wedge t=f\wedge (t-t')+g\wedge t'$ where both summands belong to $N$ for every $t\in T$. Hence $(f+F)\wedge T^\infty=f\wedge 0+F\wedge T^\infty\subseteq N$.
\pre

\lem
\label{include Fx0}
Let $N\subseteq F\wedge T^\infty$ be a finitely generated, semipositive, cancellative, integral binoid and $M=N[f_1\wedge 0,\dots,f_k\wedge 0]\subseteq F\wedge T^\infty$, where $f_i\wedge t_i$ are generators of $N$. Then $F\wedge 0\subseteq M$.
\leme
\pr
Let $f\wedge 0\in F\wedge 0$. Then there exists $t\in T$ such that $(f,t)\in N$. So we have $$f\wedge t=\sum_{i=1}^k \alpha_i(f_i\wedge t_i)=\sum_{i=1}^k \alpha_i f_i\wedge\sum_{i=1}^k \alpha_i t_i,$$ which means $f=\sum_{i=1}^k \alpha_if_i$. Hence $$f\wedge 0=\sum_{i=1}^k \alpha_i(f_i\wedge 0)\in M.$$
\pre
\lem
\label{bound for M_t}
Let $M$ be a finitely generated, semipositive, cancellative, integral binoid with $F\wedge 0\subseteq M$. Let $f\in F^\bullet$ be such that $f\wedge 0+F\wedge T^\infty\subseteq M$. Then for every $t\in T$ and $q\gg 0$ we have the subset inclusions
$$[(f)/(f)+[q]F_+\setminus(f)\cap[q]F_+/(f)+[q]F_+]\wedge t\subseteq M_t/(M_t\cap[q]M_+)\subseteq F_t/([q]F_++F_t)\wedge t.$$
\leme
\pr
We know by Lemma \ref{element about inclusion} that $(f+F)\wedge t\subseteq M_t$. Now if $$h\in (f+F/f+[q]F_+)\setminus(f+F\cap[q]F_+)/f+[q]F_+$$ then $h=f+g\notin [q]F_+,g\in F,g\notin [q]F_+$, so  $(f+g)\wedge t\notin M_t\cap[q]M_+$, which means $(f+g)\wedge t\in M_t/(M_t\cap[q]M_+)$.

If $x\in [q]F_++F_t$ then $x=h+qf+g$, where $h\in F_t,f\in F_+,g\in F$, so $x\wedge t=(h+qf+g)\wedge t=h\wedge t+q(f\wedge 0)+g\wedge 0$. By Lemma \ref{include Fx0} we know $f\wedge 0\in M_+, g\wedge 0\in M$ and by definition $h\wedge t\in M_t$, hence $x\wedge t\in M_t\cap [q]M_+$. This means $$M_t/(M_t\cap[q]M_+)\subseteq F_t/([q]F_++F_t)\wedge t.$$
\pre

\lem
\label{F_t non-zero ideal}
Let $M$ be a finitely generated, semipositive, cancellative, integral binoid with $F\wedge 0\subseteq M$. Then $F_t$ is a non-zero ideal of $F$.
\leme
\pr
Let $g\in F_t,h\in F$ then $g\wedge t\in M$ and $h\wedge 0\in M$. Hence $g\wedge t+h\wedge 0=(g+h)\wedge t\in M_t$, which means $g+h\in F_t$. We know from Lemma \ref{element about inclusion} there is $f\in F_t,f\neq\infty,$ so not zero.
\pre
\lem
Let $M$ be a finitely generated, semipositive, cancellative, integral binoid with $F\wedge 0\subseteq M$. Then $$e_{HK}(M)=|T|\cdot e_{HK}(F),$$ where $T$ is the torsion subgroup of $(\diff M)^\times$ and $F=M_{\tf}$.
\leme
\pr
We know $M\subseteq F\wedge T^\infty$, where $F$ is finitely generated, positive, cancellative, integral, torsion-free binoid. Let $f$ be chosen as in Lemma \ref{element about inclusion}. So by Proposition \ref{strongly exact sequence} (with $N=F$, $S=I=f+F$, $J=[q]F_+$) we have a strongly exact sequence $$\infty\rightarrow (f+F)\cap [q]F_+/(f+[q]F_+) \rightarrow (f+F)/(f+[q]F_+) \rightarrow F/[q]F_+\rightarrow F/((f+F)\cup [q]F_+)\rightarrow\infty.$$
We know by Lemma \ref{ideal case} $$(f+F)/(f+[q]F_+)=e_{HK}(F)\cdot q^d+O(q^{d-1})$$ and $$F/[q]F_+=e_{HK}(F)\cdot q^d+O(q^{d-1}),$$
where $d=\dim F$. Also $$F/((f+F)\cup [q]F_+)=(F/f+F)/[q](F/f+F)_+$$ and by Corollary \ref{bound for N/I} we know $\hkf(F/(f+F),q)\leqslant Dq^{d'}$, where $D$ is some constant and $d'=\dim (F/f+F)<d$. Hence by Corollary \ref{equality of e.s. corollary} we have $$\# (f+F)\cap [q]F_+/(f+[q]F_+)+\# F/[q]F_+=\# (f+F)/(f+[q]F_+)+\# F/((f+F)\cup [q]F_+).$$
From here we can conclude that $$\# (f+F)\cap [q]F_+/(f+[q]F_+)\leqslant Cq^{d-1},$$ where $C$ is some constant.
From lemma \ref{bound for M_t} we know $$[(f+F/f+[q]F_+)\setminus(f+F\cap[q]F_+)/(f+[q]F_+)]\wedge t\subseteq M_t/(M_t\cap[q]M_+)\subseteq F_t/([q]F_++F_t)\wedge t,$$ where $t\in T,M_t=F_t\wedge t=M\cap (F\wedge t), F\wedge t=\{(f\wedge t)\mid f\in F\}$. We know by Lemma \ref{F_t non-zero ideal},  $f+F$ and $F_t$ are the non zero ideals of $F$. So we have
$$\hkf(f+F,q)-\# (f+F)\cap [q]F_+/(f+[q]F_+)\leqslant \# M_t/(M_t\cap[q]M_+)\leqslant \hkf(F_t,q),$$ for every $t\in T$. We can write $M/[q]M_+=\bigcupdot_{t\in T} M_t/(M_t\cap[q]M_+)$, hence $$\hkf(M,q)=\# M/[q]M_+=\sum_{t\in T} \# M_t/(M_t\cap[q]M_+).$$ Therefore $$|T|\cdot[\hkf(f+F,q)-\# (f+F)\cap [q]F_+/(f+[q]F_+)]\leqslant \hkf(M,q)\leqslant \sum_{t\in T} \hkf(F_t,q).$$ Dividing by $q^d$ and using Lemma \ref{ideal case} we get $e_{HK}(M)=|T|\cdot e_{HK}(F).$
\pre

\section{Hilbert-Kunz function of binoid rings}

We want to relate the Hilbert-Kunz function of a binoid $N$ to the Hilbert-Kunz function of its binoid algebra $K[N]$. However, $K[N]$ is not a local ring. If $N$ is positive, then $K[N]$ contains the unique combinatorial maximal ideal $K[N_+]$ and we work with the localization $K[N]_{K[N_+]}$. In this setting we get via Proposition \ref{finite set and dimension} and Proposition \ref{K algebra of quotient}, immediately $$HKF^N(\mathfrak{n},S,q)=\# \big(S/(S+[q]\mathfrak{n})\big)=\dim_K K[S]/(K[S]\cdot K[\mathfrak{n}]^{[q]})=HKF^{K[N]}(K[\mathfrak{n}],K[S],q).$$
So the rationality results of the previous sections transform directly to results on the Hilbert-Kunz multiplicity of the localization of a binoid algebra.

This translation is more involved for semipositive binoids. We have seen in Lemma \ref{semipositive binoid algebra} that $K[N_+]$ is the intersection of finitely many maximal ideals $\mathfrak{m}_1,\dots,\mathfrak{m}_r$. Then $T=K[N]\setminus\mathfrak{m}_1\cap\cdots\cap K[N]\setminus\mathfrak{m}_r$ is a multiplicatively closed subset of $K[N]$ and the localization $K[N]_T$ is a semilocal ring with maximal ideals $\mathfrak{m}_i \cdot K[N]_T$ and $K[N_+]\cdot K[N]_T=\bigcap_{i=1}^r \mathfrak{m}_i\cdot K[N]_T$.

Now, for a semilocal Noetherian ring $R$ with Jacobson ideal $\mathfrak{m}=\bigcap_{i=1}^r \mathfrak{m}_i$ containing a field of positive characteristic, we can define the Hilbert-Kunz function as before.
For a finite $R$-module $M$ and an $\mathfrak{m}$-primary ideal $\mathfrak{n}$ we set
$$\hkf^R(\mathfrak{n},M,q)=\len(M/\mathfrak{n}^{[q]}M).$$
If $J$ is an ideal in a Noetherian ring $R$ with $V(J)=\{\mathfrak{m}_1,\dots,\mathfrak{m}_r\}$, then $$\len^R(M/JM)=\len^{R_T}(M_T/JM_T)$$ for $T=\bigcap_{i=1}^r R\setminus \mathfrak{m}_i$. In this way we consider $K[N_+]$-primary ideals in $K[N]$ for a semipositive binoid $N$, and we  write $\hkf^{K[N]}(K[\mathfrak{n}],K[S],q)$ instead of $\hkf^{K[N]_T}(K[\mathfrak{n}]\cdot K[N]_T,K[S]_T,q)$.
\tr
\label{HK binoid ring}
Let $K$ be a field of characteristic $p$, $N$ a finitely generated, semipositive binoid, $S$ an $N$-set, $\mathfrak{n}$ an $N_+$-primary ideal and $q=p^e$. Then we have $$\hkf^N(\mathfrak{n},S,q)=\hkf^{K[N]}(K[\mathfrak{n}],K[S],q).$$
\tre
\pr
We know, from Proposition \ref{K algebra of quotient}, that $$K[S/(S+[q]\mathfrak{n})]\cong K[S]/(K[S]\cdot K[\mathfrak{n}]^{[q]}),$$ and by Proposition \ref{finite set and dimension}, we can conclude that $$\# (S/(S+[q]\mathfrak{n}))=\dim_K K[S/(S+[q]\mathfrak{n})]=\dim_K K[S]/(K[S]\cdot K[\mathfrak{n}]^{[q]}).$$
\pre
\tr
\label{HKM binoid ring and binoid}
Let $K$ be a field of characteristic $p$, $N$ a finitely generated, semipositive binoid. Suppose that $\dim N=\dim K[N]$ and $e_{HK}(\mathfrak{n},N)$ exists. Then
 $$e_{HK}^{K[N]}(K[\mathfrak{n}],K[N])=e_{HK}(\mathfrak{n},N)$$ and it is independent of the (positive) characteristic of $K$.
\tre
\pr
If we take $S=N$ in Theorem \ref{HK binoid ring} then we have $$\hkf^{K[N]}(K[\mathfrak{n}],K[N],q)=\hkf^N(\mathfrak{n},N,q).$$ By assumption $(q=p^e,d=\dim N)$ $$e_{HK}(\mathfrak{n},N)=\lim_{q\rightarrow \infty}\dfrac{\hkf^N(\mathfrak{n},N,q)}{q^d}=\lim_{e\rightarrow \infty}\dfrac{\hkf^{K[N]}(K[\mathfrak{n}],K[N],q)}{q^d}=e_{HK}^{K[N]}(K[\mathfrak{n}],K[N]).$$
\pre
\tr[Miller conjecture for cancellative binoid rings]
Let $K$ be a field of characteristic $p$, $N$ a binoid.
Suppose that $\dim N=\dim K[N]$ and $e_{HK}(\mathfrak{n},N)$ exists. Then
 $$\lim_{p\rightarrow\infty} \frac{\hkf^{(\z/p)[N]}((\z/p)[\mathfrak{n}],(\z/p)[S],p)}{p^{\dim N}}$$ exists and equals $$\lim_{p\rightarrow\infty} e_{HK}^{(\z/p)[N]}((\z/p)[\mathfrak{n}],(\z/p)[S]).$$
\tre
\pr
This follows from the identity $\hkf^{(\z/p)[N]}((\z/p)[\mathfrak{n}],(\z/p)[S],p)=\hkf^N(\mathfrak{n},S,p)$ from Theorem \ref{HK binoid ring} and the existence of the limit.
\pre
\tr
Let $K$ be a field of characteristic $p$, $N$ be a finitely generated, semipositive, cancellative, reduced binoid and $\mathfrak{n}$ be an $N_+$-primary ideal of $N$. Then $$e_{HK}^{K[N]}(K[\mathfrak{n}],K[N])$$ exists, is independent of the characteristic of $K$ and it is a rational number.
\tre
\pr
By Theorem \ref{f.g,s.p,c,r binoid} we know that $e_{HK}(\mathfrak{n},N)$ exists and is a rational number. So by Theorem \ref{HKM binoid ring and binoid} we have the result.
\pre

\rem
Note that $K[F\wedge T^\infty]$ over an algebraically closed field $K$ where $\chara K$ does not divide $|T|$ has an easy structure, it is the product of $|T|$ copies of $K[F]$ and $\spec K[F\wedge T^\infty]$ is the disjoint union of $|T|$ copies of $\spec K[F]$.
\reme

\chapter{Partition of \protect\qnn}

In this chapter we will study the partitions of the $N$-sets \qnn. K.-I.Watanabe showed in his paper (\cite{Watanabe}), that for toric rings there exists finitely many non-isomorphic irreducible parts in the partitions of \qnn. We will prove the same result for finitely generated, positive, integral, cancellative and torsion-free binoids.

\section{Generalities}

Let $S$ be an $N$-set and $\phi:N\rightarrow N$ a binoid homomorphism. Then $S$ with the new operation $$n\dotplus s:=\phi(n)+s$$ is again an $N$-set, denoted by $^\phi S$.

For every commutative binoid $N$ and every $q\in \n$ we have a binoid homomorphism $[q]$, where
$$N\xrightarrow{\; [q]\; } N,~n\longmapsto qn.$$
So we have a set \qnn whose elements are the same as $N$'s and $N$ acts via the homomorphism $[q]$ as above. From Lemma \ref{homomorphism gives N set} we can think of \qnn as a pointed $N$-set. We also set $qN:=\{qn\mid n\in N\}$, $b+qN:=\{b+qn\mid n\in N\}$,  for $b\in N$ and $[q]N_+$ is the extended ideal.
\lem
\label{f.g.N-set}
Let $N$ be a finitely generated binoid. If $q$ is a positive integer then \qnn is a finitely generated $N$-set.
\leme
\pr
If $n_1,\dots,n_k$ are the generators of $N$ then the set $\{a_j n_i\mid 0\leqslant  a_j<q, 1\leqslant  i\leqslant  k\}$ will generate \qnn, i.e. $\qn=\bigcup_{i,j}(a_j n_i+qN)$.
\pre
\cor
\label{HKF by smash}
Let $N$ be a commutative binoid, then $^qN\wedge_N N/N_+ \cong N/{[q]N_+}$.
\core
\pr
Take $I=N_+$ and $S=\qn$ in Proposition \ref{quotient}.
\pre

This gives us another way to compute the Hilbert-Kunz function.
We know that \qnn is an $N$-set and $\dotplus$ is an operation on \qnn. So we have the equivalence relation $\sim_N$ on \qnn.
\prop
\label{HK function is the number of N generators}
Let $N$ be a finitely generated, semipositive binoid and $k=k(q)$ be the number of minimal $N$-set generators of  \qnn. Then $\hkf(N,q)=k\cdot|N^\times|$.
\prope
\pr
Let $(S,p$) be a finitely generated $N$-set and $\{s_1,\dots,s_k\}$ be a minimal generating set of $S$, 
which is determined up to multiplication by units. Then by Proposition \ref{quotient} we have 
\[S\wedge_N N/N_+\cong S/(N_++S).\] But we know $S/(N_++S)=\{n+s_i\mid n\in N^\times,1\leqslant i\leqslant k\}\cup\{p\}$, 
because $n_1+s_i\neq n_2+s_j$.
Hence 
\[\# (S\wedge_N N/N_+)=\# \big(S/(N_++S)\big)=|\{n+s_i\mid n\in N^\times,1\leqslant i\leqslant k\}|=k\cdot|N^\times|.\]
So if $S=\qn$ then 
\[\hkf(N,q)=\# N/[q]N_+=\#(\qn\wedge_N N/N_+)=k\cdot |N^\times|. \qedhere\]
\pre

Let $N$ be a commutative binoid and $q\geqslant1$ be some integer. For every $b\in N$ we define the sets
\begin{align*}
\qp(b)&:=\{a\in N\mid  \exists n\in N^\bullet,a+qn\in b+qN\}\subseteq N,\\
\qn(b)&:=\{a\in N\mid  a+b\in qN\}.
\end{align*}
If $a\in$ $\qp(b)$ then there exist $n_1\in N^\bullet,n_2\in N$ such that $a+qn_1=b+qn_2$. For $n\in N$ we have $n\dotplus a=qn+a$, so $qn+a+qn_1=b+qn+qn_2$ which means $n\dotplus a\in$ $\qp(b)$.
So $\qp(b)$ is an $N$-subset of \qnn. Similarly we can show that $\qn(b)$ is an $N$-set.
\lem
\label{inclusion}
Let $N$ be an integral binoid. If $q$ is some positive integer then $\qp(b)=$ $\qp(b+qc)$ and $\qn(b)\subseteq$ $\qn(b+qc)$ for every $b\in N,c\in N^\bullet$.
\leme
\pr
If $a\in$ $\qp(b)$ then there exists $n\in N^\bullet$ such that $a+qn=b+qn'$ so $$a+q(n+c)=b+qn'+qc=b+q(n'+c)$$ for every $c\in N^\bullet$. By integrality, we have $n+c\in N^\bullet$. Hence $a\in$ $\qp(b+qc)$. If $a\in$ $\qp(b+qc)$ then there exists $n\in N^\bullet$ such that $a+qn=b+qc+qn'=b+q(c+n')$ and we know that $c+n'\in N^\bullet$ so $a\in$ $\qp(b)$.

Let $d\in \qn(b)$, then $d+b=qm$ for some $m\in N$. So $d+(b+qc)=qm+qc=q(m+c)\in qN$, which means $\qn(b)\subseteq$ $\qn(b+qc)$.
\pre
\lem
\label{irred qn}
Let $N$ be an integral binoid. If $q$ is some positive integer and $b\in N^\bullet$ then $\qp(b)$ is an irreducible component of \qnn.
\leme
\pr
 We know that $b\in$ $\qp(b)$ and $b\in N^\bullet$, so $\qp(b)\neq \{\infty\}$. If we take two  elements $x,y\in$ $\qp(b)^\bullet$ then there exist $n_1,n_2,n_3,n_4\in N^\bullet$ such that
 $$x+qn_1=b+qn_2~\operatorname{and}~y+qn_3=b+qn_4.$$ Hence from $$x+qn_1+qn_4=b+qn_2+qn_4=y+qn_3+qn_2$$ and $x,y\in\qp(b)^\bullet$, we have $x\sim_N y$.
 Also, if $x\in$ \qnn and $n\dotplus x\in$ $\qp(b),n\dotplus x\neq \infty$, then there exist $n_3\in N^\bullet, n_4\in N$ such that $(n\dotplus x)+qn_3=b+qn_4$, so $x+q(n+n_3)=b+qn_4$ and this means $x\in$ $\qp(b)$. Hence $\qp(b)$ is one of the equivalence classes of \qnn for the  relation $\sim_N$. So by Proposition \ref{equivalence class}, it follows that $\qp(b)$ is an irreducible component of \qnn.
\pre
\cor
Let $N$ be an integral binoid. If $q$ is a positive integer then $$\qn=\bigcupdot_{b\in T}~\qp(b)$$ for some subset $T\subseteq N$.
\core
\pr
By Lemma \ref{irred qn} we know that $\qp(b)$ is an irreducible component of \qnn and for every element $c\in N^\bullet$ we know that $c\in$ $\qp(c)$. So \qnn is the union of the sets $\qp(b)$. Now it is enough to take $T$ to be a set of representatives for the different $\qp(b)$'s.
\pre

\cor
\label{finite union}
Let $N$ be a finitely generated and integral binoid. If $q$ is a positive integer then there exist $b_1,\dots,b_s\in N$ such that $\qn=\bigcupdot_{i=1}^s~\qp(b_i)$.
\core
\pr
Let $e_1,\dots,e_k$ be the generators of $N$. Take any element $a\in$ \qnn. Then we can write $a=\sum_{i=1}^k a_ie_i+qc$ where $0\leqslant a_i\leqslant q-1,c\in N^\bullet$. Hence $a\in$ $\qp(\sum_{i=1}^k a_ie_i)$ and it is clear that $\{\sum_{i=1}^k a_ie_i\mid 0\leqslant a_i\leqslant q-1\}$ is finite.
\pre
\lem
\label{make easier}
Let $N$ be a finitely generated and integral binoid. If $q$ is a positive integer and $b\in N$ then there exists $b'\in N$ such that $\qp(b)=$ $\qn(b')$.
\leme
\pr
We know by Lemma \ref{f.g.N-set} and Proposition \ref{Subset of f.g.} that $\qp(b)$ is a finitely generated $N$-set. So if $a_1,\dots,a_s$ are the generators of $\qp(b)$ then there exist $n_1,\dots,n_s\in N^\bullet$ such that $a_i+qn_i\in b+qN$ for $i\in \{1,\dots,s\}$. Hence $a_i\in$ $\qn((q-1)b+qn_i)$ for $i\in \{1,\dots,s\}$. If we take $n=n_1+\cdots+n_s$ then by Lemma \ref{inclusion} we have that $a_i\in$ $\qn((q-1)b+qn)$ for every $i\in \{1,\dots,s\}$, which means $\qp(b)\subseteq$ $\qn((q-1)b+qn)$. The other inclusion is obvious, so $\qp(b)=$ $\qn((q-1)b+qn)$.
\pre
\cor
\label{finite qN partition}
Let $N$ be a finitely generated and integral binoid. If $q$ is a positive integer then there exists $b_1,\dots,b_s\in N$ such that $\qn=\bigcupdot_{i=1}^s~\qn(b_i)$.
\core
\pr
This follows from Corollary \ref{finite union} and Lemma \ref{make easier}.
\pre

\prop
\label{gen_n(S)}
Let $N$ be a finitely generated, positive and integral binoid.
If $S$ is a finitely generated $N$-set then $S\wedge_N N/N_+=\{a\wedge 0\mid a\in \gen_N(S)\}\cup \{\infty\wedge \infty\}$.
In particular this is true for the $N$-sets $\qn(b)$.
\prope
\pr
If $a\in S^\bullet$ is not a generator then there exist $a'\in \gen_N(S)$ and $0\neq n\in N^\bullet$ such that $a=a'+n$. $N$ is a positive binoid, so $N/N_+=\{0,\infty\}$. This means we have $$a\wedge 0=(a'+n)\wedge 0=a'\wedge n=a'\wedge \infty=\infty\wedge \infty.$$  Also $$S\wedge_N N/N_+\subseteq \{a\wedge 0\mid a\in \gen_N(S)\}\cup \{\infty\wedge \infty\}.$$ Let $a,a'\in \gen_N (S)$. By Proposition \ref{minimal generating set of N-set} we can not write $a=b+n,a'=b'+n'$, where $b,b'\in S, 0\neq n,n'\in N$. So $a\wedge 0\neq a'\wedge 0$. The other inclusion is clear.
\pre

\section{Relation to normalization}

\df
let $N$ be a finitely generated, integral and cancellative binoid and $(\diff N)^\times=\z^d\times T$, where $d=\rk(N)$. An element of the normalization $\nn$ is a \emph{gap} if it does not belong to $N$. A gap $b$ is called \emph{primary} if we can not write $b=a+m,m\in N_+$, where $a$ is a gap. 
\dfe
The normalization $\nn$ is generated as a binoid over $N$ by 0 and the primary gaps.

Let $A = \{n_1,n_2,\dots, n_k\}$ be a nonempty set of positive integers. The set of all integers of the form $x_1 n_1 + x_2 n_2 + \dots + x_r n_r$ is the subset of $\n$ generated by $A$ and is denoted by $\langle A \rangle$. 
\df
Let $a_1,\dots,a_k$ be a set of positive integers such that $\gcd(a_1,\dots,a_k)=1$. The \emph{numerical semigroup} generated by $a_1,\dots,a_k$ is the set $\langle a_1,\dots,a_k\rangle$.
\dfe
The following statement follows also from Theorem \ref{not normal}, since the normalization of $N$ is $\n$ and the extended ideal of $N_+$ is $\n_{\geqslant n_1}$. The smallest generator of a numerical semigroup is called its \emph{multiplicity}.
\prop
\label{numerical semigroup}
Let $N$ be a numerical semigroup. If $q$ is a positive integer greater than any gap of $N$ then $\# N/[q]N_+=qn_1$, where $n_1$ is the smallest generator of $N$.
\prope
\pr
Let $$n_1<n_2<\cdots<n_k\in \n_+$$ be the generators of $N$, then we have $$N=\Big\{\sum_{i=1}^k x_in_i\mid x_1,\dots,x_k\in \n\Big\}.$$
From $\gcd(n_1,\dots,n_k)=1$ there exist $d_1,\dots,d_k\in\n$ such that $d_1n_1+\cdots+d_kn_k=1$. Let $D=\max |d_i|$. Let $s=Da_1\sum_{i=1}^k a_i$. Then for any $0\leqslant r < a_1$, $$s+r=\sum_{i=1}^k (Da_i+rd_i)a_i\in N.$$ So any integer $x\geqslant s$ can be written as $x=ma_1+s+r$, where $0\leqslant r < a_1$ and $m\in\n$. Hence $x\in N$ which means that we have finitely many gaps, so call them $b_1,\dots,b_s$. 

Let $q>b_i$ for all $i$, which means in particular $q\in N$. Hence it is easy to see that
\begin{enumerate}
\item if $b\in N$ then $\qp(b)=b+q\n$
\item if $b\notin N$ then $\qp(b+q)=b+q+q\n$
\end{enumerate}
and $$\qn=\Big(\bigcupdot_{q>b\in N}\, \qp(b)\Big)\cupdot\Big(\bigcupdot_{q>b\notin N}\, \qp(b+q)\Big).$$
Which means they are all isomorphic to the normalization $\hat{N}\cong\n$. But we know that $\n$ is generated over $N$ as an $N$-set by the elements $0,1,\dots, n_1-1$. Hence by Proposition \ref{gen_n(S)} and Corollary \ref{HKF by smash} we have $\# N/[q]N_+=\#\gen_N(\qn)=qn_1$.
\pre

\lem
\label{irreducibility}
Let $N$ be a finitely generated, integral and cancellative binoid. Let $L$ be an irreducible $\nn$-set (or $\diff N$-set). Then $L$ is also an irreducible $N$-set.
\leme
\pr
That $L$ is an irreducible $\nn$-set means $x\sim_{\nn} y$ for every $x,y\in L$ and we have to show that $x\sim_N y$. It is enough to see that if we have $n,m\in \nn^\bullet$ such that $n+x=m+y$ then $x\sim_N y$. By definition $n,m\in \nn^\bullet$ means $n=a-b$ and $m=c-d$, where $a,b,c,d\in N^\bullet$. So $n+x=m+y$ implies that $(a+d)+x=(b+c)+y$. Hence by integrality $a+d,b+c\in N^\bullet$ and  $x\sim_N y$.
\pre
\lem
\label{irr to irr}
Let $N$ be a finitely generated, integral and cancellative binoid and $q\in\n$. Let $L\subseteq$ $^q\nn$ be irreducible as $\nn$-set (via the $[q]$ map). Then $L\cap N\subseteq$\qnn is irreducible as an $N$-set.
\leme
\pr
Let $x,y\in L\cap N$, then $n\dotplus x=m\dotplus y$, where $n,m\in \nn$. This means that $qn+x=qm+y$ and we can write $n-m=n'-m'$, where $n',m'\in N$. So $$qm'+qm+y=q(n+m')+x=q(n'+m)+x=qn'+qm+x$$ and by cancellation we get $qm'+y=qn'+x$, which means $x\sim_N y$. Hence, by Proposition \ref{irr set has one component}, $L\cap N$ is an irreducible $N$-set.
\pre
\lem
\label{irr comp to irr comp}
Let $N$ be a finitely generated, integral and cancellative binoid. If $L$ is an irreducible component of $^q\nn$ then $L\cap$ \qnn is an irreducible component of \qnn.
\leme
\pr
By Proposition \ref{irr to irr} $L\cap$\qnn is an irreducible $N$-set. It is easy to see that if $L\cap$\qnn is not an irreducible component of \qnn then $L$ is not an irreducible component of $^q\nn$.
\pre
\prop
\label{easy normalization}
Let $N\subseteq \n^k$ be a binoid with $\hat{N}=\n^k$ and let $q$ be a (large enough) positive integer. Then the irreducible components of \qnn are isomorphic to $\n^k_c\cap N$, where $\n^k_c=c+q\n^k$ is an irreducible component of the $\n^k$-set $^q\n^k$.
\prope
\pr
Let $I^k_q=\big\{(i_1,\dots,i_k)\mid i_j\in \{0,1,\dots,q-1\}\big\}$.
It is easy to see that $$^q\n^k=\bigcupdot_{c\in I^k_q} c+q\n^k.$$ We denote $\n^k_c=c+q\n^k$. It is an equivalence class of $^q\n^k$ and isomorphic to $\n^k$, which means that it is an irreducible component of  $^q\n^k$. So by Lemma \ref{irr comp to irr comp}, $\n^k_c\cap N$ is an irreducible component of \qnn.
\pre

Let $N$ be a binoid that satisfies all conditions in Proposition \ref{easy normalization} and let  $c=(c_1,\dots,c_k)$ be the generator of $\n^k_c$. We consider the following diagram.
$$\begin{tikzpicture}[node distance=2cm, auto]
  \node (N) {$N\times \n^k$};
  \node (NK) [right of=N]  {$\n^k$};
  \node (N1) [below of=N]   {$N\times\n^k_c$};
  \node (NK1) [right of=N1]{$\n^k_c$};
  \draw[->] (N) to node {$+$} (NK);
  \draw[->] (N.220) to node [swap] {$Id$} (N1.138);
  \draw[->] (N.310) to node  {$\psi$} (N1.53);
  \draw[->] (NK) to node [swap] {$\psi$} (NK1);
  \draw[->] (N1) to node {$+$} (NK1);
\end{tikzpicture}
$$
Here if $(m_1,\dots,m_k)\in N$ and $(n_1,\dots,n_k)\in \n^k$ then $\psi(n_1,\dots,n_k)=(c_1+qn_1,\dots,c_k+qn_k)$,
$((m_1,\dots,m_k),(n_1,\dots,n_k))\mapsto (m_1+n_1,\dots,m_k+n_k)$ and $\psi(m_1+n_1,\dots,m_k+n_k)=(c_1+q(m_1+n_1),\dots,c_k+q(m_k+n_k))$.\\
On the other hand $((m_1,\dots,m_k),\psi(n_1,\dots,n_k))=((m_1,\dots,m_k),(c_1+qn_1,\dots,c_k+qn_k))\mapsto (c_1+qn_1+qm_1,\dots,c_k+qn_k+qm_k)$. So this diagram commutes. Hence we can assume the irreducible components $\n^k_c$ to be $\n^k$ with the usual addition action of $N$.

\prop
\label{finite gaps}
Let $N$ be a finitely generated binoid of the form $(\n^k\setminus \{a_1,\dots,a_n\})^\infty$. If $q\gg 0$ is an integer then $$\# N/[q]N_+=(p+1)(q^k-n)+dn,$$ where $p$ is the number of primary gaps and $d$ is the number of $N$-generators of $\n^k\setminus \{0\}$. Moreover, we have the bound $$k+p\leqslant d\leqslant k(p+1).$$
\prope
\pr
It is easy to check that the normalization is $\nn=\n^k$, so the elements $a_i$ are the gaps. So by Lemma \ref{irr comp to irr comp} the irreducible components of \qnn are $(c+q\n^k)\cap N$, where $c=(c_1,\dots,c_k)\in \n^k$ and $c_i<q$ for all $i$.
So for large $q$ we have $$(c+q\n^k)\cap N\cong \begin{cases} \n^k, & \mbox{if } c\mbox{ is not a gap,} \\ \n^k\setminus \{0\}, & \mbox{if } c\mbox{ is a gap,} \end{cases}$$ and these irreducible sets have in the first case one and in the second case $k$ generators, namely the standard vectors $e_1,\dots,e_k$.
\begin{enumerate}
\item If $c$ is not a gap then $(c+q\n^k)\cap N\cong\n^k$, so 0 and all primary gaps will generate $\n^k$. Hence $$\# \gen_N(\n^k)=p+1.$$
\item If  $c$ is a gap then $(c+q\n^k)\cap N=\n^k\setminus \{0\}$, so we have $\#\gen_N(\n^k\setminus \{0\})=d$.
\end{enumerate}
Hence $\# \gen_N(\qn)=(p+1)(q^k-n)+dn$.

Now we will show that $k+p\leqslant d\leqslant k(p+1)$. Let $e_1,\dots,e_k$ be the standard basis of $\n^k$. We know that the primary gaps are $N$-generators of $\n^k$. If $e_i,2e_i,\dots,k_ie_i$ are primary gaps and $(k_i+1)e_i$ is not, then $e_i,2e_i,\dots,k_ie_i,(k_i+1)e_i$ belong to the $N$-generators of $\n^k\setminus \{0\}$, where $k_i\in\n$. 
Hence $d\geqslant k+p$. It is also easy to see that the set $\{e_i,e_i+g_1,\dots,e_i+g_p\mid 1\leqslant i \leqslant k\}$ will generate $\n^k\setminus \{0\}$, where $g_1,\dots,g_p$ are the primary gaps.
\pre

By Proposition \ref{finite gaps} we can deduce the following results. Let $k=1$. Then $$k+p\leqslant d\leqslant k(p+1)$$ implies that $d=p+1$ and $\# N/[q]N_+=(p+1)(q^k-n)+dn=(p+1)(q-n)+(p+1)n=(p+1)n$, which means we recover Proposition \ref{numerical semigroup}, since the primary gaps are the gaps below the multiplicity.
\ex
\begin{enumerate}
\item Let $N=\n^2\setminus \{(1,0),(2,0),(0,1),(1,1),(2,1),(0,2),(1,2),(2,2)\}$. Then it is easy to see that $d=10=k+p$.
\item Let $N=\n^2\setminus \{(1,0),(2,0),(0,1),(1,1),(0,2)\}$ then $2+5< d=9< 2(4+1)=10$.
\end{enumerate}
\exe
\ex
Let $n\in \n$ and $N=\langle (k,0),(1,1),(0,k),(0,k+1),\dots,(0,2k-1)\rangle$ be a binoid. If $q\gg 0$ is an integer then the number of generators of $\qp(0)$ is equal to $2k-2$, so the number of generators can be large.
\exe

\section{The toric case}

Let $N$ be a finitely generated, positive, integral, cancellative and torsion-free binoid. Then by \cite[Proposition 2.17]{BrunsGubeladze} we can assume that $N\subseteq \n^d$ and $\diff N\cong\z^d$, where $d=\rk(N)$. Let $I$ be an ideal of $N$. By \cite[Proposition 2.1.14]{Simone} we know that it is finitely generated, which means that $I=(b_1+N)\cup\cdots\cup(b_s+N)$.
Let $C=\r_+N^\bullet\subset \r^d$ be the rational cone generated by $N^\bullet$. This cone has a representation $$C=H^+_{\sigma_1}\cap\cdots\cap H^+_{\sigma_s}$$ as an irredundant intersection of halfspaces defined by linear forms $\sigma_i$ on $\r^d$.
Each of the hyperplanes $H_{\sigma_i}$ is generated as a vector space by integral vectors. Therefore we can assume that $\sigma_i$ is a linear form on $\z^d$ with nonnegative values on $C$.

We define the sets $$L_q(x,I)=\{m\in \diff N\mid  x+qm\in I\}$$ and
$$\qp(x)=\{a\in N\mid \exists n\in N^\bullet, a+qn\in x+qN\}$$ where $x\in \z^d,~q\in \n$. It is easy to see that $\qp(x)$ is an $N$-subset of \qnn.
If $s\in L_q(x,I)$ and $n\in N$, then $n+s\in \diff N$. So $x+q(n+s)=(x+qs)+qn\in I$ which means that $n+s\in L_q(x,I)$. Hence $L_q(x,I)$ is also an $N$-set.
\ex
Let $N=\langle (2,0),(3,2),(3,3),(2,3),(0,2)\rangle\subseteq\n^2$ then $\hat{N}=(\n^2)^\infty$ and if $a,b\geqslant 2$ then $(a,b)\in N$. Let $1\leqslant q\in\n$ be a positive integer and if it is even then $\qp(0)\cong \hat{N}=(\n^2)^\infty$ else  $\qp(0)\cong\langle(2,0),(0,2),(1,1),(2,1),(1,2)\rangle$.
\exe
\prop
\label{isomorph}
Let $N$ be a finitely generated, positive, integral, cancellative and torsion-free binoid. Let's assume $N\subseteq \n^d$ and $\diff N=\z^d$. If $q$ is some positive number and $x\in \z^d$ then $L_q(x,N)\cong$ $\qp(x)$ as $N$-sets.
\prope
\pr
If we take any element $\infty\neq a\in$$\qp(x)$ then $x+qn=a+qm$ for some $n,m\in N^\bullet$, so $a=x+qm_a$, where $m_a=n-m=\frac{a-x}{q}\in \diff N=\z^d$.
We have an $N$-set homomorphism $\phi:$$\qp(x)\rightarrow L_q(x,N)$, where $\phi(a)=m_a$, since $a+nq\mapsto \frac{a+qn-x}{q}=n+\frac{a-x}{q}$. If
$a\neq b$ then $$x+qm_a\neq x+qm_b\Rightarrow \phi(a)\neq \phi(b),$$ which means that 
$\phi$ is injective. If $m\in L_q(x,N)$ then $x+qm\in N$ and also by
definition $x+qm\in$$\qp(x)$. We know $\phi(x+qm)=m$ so this means that $\phi$ is
surjective.
\pre
\rem
$L_q(x,I)\cong (x+q\diff N)\cap I$ which is an irreducible $N$-subset occurring in the $N$-set $^qI$.
\reme
If we prove that there are only finitely many non-isomorphic $N$-sets $L_q(x,N)$ then by Corollary \ref{finite union} we obtain also the same number of non-isomorphic irreducible $N$-subsets of \qnn for all $q$. To prove this it is enough to show that the coordinates of every minimal $N$-generator of $L_q(x,N)$ are bounded by some constant (not depending on $q$).

\lem
\label{plane and binoid intersection}
Let $N$ be a finitely generated, positive, integral, cancellative and torsion-free binoid. If $H_{\sigma_i}$ is a supporting hyperplane of the cone $\r_+ N^\bullet$ then $N_i:=(H_{\sigma_i}\cap N)^\infty$ is also a finitely generated, positive, integral, cancellative and torsion-free binoid.
\leme
\pr
We can assume $\diff N=\z^d$, where $d=\dim N$.
Let $\{v_1,\dots,v_k\}\subseteq N$ be the minimal set of generators of $N$. For $1\leqslant i\leqslant s$ $$(H_{\sigma_i}\cap N)^\infty=\Big\{\sum c_jv_j\in N\mid \sigma_i\Big(\sum c_jv_j\Big)=\sum c_j\sigma_i(v_j)=0,c_j\in \n\Big\}\cup\{\infty\}.$$
Hence if $\sigma_i(v_j)\neq 0$ then $c_j=0$, which means that $N_i$ is generated by $\{v_j\mid \sigma_i(v_j)=0\}$ as a binoid. Also positivity, integrality, cancellativity and torsion-freeness are true for these  subbinoids.
\pre

\lem
\label{planes and binoid intersection}
Let $N$ be a finitely generated, positive, integral, cancellative and torsion-free binoid. Let $H_{\sigma_1},\dots,H_{\sigma_s}$ be the supporting hyperplanes of $\r_+N^\bullet$ and $H_{\alpha,i}=\{y\in \z^d\mid \sigma_i(y)=\alpha\}$. Then $\big(N\cap(\bigcap_{i\in I}H_{\alpha_i,i})\big)^\infty$ is a finitely generated $(\bigcap_{i\in I}N_i)$-set, where $I\subseteq \{1,\dots,s\},$ and $N_i:=(H_{\sigma_i}\cap N)^\infty$.
\leme
\pr
Let $\{v_1,\dots,v_k\}\subseteq N$ be the minimal set of generators of $N$. We know that $\hat{N}=(\r_+N^\bullet\cap\z^d)^\infty$ and $\r_+N^\bullet=H^+_{\sigma_1}\cap\cdots\cap H^+_{\sigma_s}$.
Let $$J:=\{j\in\{1,\dots,s\}\mid \sigma_i(v_j)=0 \;\operatorname{for\; every}\; i\in I\}$$ and $$N_I:=\bigcap_{i\in I}N_i.$$ is a subbinoid of $N$ by Lemma \ref{plane and binoid intersection}. If $y\in N_I$ then $y=\sum_{j=1}^s c_jv_j,c_j\in\n$ and for every $i\in I$ we have $0=\sigma_i(y)=\sum_{j=1}^s c_j\sigma_i(v_j)$.
So if $\sigma_i(v_j)>0$ then it has to be $c_j=0$, which means that $y=\sum_{j\in J}c_jv_j$. Hence $N_I$ is a binoid generated by $\{v_j\mid j\in J\}$.
If $N\cap(\bigcap_{i\in I}H_{\alpha_i,i})=\varnothing$ then there is nothing to prove so we can assume $N\cap(\bigcap_{i\in I}H_{\alpha_i,i})\neq\varnothing$.
Let $y\in N\cap(\bigcap_{i\in I}H_{\alpha_i,i})$ and $n\in N_I$ then $$\sigma_i(n+y)=\sigma_i(n)+\sigma_i(y)=\sigma_i(y)=\alpha_i,$$ for every $i\in I$, which means that $n+y\in N\cap(\bigcap_{i\in I}H_{\alpha_i,i})$. Hence $\big(N\cap(\bigcap_{i\in I}H_{\alpha_i,i})\big)^\infty$ is an $N_I$-set.

Now for $z\in N\cap(\bigcap_{i\in I}H_{\alpha_i,i})$ we can write $$z=\sum_{j\in J}c_jv_j+\sum_{j\notin J}c_jv_j$$ and so $$\sigma_i(z)=\sum_{j\in J}c_j\sigma_i(v_j)+\sum_{j\notin J}c_j\sigma_i(v_j)=\sum_{j\notin J}c_j\sigma_i(v_j)=\alpha_i,$$ for every $i\in I$.
For every $j$ not belonging to $J$ there exists $l\in I$ such that $\sigma_l(v_j)>0$. Hence  $\sigma_l(z)=\sum_{j\notin J}c_j\sigma_i(v_j)=\alpha_l$ implies that $c_j$ is a bounded positive integer. Hence there are only finitely many $N_I$-generators of $\big(N\cap(\bigcap_{i\in I}H_{\alpha_i,i})\big)^\infty$.
\pre

\lem
\label{inner cone of f.g. set}
Let $N$ be a finitely generated, positive, integral, cancellative and torsion-free binoid. Let $H_{\sigma_1},\dots,H_{\sigma_s}$ be the supporting hyperplanes of $\r_+N^\bullet$ and $H_{\alpha,i}=\{y\in \z^d\mid \sigma_i(y)=\alpha\}$. Let $K:=\big(N\cap(\bigcap_{i\in I}H_{\alpha_i,i})\big)^\infty$ and $N_I:=\big(\bigcap_{i\in I}H_{\sigma_i}\cap N\big)^\infty$ for $I\subseteq \{1,\dots,s\}$.
Then there exist $n_1,\dots,n_t\in K$ such that $\bigcup_{i=1}^t (n_i+\hat{N_I})\subseteq K$ and $K\setminus \bigcup_{i=1}^t (n_i+\hat{N_I})$ is contained in finitely many hyperplanes parallel to the facets of $\hat{N_I}$.
\leme
\pr
Let $g_1,\dots,g_t$ be the $N_I$-generators of $K$. Then we have $K=\bigcup_{i=1}^t (g_i+N_I)$. Also by Proposition \ref{m} we have $m\in N_I$ such that $m+\hat{N_I}\subseteq N_I$, so $g_i+m+\hat{N_I}\subseteq g_i+N_I$ for every $i\in\{1,\dots,t\}$. Setting $n_i=g_i+m$ we get $\bigcup_{i=1}^t (n_i+\hat{N_I})\subseteq K$. From Proposition \ref{c}, we have $(g_i+N_i)\setminus (g_i+m +\hat{N_I})$ is contained in finitely many hyperplanes parallel to the facets of $\hat{N_I}$ for every $i\in\{1,\dots,t\}$. Hence $K\setminus \bigcup_{i=1}^t (n_i+\hat{N_I})\subseteq \bigcup_{i=1}^t (g_i+N_i)\setminus (n_i+\hat{N_I})$ is contained in finitely many hyperplanes parallel to the facets of $\hat{N_I}$.
\pre
\ex
Let $N=\langle (1,0),(2,1),(0,2),(1,3)\rangle$, then $\hat{N}=\n^2$. Let us take $$K=N\cap H_{x=1}=\langle (1,0),(1,2),(1,3)\rangle$$ and $$N_1=N\cap H_{x=0}=\langle(0,2)\rangle.$$
Then we have $K=\big((1,0)+N_1\big)\cup\big((1,3)+N_1\big)$.
\exe
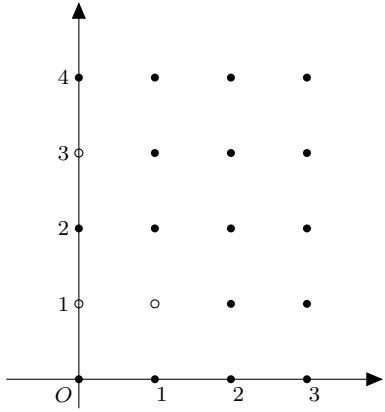
\begin{tikzpicture}[line cap=round,line join=round,>=triangle 45,x=1.0cm,y=1.0cm]
\draw[->,color=black] (-0.95,0) -- (4,0);
\draw[->,color=black] (0,-0.38) -- (0,5);
\begin{scriptsize}
\draw [color=black] (0,1) circle (1.5pt);
\draw[color=black] (-0.2,1) node {$1$};
\fill [color=black] (0,2) circle (1.5pt);
\draw[color=black] (-0.2,2) node {$2$};
\draw [color=black] (0,3) circle (1.5pt);
\draw[color=black] (-0.2,3) node {$3$};
\fill [color=black] (0,4) circle (1.5pt);
\draw[color=black] (-0.2,4) node {$4$};
\fill [color=black] (1,0) circle (1.5pt);
\draw[color=black] (1.09,-0.2) node {$1$};
\fill [color=black] (2,0) circle (1.5pt);
\draw[color=black] (2.1,-0.2) node {$2$};
\fill [color=black] (3,0) circle (1.5pt);
\draw[color=black] (3.1,-0.2) node {$3$};
\draw [color=black] (1,1) circle (1.5pt);
\fill [color=black] (1,2) circle (1.5pt);
\fill [color=black] (1,3) circle (1.5pt);
\fill [color=black] (1,4) circle (1.5pt);
\fill [color=black] (2,1) circle (1.5pt);
\fill [color=black] (2,2) circle (1.5pt);
\fill [color=black] (2,3) circle (1.5pt);
\fill [color=black] (3,1) circle (1.5pt);
\fill [color=black] (3,2) circle (1.5pt);
\fill [color=black] (3,3) circle (1.5pt);
\fill [color=black] (3,4) circle (1.5pt);
\fill [color=black] (2,4) circle (1.5pt);
\fill [color=black] (0,0) circle (1.5pt);
\draw[color=black] (-0.2,-0.2) node {$O$};
\end{scriptsize}
\end{tikzpicture}
\tr
\label{finitely many irr subset}
Let $N$ be a finitely generated, positive, integral, cancellative and torsion-free binoid. Let $q$ be some positive integer. Then the coordinates of a minimal $N$-generator of $L_q(x,N)$ are bounded by some constant depending neither on $q$ nor on $x=(x_1,\dots,x_d)$, ($0\leqslant x_i\leqslant q-1$ for all $i$).
\tre
\pr
We can assume that $N\subseteq \n^d$ and $\diff N=\z^d$, where $d=\dim N$.
Let $\{v_1,\dots,v_k\}\subseteq N$ be the minimal set of generators of $N$ and let $\sigma_i,~i=1,\dots,s$, be the describing linear forms of $\r_+N^\bullet$ with hyperplanes $H_i$. Let $a=(a_1,\dots,a_d)$ be one of the minimal $N$-generators of an irreducible $N$-subset $L_q(x,N)$, where $x=(x_1,\dots,x_d)$ and $0\leqslant x_i\leqslant q-1$ for all $i$. By Proposition \ref{m} there exists $n\in N$ such that $n+\hat{N}\subseteq N$. So $x+qa\in N=(n+\hat{N})\cup \big(N\setminus(n+\hat{N})\big)$. We consider two cases.

Let $x+qa\in n+\hat{N}$. Then $$x+qa\in \{n+\beta_1v_1+\cdots+\beta_kv_k\mid 0\leqslant \beta_j< q,\beta_j\in \n\},$$ because if $x+qa=n+(w+qv_i)$, where $w\in N$ then $x+q(a-v_i)\in N$, which contradicts the assumption that $a$ is a minimal generator.
So $x+qa=n+\beta_1v_1+\cdots+\beta_kv_k$ and $a=\frac{1}{q}(n-x+\beta_1v_1+\cdots+\beta_kv_k)$. The coordinates of $n$ are constant and the coordinates of $x$ and of $\beta_j$'s are bounded by $q$, so the coordinates of $a$ are bounded by a constant not depending on $q$.

Now, if $x+qa\in N\setminus(n+\hat{N})$ then by Proposition \ref{c} there exists $i_1\in\{1,\dots,s\}$ and $\alpha_{i_1}\in\n$ such that $$x+qa\in N\cap H_{\alpha_{i_1},i_1}=:K_1,$$ where $H_{\alpha_{i_1},i_1}=\{y\in \z^d\mid \sigma_{i_1}(y)=\alpha_{i_1}\}$. We set $M_1=N_{i_1}=N\cap H_{0,i_1}$.

By Lemma \ref{inner cone of f.g. set} there exists $n_{11},\dots,n_{1t_1}\in K_1$ such that  $\bigcup_{i=1}^{t_1} (n_{1i}+\hat{M_1})\subseteq K_1$ and $K_1\setminus \bigcup_{i=1}^{t_1} (n_{1i}+\hat{M_1})$ is contained in finitely many hyperplanes parallel to the facets of $\hat{M_1}$. 
If $x+qa\in \bigcup_{i=1}^{t_1} (n_{1i}+\hat{M_1})$ then we can write $$x+qa=n_{1r_1}+\sum_{j\in J_1} \beta_jv_j$$ for some $r_1\in\{1,\dots,t_1\}$ and $$a=\frac{1}{q}(n_{1r_1}-x+\sum_{j\in J_1} \beta_jv_j),$$ where $J_1=\{j\mid \sigma_{i_1}(v_j)=0\},0\leqslant\beta_j< q$.
So the coordinates of $a$ are bounded by a constant.
If $x+qa\in K_1\setminus \bigcup_{i=1}^{t_1} (n_{1i}+\hat{M_1})$ then there exists $i_2\neq i_1$ and $\alpha_{i_2}$ such that $$x+qa\in K_1\cap H_{\alpha_{i_2},i_2}=:K_2,$$ where $H_{\alpha_{i_2},i_2}=\{y\in \z^d\mid \sigma_{i_2}(y)=\alpha_{i_2}\}$.
Now by Lemma \ref{inner cone of f.g. set} there exists $n_{21},\dots,n_{2t_2}\in K_2$ such that  $\bigcup_{i=1}^{t_2} (n_{2i}+\hat{M_2})\subseteq K_2$ and $K_2\setminus \bigcup_{i=1}^{t_2} (n_{2i}+\hat{M_2})$ is contained in finitely many hyperplanes parallel to the facets of $\hat{M_2}$, where $M_2=N_{i_1}\cap N_{i_2}$. So we go on until either $x+qa\in n_{jr_j}+\hat{M_j}$ for some $j\in\n$. This situation occurs eventually since the dimension of $M_j$ drops and finally $K_j$ consists of only one point. So in each case $a$ is bounded by some constant depending neither on $q$ nor on $x=(x_1,\dots,x_d)$.
\pre
\tr
\label{finitely many isomorphic types}
Let $N$ be a finitely generated, positive, integral, cancellative and torsion-free binoid. Then there are only finitely many non-isomorphic irreducible $N$-subsets of \qnn for all $q$.
\tre
\pr
This follows from Corollary \ref{finite union} and Theorem \ref{finitely many irr subset}.
\pre
\prop
Let $n\in \n$ and $N=(\n^3)^{\infty}/(e_1+e_2=ne_3)$ be a binoid. If $q\gg 0$ is an integer then
$$\# N/[q]N_+=(2-\frac{1}{n})q^2-\frac{r(n-r)}{n}$$ where $q=ns+r,0\leqslant r<n,s,r\in\n$.
\prope
\pr
We can assume that our binoid is $N=\langle v_1,v_2,v_3,\infty\rangle\subseteq (\n^2)^\infty$, where $v_1=(n,0),v_2=(0,n),v_3=(1,1)\in \n^2$. So the elements of our binoid is the points $(ni+k,nj+k)\in \n^2$ where $i,j,k\in \n$. We will show that the irreducible parts (equivalence classes) of \qnn are isomorphic to the following sets (which are ideals):
$$S_t\cong (0,n-t)+N\cup(t,0)+N,$$ where $t=1,\dots,n$.

Fix $q$ and let $S$ be an irreducible part of \qnn. If $(x,y)\in N$ and $x\geqslant q,y\geqslant q$ then $(x,y)\sim_N(x-q,y-q)\in N$. Hence if $(x,y)\in S$ is a generator of $S$ then $x<q$ or $y<q$. Let $(x_1,y_1)$ be a generator of $S$ and $x_1<q$. Then also $y_1<qn$, since otherwise we would have $(x_1,y_1)\sim (x_1,y_1-qn)$. Let $(x_1',y_1')$ be another generator with $x_1'<q$ and $y_1'<qn$. Then $(x_1,y_1)\sim_N (x_1',y_1')$ implies that $x_1=x_1'$ and $y_1=y_1'$.
So there is no further generator where the $x$ coordinate is less than $q$. The same is true for a generator $(x_2,y_2)$ with $y_2<q$ and $x_2<qn$.
Hence $S$ is generated by two points $(x_1,y_1),(x_2,y_2)$ with the above property. 
Now we want to show that such an irreducible part $S$ is isomorphic to $S_t$ for some $t\in\{1,\dots,n\}$.

If these two points are different $(x_1<x_2,y_1>y_2)$ then we replace $S$ by $S-(x_1,y_2)$. Then these generators become $(0,y_1-y_2)$ and $(x_2-x_1,0)$.
So these points are equivalent, which means $(0,y_1-y_2)+qm_1=(x_2-x_1,0)+qm_2$ for some $m_1,m_2\in N$.
Hence with $m_1=i_1v_1+j_1v_2+k_1v_3$, $m_2=i_2v_1+j_2v_2+k_2v_3$ we get $$y_1-y_2=qn(j_2-j_1)+q(k_2-k_1),$$ $$x_2-x_1=qn(i_1-i_2)+q(k_1-k_2),$$ where $i_1,i_2,j_1,j_2,k_1,k_2\in \n$. From $y_1,x_2<qn$, $q\mid (y_1-y_2)$, $q\mid (x_2-x_1)$ and $qn\mid (y_1-y_2+x_2-x_1)$ we have $$y_1-y_2=q(n-t),$$ $$x_2-x_1=qt,$$ for some $0<t<n$. So $S$ is isomorphic to $S_t$.

If these two points are the same then $x_1=x_2<q,y_1=y_2<q$ and in this case it is easy to see that $S\cong S_n$.
So we can compute the number of parts isomorphic to $S_n$, which have only one generator. It is equal to all the non-isomorphic binoid points strictly inside the $q\times q$ box, which means $$q+2(q-n)+2(q-2n)+\dots+2(q-sn)=\frac{q^2+r(n-r)}{n}.$$ 
Hence we can conclude that the number of all parts which have 2 generators is equal to $q^2-\frac{q^2+r(n-r)}{n}$. So by Corollary \ref{finite qN partition} and Proposition \ref{gen_n(S)} we can conclude that
\[\# N/[q]N_+=2\bigg(q^2-\frac{q^2+r(n-r)}{n}\bigg)+\frac{q^2+r(n-r)}{n}=\bigg(2-\frac{1}{n}\bigg)q^2-\frac{r(n-r)}{n}.\qedhere\]
\pre

\end{document}